\documentclass[hidelinks,onefignum,onetabnum]{siamart220329}

\ifpdf
\hypersetup{
  pdftitle={Calibration-Based ALE Model Order Reduction for Hyperbolic Problems with Self-Similar Travelling Discontinuities},
  pdfauthor={Monica Nonino and Davide Torlo}
}
\fi

\newsiamremark{remark}{Remark}
\newsiamremark{hypothesis}{Hypothesis}
\crefname{hypothesis}{Hypothesis}{Hypotheses}
\newsiamthm{claim}{Claim}

\headers{Calibration ALE MOR for Hyperbolic Problems}{M. Nonino, D. Torlo}

\title{Calibration-Based ALE Model Order Reduction for Hyperbolic Problems with Self-Similar Travelling Discontinuities\thanks{Submitted to the editors on the 18th of March 2024.
		\funding{M. Nonino has been funded by the Austrian Science Fund (FWF) project 10.55776/F65, project 10.55776/P33477 and project 10.55776/ESP519. D.T. has been funded by a SISSA Mathematical Fellowship.}}}

\author{Monica Nonino\thanks{Faculty of Mathematics, University of Vienna, Vienna, AT
		(\email{monica.nonino@univie.ac.at}).}
	\and Davide Torlo\thanks{Università degli Studi di Roma, La Sapienza, Roma, IT
		(\email{davide.torlo@uniroma1.it}).}}

\usepackage{amsopn}

\usepackage{graphicx}%
\usepackage{multirow}%
\usepackage{amsmath,amssymb,amsfonts}%
\usepackage{lmodern}%
\usepackage{anyfontsize}
\usepackage{mathrsfs}%
\usepackage[title]{appendix}%
\usepackage{xcolor}%
\usepackage{textcomp}%
\usepackage{manyfoot}%
\usepackage{booktabs}%
\usepackage{algorithm}%
\usepackage{algorithmicx}%
\usepackage{algpseudocode}%
\usepackage{listings}%
\usepackage{bm}
\usepackage{mathabx}
\usepackage{tikz}
\usepackage{listofitems}
\usepackage{comment}

\usepackage{pgfplots}
\usepackage{soul}
\usepackage{adjustbox}
\usepackage{nomencl}
\makenomenclature

\renewcommand{\nompreamble}{The next list describes all the symbols that have been used in this manuscript.}
\usetikzlibrary{positioning,calc}
\usetikzlibrary{arrows.meta}
\tikzset{>=latex}
\tikzstyle{node}=[thick,circle,minimum size=22,inner sep=0.5,outer sep=0.6]
\tikzstyle{node in}=[node]
\tikzstyle{node hidden}=[node]
\tikzstyle{node convol}=[node]
\tikzstyle{node out}=[node]
\tikzstyle{connect}=[thick] 
\tikzstyle{connect arrow}=[-{Latex[length=4,width=3.5]},thick,shorten <=0.5,shorten >=1]
\tikzset{ 
  node 1/.style={node in},
  node 2/.style={node hidden},
  node 3/.style={node out},
}

\usetikzlibrary{fadings}
\pgfplotsset{compat=1.18}

\newcommand{\rhohat}{{\hat{\rho}}}
\newcommand{\rhoref}{\overline{\rho}}
\newcommand{\bbxhat}{{\hat{\bm{x}}}}
\newcommand{\bbmu}{\bm{\mu}}
\newcommand{\bbnu}{\bm{\nu}}
\newcommand{\bbtheta}{\bm{\theta}}
\newcommand{\thetavec}{\overrightarrow{\theta}}

\newcommand{\tildebbmu}{\Tilde{\bbmu}}
\newcommand{\phiref}{\hat{\Phi}}
\newcommand{\Omegaref}{\mathcal{R}}
\newcommand{\xhat}{\hat{x}}
\newcommand{\yhat}{\hat{y}}
\newcommand{\bbw}{\bm{w}}
\newcommand{\bbv}{\bm{v}}
\newcommand{\bbx}{\bm{x}}
\newcommand{\bbwref}{\overline{\bm{w}}}
\newcommand{\POD}{\textrm{POD}}
\newcommand{\bbalpha}{\bm{\alpha}}
\newcommand{\red}[1]{{\color{red}#1}}

\begin{document}

\maketitle

\begin{abstract}
We propose a novel Model Order Reduction framework that is able to handle solutions of hyperbolic problems characterized by multiple travelling discontinuities.
By means of an optimization based approach, we introduce suitable calibration maps that allow us to transform the original solution manifold into a lower dimensional one. {The novelty of the methodology is represented by the fact that the optimization process does not require the knowledge of the discontinuities location. The optimization can be carried out simply by choosing some reference control points, thus avoiding the use of some implicit shock tracking techniques, which would translate into an increased computational effort during the offline phase.}
In the online phase, {we rely on a non-intrusive approach,} where the coefficients of the projection of the reduced order solution onto the reduced space are recovered by means of an Artificial Neural Network. 
To validate the methodology, we present numerical results for the 1D Sod shock tube problem, for the 2D double Mach reflection problem, also in the parametric case, {and for the triple point problem}.
\end{abstract}

\begin{keywords}
hyperbolic problems, multiple travelling discontinuities, calibration map, neural network, model order reduction
\end{keywords}

\begin{MSCcodes}
65M08,35L60,35L67,76L05
\end{MSCcodes}

\section{Introduction}
The goal of MOR techniques \cite{benner2017model, stamm_rozza_hesthaven, haasdonk_ohlberger_2008}, which are particularly suited for the real-time computations and many-query context, is to obtain efficient and reliable approximations of solutions of high dimensional systems of partial differential equations (PDEs). 
Let us consider the approximate solution $u(t; \mu)\in L^2(\Omega)$ of a parametrized PDE, with $\Omega\subset \mathbb R^d$, with the parameter $\mu\in\mathcal{P}_{\text{phys}}\subset \mathbb R^s$ and with time $t\in[0,t_f]$: for the spatial discretization one can consider, for instance, the Finite Volume (FV) discretization. 
We introduce the \emph{solution manifold} related to this parametric PDE:
$\mathcal{M} = \{u(t; \mu)\in V_{\mathcal{N}},\text{ }\mu\in{\mathcal{P}_{\text{phys}}}, \text{ }t\in[0, t_{f}]\}, $ where $V_{\mathcal{N}}$ is a suitable functional space defined by the chosen spatial discretization.
The key idea behind MOR is to represent $\mathcal{M}$ with a finite dimensional linear space $V_N$, such that $N\ll\mathcal{N}$, where $N = \text{dim}V_N$ and $\mathcal{N}=\text{dim}V_{\mathcal{N}}$. 
To find the lower dimensional space $V_N$, one can use the well known Proper Orthogonal Decomposition (POD) strategy that, given in input a set of discrete solutions (obtained, for example, with the FV method), is able to extract a set of small cardinality $N$, which contains the so-called reduced basis functions that best approximate the manifold. 
A pivotal aspect for the efficiency of the MOR is the ability of the POD of compressing the discrete solution manifold $\mathcal{M}$: this concept is strictly related to the definition of Kolmogorov $N$-width of $\mathcal{M}$ and, ultimately, to the reducibility of the problem of interest.

The \emph{Kolmogorov $N$-width} $d_N$ of $\mathcal{M}$ is defined as
\begin{equation}
\label{eq:Kolmogorov n-width}
    d_N(\mathcal{M}, V_{\mathcal{N}}) = \inf_{\substack{V_N\subset V_{\mathcal{N}}\\ \text{dim}V_N=N}}\sup_{f\in\mathcal{M}}\inf_{g\in V_N}\lvert\lvert f-g\rvert\rvert,
\end{equation}
where $\lvert\lvert\cdot\rvert\rvert$ is a suitable norm in $V_{\mathcal{N}}$.
Definition \eqref{eq:Kolmogorov n-width} describes in a rigorous mathematical setting the capability of finite dimensional linear subspaces $V_N\subseteq V$ of reproducing any element in $\mathcal{M}$, that is, any discrete solution of the problem of interest. 
Therefore, the faster $d_N$ decays the more efficient a linear MOR will be for such a problem, as $N$ grows.
Some rigorous bounds for $d_N$, for particular classes of problems, are available in literature \cite{CohenDeVore, MELENK2000272}; as an alternative, one can look at the rate of decay of the eigenvalues $\lambda_k$ returned by the POD on $\mathcal{M}$.

Despite the capability of standard MOR techniques to handle a vast number of applications, problems advecting local structures still represent a challenge for the MOR community. 
Indeed, for such problems the decay of the Kolmogorov $N$-width is slow, see for example \cite{GREIF2019216}. As a result, standard MOR struggles to suitably reproduce steep features, such as solutions with (multiple) travelling shocks. 
For this reason, in the last decade a great number of works appeared in the literature, offering numerous approaches to deal with advection dominated problems. 
We mention the method of freezing \cite{OHLBERGER2013901}, the shifted POD \cite{shifted_pod_reiss} ({also in combination with the use of neural networks \cite{nn_shifted_pod}}), the generation of advection modes by means of an optimal mass transport problem \cite{optimal_mass_transfer_iollo,iollo2022mapping,khamlich2023optimal,battisti2023wasserstein}, $L^1$ minimization \cite{abgrall2016robust}, the calibration method introduced in \cite{cagniart_calibration, cagniart_hal}, Lagrangian based MOR techniques \cite{mojgani_pinns, Mojgani2017}, the preprocessing of the snapshots used in \cite{karatzas, nonino}, the registration method \cite{taddei_sisc, taddei_jcp,barral2023registration}, adaptive basis methods \cite{peherstorfer2020model}, implicit feature tracking \cite{mirhoseini2023model} and displacement interpolation \cite{trasport_greedy,rim18displacementInterpolation}. 
Next to these more classical techniques, some nonlinear approaches have been lately studied starting from convolutional autoencoders neural network based approaches for learning the solution manifold \cite{carlberg_lee_jcp, fresca_dede_manzoni}, passing through graph neural networks autoencoders \cite{pichi2023graph} to graph neural network to perform the limit to vanishing viscosity \cite{romor2023friedrichs} {and entropy-stable rational quadratic manifolds \cite{klein2024entropy}}.
Motivated by the interest that MOR for transport dominated problems sparks in the applied mathematics community, the goal of this work is to propose a calibration based reduced order algorithm that can be used to gain significant speedup in the simulations of solutions of hyperbolic PDEs. 
{The proposed methodology will avoid the use of implicit feature tracking \cite{mirhoseini2023model}, with the goal of not increasing the computational costs of the offline phase that is performed explicitly.
The framework is very similar to the calibration/registration approach proposed in \cite{taddei_sisc,taddei_jcp,cagniart_calibration} with novelties in the calibration process and in the range of applicability of the method. Indeed,
our methodology has a more physical intuition as it is based on the interpolation of some calibration points rather than on more abstract interpolation functions.
Moreover, the proposed applications range on a novel paradigm, where parametric time-dependent solutions widely vary, with large deformations of the solutions structures.
To the knowledge of the authors, this is the first time that solutions with such a varying behavior were successfully tackled with a model order reduction technique.}

In this work, we will focus on time-dependent hyperbolic problems, whose solutions are (quasi) self-similar: the formal definition of (quasi) self-similar solution will be given in Section~\ref{sec:motivation}. 
In particular, we want to study problems where multiple structures travel along the domain with different speeds. 
This is typical for hyperbolic problems, where shocks, rarefactions and other discontinuities are generated and travel along the domain.
The novelties of this work, in comparison to the state of the art \cite{taddei_jcp}, lie in two key aspects. First, our optimization process operates independently of the solution structure, eliminating the need for shock detectors or similar tools. Secondly, our method demonstrates a broader range of applicability, encompassing problems featuring multiple shocks whose positions undergo significant variation, sometimes nearly colliding with each other.

The rest of the manuscript is structured as follows.
In Section~\ref{sec:motivation}, we introduce the problems of interest and the definition of self-similarity. In Section~\ref{sec:calibration}, we define the calibration procedure and the  optimization algorithms.
In Section~\ref{sec:geom_trans}, we present some geometrical transformations that interpolate the calibrated points and allow to define the original problem onto a reference domain where all the structures do not move.
In Section~\ref{sec:MOR with calibration}, we describe the combination of classical MOR techniques with the calibration process. 
In Section~\ref{sec:numerical}, we show the good performances of the proposed reduced order model (ROM) onto one and two dimensional parametric time--dependent problems and, in Section~\ref{sec:conclusion}, we draw some conclusions. {To improve the readibility of the manuscript, a table with the list of all the mathematical symbols used is provided at the end of the paper.}
 
\subsection{Motivation}\label{sec:motivation}
We begin by introducing the problem of interest that we will tackle in this manuscript: we focus here on hyperbolic time--dependent conservation laws. As an example, we turn our attention to Euler equations, but the same framework can be applied to other conservation and balance laws. 
Let $\Omega\subset\mathbb{R}^d$, $d\geq 1$, be our physical domain. We restrict ourselves to rectangular domains of the type $\Omega=[a^1,b^1]\times\dots\times[a^d, b^d]$, with $a^i, b^i\in\mathbb{R}$ for $i=1,\dots, d$. 
The generalization for more complex domains can be performed as in \cite{taddei2023compositional}. 
Let $[0,t_f]\subset \mathbb R$ be the time span of the problem and let $\bbmu\in\mathcal{P}\subset\mathbb{R}^{s+1}$, $s\geq 0$, be the collection of all parameters (including time). 
From now on, we will assume that $s=0$ in the non parametric regime, i.e., $\bbmu = t$ and $\mathcal{P}=[0, t_f]$, or $s\geq 1$ in the parametric regime, i.e., $\bbmu = (\mu,t)$ and $\mathcal{P}=\mathcal{P}_{\text{phys}}\times [0, t_f]$.  
The parameteric Euler equations of gas dynamics, in conservative form, read as follows: find the density $\rho\colon\mathcal{P}\times\Omega\mapsto\mathbb{R}$,  the momentum $\boldsymbol{m} \colon\mathcal{P}\times\Omega\mapsto\mathbb{R}^d$ and the total energy $E\colon\mathcal{P}\times\Omega\mapsto\mathbb{R}$ such that
\begin{equation}
\label{eq:euler conservative}
    \begin{cases}
        \partial_t\rho + \nabla_{\bm{x}} \cdot \boldsymbol{m} = 0 &\text{ in }\mathcal{P}\times\Omega,\\
        \partial_t\boldsymbol{m} + \nabla_{\bm{x}} \cdot \left(\frac{\boldsymbol{m} \otimes \boldsymbol{m}}{\rho} + p\,I\right) = 0 &\text{ in }\mathcal{P}\times\Omega,\\
        \partial_tE + \nabla_{\bm{x}} \cdot \left(\frac{\boldsymbol{m}}{\rho}(E + p)\right) =  0 &\text{ in }\mathcal{P}\times\Omega,
    \end{cases}
\end{equation}
where $\nabla_{\bm{x}}\cdot $ is the divergence with respect to $\bm{x}\in\Omega$, $I\in\mathbb R^{d\times d}$ is the identity matrix and the pressure $p$ is defined through the following equation of state $ p = (\gamma-1)(E-0.5 |\boldsymbol{m}|^2/\rho )$, with $\gamma = 1.4$ being the adiabatic constant. System \eqref{eq:euler conservative} is then completed by some proper initial conditions (IC) and boundary conditions (BC). We will consider as IC some Riemann problems both in one and two dimensional problems: the Sod shock tube problem \cite{toro}, the 2D double Mach reflection problem \cite{woodward} {and the triple point problem}. 
In these examples, the solution of \eqref{eq:euler conservative} turns out to be self-similar, with features as shocks, contact discontinuities and rarefaction waves traveling in the physical domain.
\begin{definition}
    Let $\Omegaref\subset\mathbb{R}^d$ be a reference domain, which is time (and parameter) independent. We call \textbf{self--similar} a solution manifold $\mathcal M$ for which there exists a reference solution $\bar{u}:\Omegaref \to \mathbb R^{d}$ and a transformation $T^{-1}[\bbmu](\cdot) : \Omega \to \Omegaref $ such that we have $u(\bbmu)(T^{-1}[\bbmu](\bbx)) \approx \bar{u}(\bbx)$ for all $\bbmu\in\mathcal{P},$ $\forall \bbx \in \Omegaref$. When  this condition is not satisfied, but still all solutions in $\mathcal M = \lbrace u(\bbmu) \rbrace_{\bbmu\in\mathcal{P}}$ present the same features, with different values of the solution in between these features, we will call such manifold {\textbf{quasi--self--similar}}. More precisely, $\mathcal M$ is  \textbf{quasi--self--similar} if there exists a transformation $T$ such that the transformed solution manifold $\hat{\mathcal{M}}:=\lbrace u(\bbmu)(T^{-1}[\bbmu](\cdot)) \rbrace$ has a fast decay of the Kolmogorov $N$-width.
\end{definition}

\begin{figure}
    \centering
    \includegraphics[width=0.45\textwidth,trim={25 20 35 40}, clip]{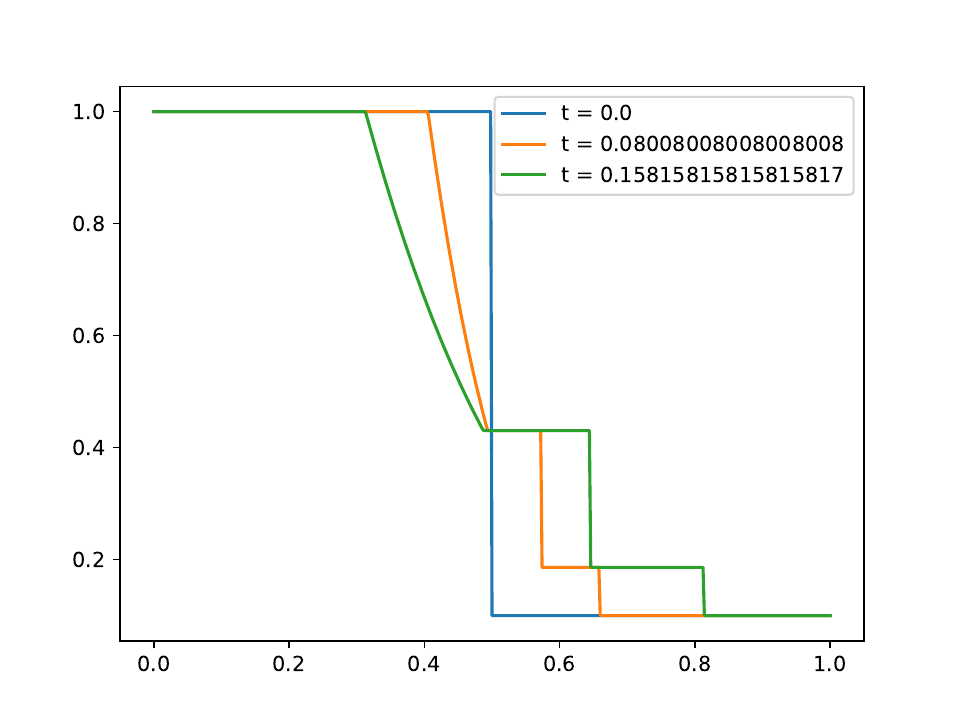}\hfill
    \includegraphics[width=0.45\textwidth,trim={25 20 35 40}, clip]{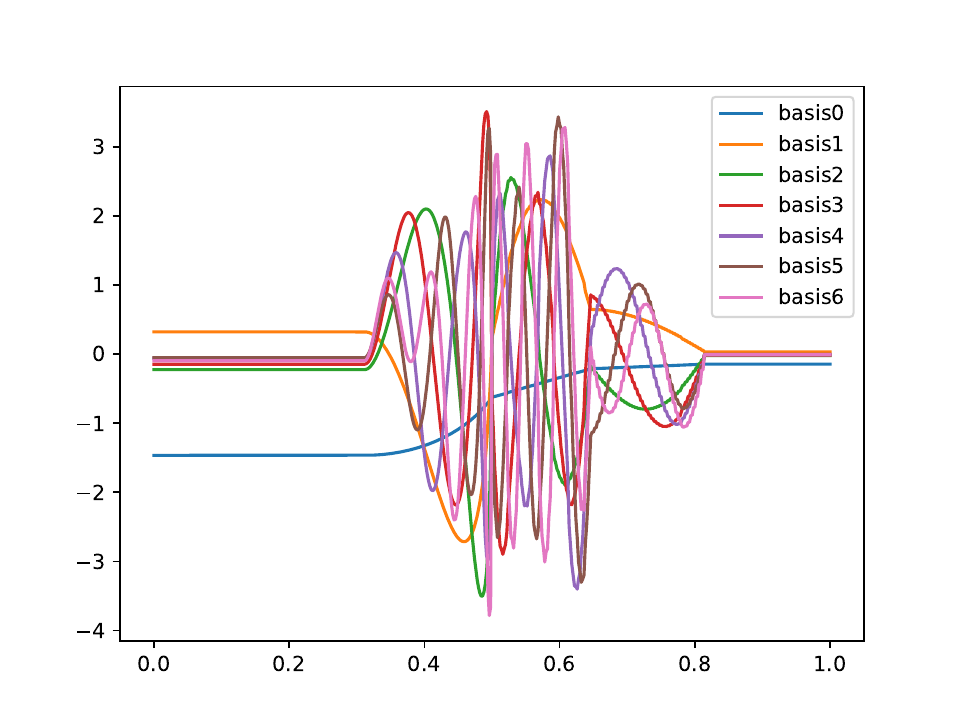}
    \caption{The density $\rho$, solution of the Sod shock tube problem at different timesteps, and the corresponding POD modes (right).}
    \label{fig:exact sod results}
\end{figure}

We start by considering a simple 1D Sod shock tube problem, in the non--parametric regime. 
Here, $\Omega = [0,1]$, and the following initial data is considered:
\begin{equation*}
    \begin{bmatrix}
    \rho & u & p 
    \end{bmatrix}(t=0)
    = 
    \begin{cases}
        [1 \quad 0 \quad 1]^T &\text{ if } x< 0.5,\\
        [0.1 \quad 0 \quad 0.125 ]^T &\text{ if } x> 0.5,\\
    \end{cases}
\end{equation*}
where $u\colon\Omega\mapsto\mathbb{R}^d$ is the velocity. The initial conditions in conservative variables $m$ and $E$ can be derived using $m=\rho u $ and $E = p/(\gamma-1) + 1/2 |m|^2/\rho$.
Figure~\ref{fig:exact sod results} shows the density $\rho(\bbmu)$ for the Sod shock tube problem (left), and the corresponding modes (right), obtained by running a POD on the solution manifold $\mathcal{M}_{\rho}$: the solution of the Sod shock tube problem is exact, and its analytical expression has been taken from \cite{toro}. 
The density presents a shock, a contact discontinuity and a rarefaction wave that travel in the domain: as a consequence, the POD modes exhibit an highly oscillatory behavior, struggling to correctly capture the position of the moving features.
Still, the solutions are \textbf{self-similar} as we just need to transport each feature onto reference positions {to make the solutions essentially linearly dependent}. 
Indeed, as observed in \cite{iollo2022mapping}, the optimal transport for the density of the Sod 1D problem would lead to the exact solution, without the need for further ROM techniques. 
Nevertheless, we are interested also in higher dimensional problems with more complicated discontinuity structures.
Hence, we will not proceed in the optimal transport direction.

\section{Calibration of the snapshots}\label{sec:calibration}
We now present the calibration technique that we use to align the different features of our snapshots, to obtain a solution manifold with a faster decaying Kolmogorov $N$-width. 
The key of the proposed calibration is that it can be used to align different travelling features (shocks, contact discontinuities, rarefaction waves), without the need to know explicitly the exact location of these features, as opposed to, for example, what is assumed in \cite{nonino,taddei_jcp}. 
Moreover, we assume to calibrate the density $\rho$ of the Euler system~\eqref{eq:euler conservative}, other scalar quantities depending on the system unknowns, e.g., the entropy, can be equivalently used.

\begin{figure}
    \centering
\adjustbox{max width = \textwidth}{
\begin{tikzpicture}[scale=1.8]
    \tikzset{dot/.style={fill=black,circle}}
    \tikzset{straightness/.style={straightness, draw, 
            minimum size=3pt, 
            inner sep=0pt, outer sep=0pt}}
    \tikzset{arrow/.style={
            color=black,
            draw=black,
            -latex,
            font=\fontsize{8}{8}\selectfont}}
    \newcommand{\xa}{1}
    \newcommand{\xb}{2.5}
    \newcommand{\xc}{4}
    \newcommand{\xd}{5.5}
    \newcommand{\xe}{7}

    \foreach\l[count=\y] in {0.5,1,1.5,2}
    {
        \draw (1,\l) -- (4,\l);
        \foreach\w/\z[count=\x] in {\xa/1,\xb/2,\xc/3}
        {
            \fill (\w,\l) circle (0.5mm) node[anchor=south west] {$(\bbwref^1_{\z,\y},\bbwref^2_{\z,\y})$};
        }
    }

    \foreach\l/\w[count=\x] in {0/\xa,1/\xb,2/\xc}
    {
        \draw (\w,0.5) -- (\w,2.);
    }

    \newcommand{\yzero}{0.25}
    \draw (\xa,\yzero) -- (\xc,\yzero);
    \foreach\w[count=\x] in {\xa,\xb,\xc}
    {
        \draw  (\w,\yzero)  node[anchor=north]  {$\bar{x}_\x$};
        \draw  (\w,\yzero+0.05) -- (\w,\yzero-0.05);
    }

    \newcommand{\xzero}{0.7}
    \draw (\xzero, 0.5 ) -- (\xzero,2.);
    \foreach\l[count=\y] in {0.5,1,1.5,2}
    {
        \draw  (\xzero,\l)  node[anchor=east]  {$\bar{y}_\y$};
        \draw  (\xzero-0.05,\l) -- (\xzero+0.05,\l);
    }

    \renewcommand{\xa}{6}
    \renewcommand{\xb}{7.5}
    \renewcommand{\xc}{9}
    \newcommand{\xba}{7.7}
    \newcommand{\xbb}{7.2}
    \newcommand{\xbc}{7.54}
    \newcommand{\xbd}{7.1}
    \newcommand{\ya}{0.5}
    \newcommand{\yb}{1}
    \newcommand{\yc}{1.5}
    \newcommand{\yd}{2}

    \newcommand{\yba}{1.1}
    \newcommand{\ybb}{0.9}
    \newcommand{\ybc}{1.}

    \newcommand{\yca}{1.55}
    \newcommand{\ycb}{1.73}
    \newcommand{\ycc}{1.67}

    \draw (\xa,\ya) -- (\xc,\ya);
    \draw (\xa,\yd) -- (\xc,\yd);
    \draw (\xa,\ya) -- (\xa,\yd);
    \draw (\xc,\ya) -- (\xc,\yd);

    \draw plot [smooth] coordinates {(\xa, \yba) (\xba, \ybb) (\xc, \ybc)};
    \draw plot [smooth] coordinates {(\xa, \yca) (\xbc, \ycc) (\xc, \ycb)};
    \draw plot [smooth] coordinates {(\xbb, \ya) (\xba, \ybb) (\xbc, \ycc) (\xbd, \yd)};

    \fill (\xa,\ya) circle (0.5mm) node[anchor=south west] {$(\bbw^1_{1,1},\bbw^2_{1,1})$};
    \fill (\xbb, \ya) circle (0.5mm) node[anchor=north west] {$(\red{\bbw^1_{2,1}},\bbw^2_{2,1})$};
    \fill (\xc, \ya) circle (0.5mm) node[anchor=south west] {$(\bbw^1_{3,1},\bbw^2_{3,1})$};

    \fill (\xa, \yba) circle (0.5mm) node[anchor=south west] {$(\bbw^1_{1,2},\red{\bbw^2_{1,2}})$};
    \fill (\xba, \ybb) circle (0.5mm) node[anchor=south west] {$(\red{\bbw^1_{2,2}},\red{\bbw^2_{2,2}})$};
    \fill (\xc, \ybc) circle (0.5mm) node[anchor=south west] {$(\bbw^1_{3,2},\red{\bbw^2_{3,2}})$};

    \fill (\xa, \yca) circle (0.5mm) node[anchor=south west] {$(\bbw^1_{1,3},\red{\bbw^2_{1,3}})$};
    \fill (\xbc, \ycc) circle (0.5mm) node[anchor=north west] {$(\red{\bbw^1_{2,3}},\red{\bbw^2_{2,3}})$};
    \fill (\xc, \ycb) circle (0.5mm) node[anchor=north west] {$(\bbw^1_{3,3},\red{\bbw^2_{3,3}})$};

    \fill (\xa, \yd)  circle (0.5mm) node[anchor=south west] {$(\bbw^1_{1,4},\bbw^2_{1,4})$};
    \fill (\xbd, \yd) circle (0.5mm) node[anchor=south west] {$(\red{\bbw^1_{2,4}},\bbw^2_{2,4})$};
    \fill (\xc, \yd)  circle (0.5mm) node[anchor=south west] {$(\bbw^1_{3,4},\bbw^2_{3,4})$};
            
\end{tikzpicture}    }
\vspace{-5mm}
\caption{Example of a control point grid in 2D with $M_1=3$ and $M_2=4$ on the reference domain (left) and on the physical domain (right). Note that the coordinates of the reference control points are the tensor product of unidirectional control points. We highlight in red the free coordinates of these control points}\label{fig:ref_grid}
\end{figure}
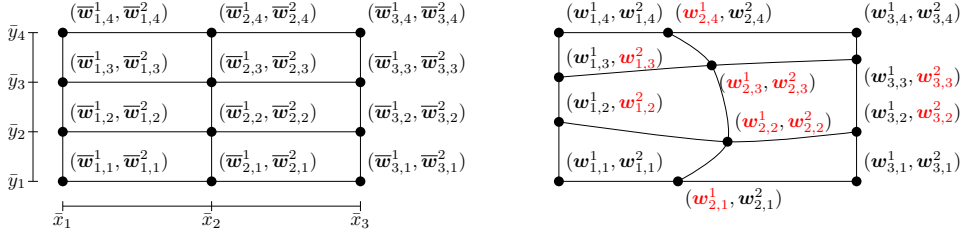

Let $\Omegaref$ be the \emph{reference domain}, similarly to what is used in the Arbitrary Lagrangian Eulerian (ALE) formalism \cite{torlo}, and let $\Omega$ be the \emph{physical domain}. 
For every $\bbmu\in\mathcal{P}\subset\mathbb{R}^{s+1}$, we introduce a grid of $M=\prod_{i=1} ^d M_i$ \textit{control points} that are collected in the vector $\bbw(\bbmu)\in\Omega^M$ that we use for the calibration. 
These control points should lead the transformation to align the different features at different $\bbmu$.
Let $\bbalpha=(\alpha_1,\dots,\alpha_d)\in\mathbb N^{d}$ be a multi-index with $\alpha_i =1,\dots, M_i$. {We consider the control point $\bbw_{\bbalpha}$ to be tensor product of points in the interval $[a^i,b^i]$ for every dimension $i=1,\dots,d$, for example in 2D: $\bbw_{\bbalpha = (\alpha_1,\alpha_2)} = (\bbw^1_{\alpha_1,\alpha_2},\bbw^2_{\alpha_1, \alpha_2})$, see Fig.~\ref{fig:ref_grid}.} Each control point $\bbw(\bbmu)_{\bbalpha}$ can either belong to the physical domain, i.e., $\bbw(\bbmu)_{\bbalpha}\in \ring{\Omega}$, or to the boundary of the domain, namely $\bbw(\bbmu)_{\bbalpha}\in\partial\Omega$. 
{If $\bbw(\bbmu)_{\bbalpha}\in\partial\Omega$, then this point is constrained for all $\bbmu$ to the boundary hyperplanes where it belongs: this, in turn, sets a constraint on one (or both) of the coordinates $(\bbw^1_{\alpha_1,\alpha_2},\bbw^2_{\alpha_1, \alpha_2})$. The constrained coordinates are, hence, not to be optimized and we can discard them in the set of coordinates that we will optimize. 
Motivated by this, we introduce $\bbtheta(\bbw(\bbmu))\in \mathbb R^{Q}$, that is the row vector of the $Q$ \emph{free coordinates} of the control points, with $Q\leq d\times M$. This vector represents all the coordinates of the control points that are free to move during the calibration. In Fig.~\ref{fig:ref_grid}, the free coordinates of the control points are highlighted in red.} We remark that there is a bijection between $\bbw$ and $\bbtheta$, by definition.\\
In order to align different features of our set of snapshots, we look for a geometrical transformation map $T\colon \Omega^M \times \Omegaref\mapsto \Omega$, such that the following properties hold true: 
\begin{itemize}
    \item $T[\cdot] \in \mathcal{C}^1(\Omega^M, \mathcal{C}^1(\Omegaref, \Omega))$;
    \item $\forall\bbmu\in\mathcal{P}$ and $\forall\bbw(\bbmu)\in{\Omega^M}$, $\exists T^{-1}[\bbw(\bbmu)]\colon \Omega \mapsto \Omegaref$ such that 
    \begin{align*}
        &T^{-1}[\bbw(\bbmu)](T[\bbw(\bbmu)](\bbxhat)) = \bbxhat \quad \forall \bbxhat\in\Omegaref,\\
        &T[\bbw(\bbmu)](T^{-1}[\bbw(\bbmu)](\bm{x})) = \bm{x},\text{ for all }\bm{x}\in\Omega;
    \end{align*}
    \item $T^{-1}[\cdot]\in\mathcal{C}^1(\Omega^M, \mathcal{C}^1(\Omega; \Omegaref))$.
\end{itemize}
The properties are imposed to setup an ALE formulation \cite{torlo,taddei_sisc}.
Some possibilities for the geometrical transformation $T$ have been presented in the literature over the past years: among the others we mention here translation maps, dilatation maps, polynomials and Gordon-Hall maps, see \cite{torlo,cagniart_calibration,taddei_sisc}. We will use some transformations based on Piecewise Cubic Hermite Interpolation Polynomials (PCHIPs) carefully described in Section~\ref{sec:geom_trans}.  
In order to use the transformation map $T$, one needs to find the calibration map $\bbw\colon\mathcal{P}\mapsto\Omega^M$, such that:
\begin{itemize}
    \item $\bbw(\cdot)\in\mathcal{C}^1(\mathcal{P}; \Omega^M)$; 
    \item $\rho(\bbmu)( T[\bbw(\bbmu)](\bbxhat)) {\approx} \rhoref(\bbxhat)$, for all $\bbxhat\in\mathcal R$, $t\in[0, t_{f}]$ and for all $\bbmu\in\mathcal{P}$, where $\rhoref(\cdot)$ is a reference solution of choice.
\end{itemize}
We are now ready to present a general calibration technique: the ultimate goal is to transform the solution manifold $\mathcal{M}_{\rho}$ so that the POD, applied to the transformed manifold, is more efficient. In order to achieve this goal, we need to perform the following steps.
First of all, we introduce the \emph{reference} density $\rhoref$, namely a solution of problem \eqref{eq:euler conservative} for some $\overline{\bbmu}\in\mathcal{P}$.
Once $\rhoref$ has been chosen, we select $M$ control points $\bbwref \in\Omegaref^M$.
Now, for any $\bbmu\in\mathcal{P}$ and $\bbw(\bbmu)\in\Omega^M$, we can define a geometrical transformation map $T[\bbw(\bbmu)]\colon\Omegaref \mapsto\Omega$ {that maps the reference control points in the reference domain onto the (parametric) control points in the physical domain, i.e.,}
\begin{equation*}
\begin{split}
   T[\bbw(\bbmu)] (\bbwref_{\bbalpha}) =  \bbw_{\bbalpha}(\bbmu) \quad\forall\, \bbalpha.
\end{split}
\end{equation*}
Once $T[\bbw(\bbmu)]$ has been defined, we can introduce the \emph{calibrated snapshot}{, which is the pullback of a solution $\rho(\bbmu)$ to the reference domain, i.e.,} $$\rhohat(\bbmu; {\bbw(\bbmu)})(\cdot):=\rho(\bbmu)(T[\bbw(\bbmu)](\cdot)):\mathcal{R}\to \mathbb R.$$ 

We want to stress that, numerically, we rely on two meshes, one on the physical and one on the reference domain. 
In our simulations, these two domains will coincide; this does not mean that each transformation map leads to a one-to-one correspondence between the degrees of freedom on the two domains. 
Hence, we will perform an interpolation of $\rho$ to evaluate $\hat{\rho}$ at its degree of freedom. 
In the numerical simulations, this procedure will bring an error of the first order of accuracy (as we will use FV approximations). 

The map $T[\bbw(\bbmu)]$ is identified by the control points $\bbw(\bbmu)$, which are sought in order to minimize the following residual function:
\begin{equation}
    \label{eq:residual definition}
    \begin{split}
        R(\rho(\bbmu),\bbtheta(\bbw(\bbmu)), \rhoref) = &\lvert\lvert \rhohat (\bbmu; {\bbw(\bbmu)}) - \rhoref\rvert\rvert^2_{L^2(\Omegaref)} + \frac{\delta}{2}\left\lVert\partial_{\bbmu}\bbw(\bbmu)\right\rVert^2_{\ell^2(\mathcal{P})}+\\
    &\frac{\alpha}{2} \max_{\bm{x}\in\Omega} \left( \max \left\lbrace \left\lVert\nabla T[\bbw(\bbmu)](\hat{\bm{x}})\right\rVert,  \left\lVert\nabla T^{-1}[\bbw(\bbmu)](\bm{x})\right\rVert \right\rbrace \right) ,
    \end{split}
\end{equation}
where $\alpha$ and $\delta$ are two penalty parameters user defined, and $\partial_{\bbmu}\bbw(\bbmu)$ will be defined more in detail in the algorithmic section. {The first term is the one we aim at minimizing, while the other two are regularization terms that penalizes discontinuities in time of the calibration points and that regularize the geometrical transformations.}
We are now ready to present the calibration technique in two different cases: the self-similar setting, and the quasi self-similar one.

\subsection{Calibration in the self--similar setting}\label{sec:calibration self similar}
To keep the presentation as general as possible, we consider here the parametric case, hence $\bbmu=(\mu, t)\in\mathcal{P}=\mathcal{P}_{\text{phys}}\times [0, t_f]$.
We select a training set of physical parameters $\lbrace \mu_1,\dots, \mu_{N_{\text{train}}}\rbrace  = \mathcal{P}_{\text{phys}}^{\text{train}} \subset \mathcal{P}_{\text{phys}}$ and a training set of times $\lbrace t_1,\dots, t_{N_\mu}\rbrace$ for each $\mu \in \mathcal{P}_{\text{phys}}^{\text{train}}$.
We then denote $\mathcal{P}^{\text{train}}:=\lbrace \bbmu = (\mu,t) : \mu \in \mathcal{P}^{\text{train}}_{\text{phys}} \text{ and } t \in [t_1,\dots, t_{N_\mu}]  \rbrace\subset \mathcal{P}$. 
For each $\bbmu\in \mathcal{P}^{\text{train}}$, we compute the full order model solution $\rho(\bbmu)$. We then choose as reference solution $\rhoref=\rho(\bar{\bbmu})$: in our numerical tests we will choose $\bar{\bbmu}~= (\mu_{N_{\text{train}}},t_f)$. In addition to this, we also fix the $M$ reference control points $\bbwref\in\Omegaref^M$ as specified above.
{Before introducing the constrained minimization problem, let us recall that there exist a bijection between $\bbtheta(\bbw(\bbmu))$ and $\bbw(\bbmu)$. For this reason, in order to keep the notation light, from now on we will denote with $\bbtheta(\bbmu)$ the free coordinates $\bbtheta(\bbw(\bbmu))$ of the control points $\bbw(\bbmu)$.}\\
Now, for all $\bbmu$ in the training set, we solve the following constrained minimization problem:
\begin{equation}
    \label{eq:calibration self similar}
    \bbtheta^{\text{opt}}(\bbmu):=\min_{\bbtheta\in \mathbb{R}^Q} R(\rho(\bbmu), \bbtheta; \rhoref),
\end{equation}
subject to the following constraints:
\begin{itemize}
    \item all the control points are within our physical domain: $\bbw_{\bbalpha}(\bbmu)\in\Omega$ for all $\bbalpha$  and for all training parameters $\bbmu\in\mathcal{P}^{\text{train}}$;
    \item $\det J[\bbw(\bbmu)]>0$, where $J[\bbw(\bbmu)]$ is the Jacobian of the map $T[\bbw(\bbmu)]$. This constraint must be checked on a quite fine grid of the physical domain $\Omega$, we have used the mesh grid;
    \item for $i=1,\dots, d$: if $\alpha_j=\beta_j$ for all $j\neq i$ and $\alpha_i < \beta_i$, then $\bbw_{\bbalpha}^i<\bbw_{\bm{\beta}}^i$ (see Fig.~\ref{fig:ref_grid}), i.e., we never switch the order of the control points on each grid line.
\end{itemize}

We approximate $\partial_{\bbmu}$ in \eqref{eq:residual definition} with the discrete derivative $D_{\bbmu}\bbw(\bbmu)$ defined as:
\begin{equation*}
D_{\bbmu}\bbw^i_{\bbalpha}(\bbmu) =
\begin{cases}
    0 &\text{ if }\bbw^i_{\bbalpha}(\bbmu) \text{ is not a free coordinate of }\bbw(\bbmu)  ; \\
    \frac{\bbtheta_q(\bbmu)-\bbtheta^{\text{opt}}_q(\bbmu_{\text{neigh}}(\bbmu))}{\bbmu - \bbmu_{\text{neigh}}(\bbmu)} &\text{with } q \text{ s.t. } \bbw^i_{\bbalpha}(\bbmu) = \bbtheta_q(\bbmu).
\end{cases}
\end{equation*}
In the previous equation, the \emph{neighboring parameter} $\bbmu_{\text{neigh}}$ is defined as follows:
\begin{equation*}
    \bbmu_{\text{neigh}}(\bbmu, \mathcal{S}) = \arg\min_{\bbnu\in\mathcal{S}} \lvert\lvert\bbmu-\bbnu\rvert\rvert_{\ell^2(\mathbb{R}^{p+1})},
\end{equation*}
where $\mathcal{S} \subset \mathcal{P}^{\text{train}}$ will contain the parameters for which we have already computed the optimal $\bbtheta$. If the minimum is not unique, we take one of the minimizers.
The definition of the discrete spatial gradient $\nabla T$ in \eqref{eq:residual definition} will be specified in Section~\ref{sec:geom_trans} for the specific transformation map we use.

Problem~\eqref{eq:calibration self similar} is solved with the Sequential Least SQuares Programming (SLSQP) method that is available within the \texttt{scipy.optimize.minimize} library. 
We solve Problem~\eqref{eq:calibration self similar} for the physical parameter $\mu\in\mathcal{P}^{\text{train}}_{\text{phys}}$ for which the solution has more developed structures (the last in our tests): we start from the final timestep $t_{N_{\mu}}$ and we proceed backwards in time. We then move onto solving Problem~\eqref{eq:calibration self similar} for the neighboring parameters. In both physical and temporal parameters, the rationale is the following: we proceed form the solutions where the structures that we want calibrate are more developed and we proceed with nearest parameters until we solve the problem for all the training set $\mathcal{P}^{\text{train}}_{\text{phys}}$.
The initial guess $\bbtheta^{(0)}(\bbmu_{(\ell,j)})$ for $\bbmu_{(\ell,j)} = (\mu_{\ell}, t_j) \in \mathcal{P}^{\text{train}}$ is the optimal output of the minimization for the closest parameter already performed. In our tests it will be defined as:
\begin{equation}
\label{eq:initial guess parametric setting}
\bbtheta^{(0)}(\bbmu_{(\ell,j)})=
    \begin{cases}     
    \bbtheta^{\text{opt}}(\bbmu_{(\ell,j+1)}) \text{ if }j=1,\dots, N_{\mu_{\ell}}-1,\\
    \bbtheta^{\text{opt}}(\bbmu_{(\ell+1,j)}) \text{ if }j=N_{\mu_{\ell}}.
    \end{cases}
\end{equation}
Algorithm~\ref{alg:calibration self similar} shows the details of the procedure.
\begin{algorithm}
\caption{Calibration for self--similar solution manifold}\label{alg:calibration self similar}
\begin{algorithmic}[1]
    \State \textbf{Input:} the reference solution $\rhoref$, the control points $\bbwref\in\mathcal{R}^M$ and the training solution manifold on the physical domain $\lbrace \rho(\mu_\ell, t_j)\text{ for } \ell=1,\dots, N_{\text{train}}, j = 1, \dots,N_{\mu_\ell} \rbrace$.
    \For{$\ell=N_{\text{train}}, \dots, 1$}
        \For{$j=N_{\mu_{\ell}},\dots, 1$}
            \State Set the initial guess $\bbtheta^{(0)}(\bbmu_{( \ell,j)})$ as in \eqref{eq:initial guess parametric setting}.
            \State Solve Problem \eqref{eq:calibration self similar}.
        \EndFor
    \EndFor
    \State \textbf{Output:} the optimal control points $\bbw^{\text{opt}}(\bbmu_{(\ell,j)})$.
\end{algorithmic}
\end{algorithm}

{The chosen strategy is not trivially favorable with respect to other ones and it is quite arbitrary. 
We have performed different trials in the tests that we have carried out, and we noticed differences in the sensitivity to this choice. 
Overall, taking as initial guess the optimal solution of a parameter close to the one we want to optimize is always a good strategy, but open questions remain on which reference configuration $\rhoref$ to choose and on how spread can the parameter of the training set be.

For the 1D test cases, the change in the choice of the initial guess was less relevant than for the 2D cases. Indeed, we could use initial guesses that are quite far away from the exact one and still obtain very accurate results.
For the 2D cases, the choice of the initial guess, and how spread the training space can be is more relevant.
For example, the timesteps used for the calibration and, hence, the distance between two consecutive parameters and the initial guess, had a strong impact. We therefore had to calibrate on timesteps that are not too far away one from the other, in order to obtain smooth and continuous results.

For the 1D test, we perform an analytical study of this choice where the exact solution is available: the outcome of this study is presented in Section~\ref{sec:validation_calibration}.}

\subsection{Calibration in the quasi--self--similar setting}\label{sec:calibration quasi self similar} 
\begin{figure}
    \centering
    \includegraphics[width=0.45\textwidth]{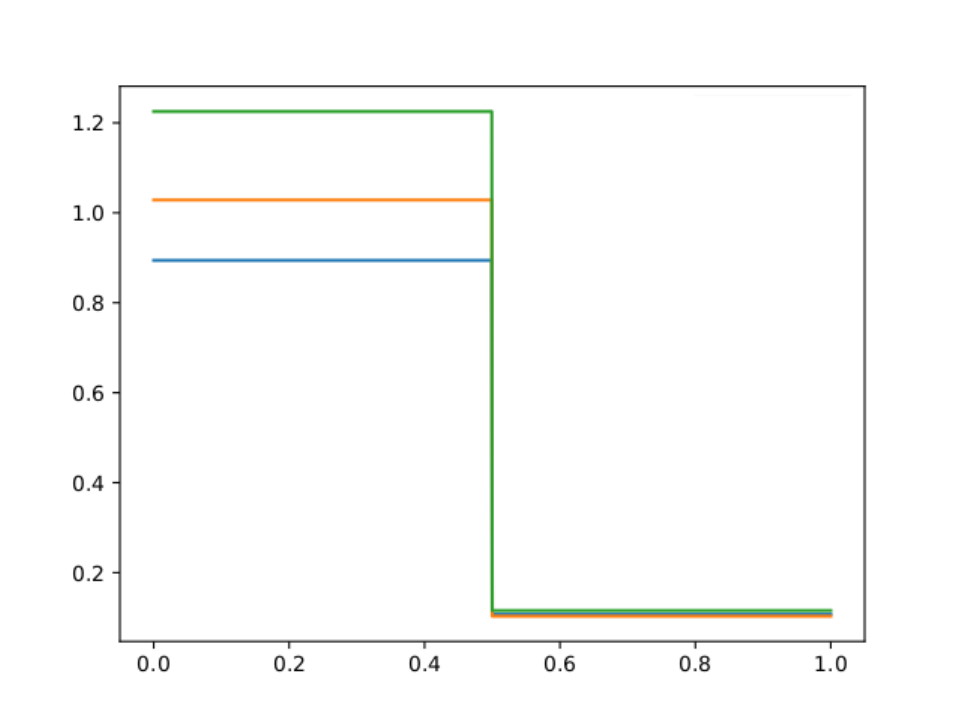}
    \vspace{-2mm}
    \caption{The density $\rho(t=0)$ for three different values of the parameter $\mu=[0.8943, 0.1075, 1.0906, 0.0572]$ (blue), $\mu=[1.0286, 0.1031, 0.7358, 0.0706]$ (orange), and $\mu=[1.2253, 0.1157, 1.1172, 0.1094]$ (green)}
    \label{fig:motivational_sod_more_params}
\end{figure}
To motivate the need of a different algorithm for quasi--self--similar solutions, we focus now on the parametric version of the 1D Sod shock tube problem \eqref{eq:euler conservative}. In this example, the physical parameter $\mu=(\mu^0,\dots,\mu^3)\in\mathcal{P}_{\text{phys}}\subset\mathbb{R}^4$ represents the IC for the Euler problem:
\begin{equation*}
    \begin{bmatrix}
    \rho & u & p 
    \end{bmatrix}^T(t=0; \mu)
    = 
    \begin{cases}
        \begin{bmatrix}
            \mu^0 & 0 & \mu^1
        \end{bmatrix}^T &\text{ if } x< 0.5,\\
        \begin{bmatrix}
            \mu^2 & 0 & \mu^3
        \end{bmatrix}^T &\text{ if } x> 0.5.\\
    \end{cases}
\end{equation*}
In Fig.~\ref{fig:motivational_sod_more_params}, we show some analytical solutions for the parametric 1D Sod shock tube, for different values of $\mu$.
Looking at Fig.~\ref{fig:motivational_sod_more_params} it is clear there is no straightforward choice for the reference solution $\rhoref$ that would lead to reasonable minimization problems, similar to the one presented in Algorithm~\ref{alg:calibration self similar}. Indeed, we would minimize the $L^2$ error between two solutions ($\rhoref$ and $\rhohat(\bbmu)$) that may have very different heights at the boundaries. When different boundaries behaviors are present, the calibration procedure needs to be reformulated in a different way.

{Instead of fixing a reference solution $\rhoref$, we consider a \textit{suitable} linear reduced space $V_{\POD}$, thus introducing a new residual for the minimization problem:}
\begin{equation}
\label{eq:optimization many parameters}
\begin{split}
\min_{\bm{\theta}({\bbmu})\in \mathbb R^{Q}}  & \lVert \rhohat ({\bbmu}; {\bbw(\bbmu})) - \Pi_{V_{\POD}}\rhohat ({\bbmu};{\bbw(\bbmu}))  \rVert_{L^2(\Omegaref)} + \frac{\delta}{2}\lVert\partial_{\bbmu}\bbw (\bbmu) \rVert^2_{\ell^2(\mathcal{P})}+\\
&
\frac{\alpha}{2} \max_{\bm{x}\in\Omega} \left( \max \left\lbrace \left\lVert\nabla T[\bbw(\bbmu)](\hat{\bm{x}})\right\rVert,  \left\lVert\nabla T^{-1}[\bbw(\bbmu)](\bm{x})\right\rVert \right\rbrace \right),
\end{split}
\end{equation}
{where we recall that $\bbtheta(\bbmu)=\bbtheta(\bbw(\bbmu))\in\mathbb{R}^Q$ is the vector of the free coordinates of the control points, and $Q\leq d\times M$ is the total number of free coordinates in the calibration step.}
In Problem~\eqref{eq:optimization many parameters}, $\bbtheta$ is constrained as in the previous section and $\Pi_{V_{\POD}}$ is the orthogonal projection onto a linear space $V_{\POD}$ obtained through a preliminary procedure that we will describe now shortly.

This minimization allows to overcome the issue of quasi--self--similar solutions thanks to the projection onto a reduced space. {To carry out this projection, we first need to build \textit{a priori} a suitable linear space $V_{\POD}$, which needs to capture the \textit{minimal} amount of information on the solution manifold: we want to be able to approximate the solution manifold with a linear space of very low dimension.}

We therefore introduce the matrix of the free coordinates $\thetavec\in\mathbb{R}^{N_{\text{few}}\times Q}$:
\begin{equation*}
   \thetavec[i, :] = \bbtheta(\tildebbmu_i) , \quad i=1, \dots, N_{\text{few}},
\end{equation*}
{where we recall once again that $\bbtheta(\tildebbmu_i)$ is the vector of the free coordinates of the control points $\bbw(\tildebbmu_i)$ associated to the parameter $\tildebbmu_i$.}
The free coordinates can be selected through another optimization process carried out in the same spirit of \eqref{eq:optimization many parameters} on the whole matrix {$\thetavec$, minimizing the projection error over all the parameters $\tildebbmu_1, \dots, \tildebbmu_{N_{\text{few}}}$, while updating the space $V_{\POD}(\thetavec)$ obtained by compressing the solutions on the reference domain for these parameters transformed using the calibration points given by $\thetavec$}.  
Therefore, we solve the following constrained minimization problem:
\begin{equation}
\label{eq:optimal calibration few parameters}
\min_{\thetavec\in \mathbb R^{N_{\text{few}}\times Q}}  \sum_{i=1}^{N_{\text{few}}}\lVert \rhohat (\tildebbmu_i; {\bbw(\tildebbmu_i)}) - \Pi_{V_{\POD}(\thetavec)}\rhohat (\tildebbmu_i; {\bbw(\tildebbmu_i)})  \rVert_{L^2(\Omegaref)}, 
\end{equation}
with $V_{\POD}:=\POD(\lbrace \rhohat(\tildebbmu_i) \rbrace_{i=1}^{N_{\text{few}}}, N^{\POD}_{\text{few}})$.
This optimization process is of larger dimension with respect to the previous ones and it requires, for each residual evaluation, the computation of a POD over the $N_{\text{few}}$ calibrated snapshots $\rhohat(\tildebbmu_i; {\bbw(\tildebbmu_i)})$, with a user defined number of modes $N_{\text{few}}^{\POD}$. {The choice of $N_{\text{few}}^{\POD}$ is, in general, not simple. In our numerical tests, we adopted the following heuristics: we want to take into account all the different travelling features. For this reason, in the Sod test case, we have chosen $N_{\text{few}}^{\POD}=3$ to keep into account the different values that there might be between the rarefaction and the contact discontinuity and between the contact and the shock. 
$N_{\text{few}}^{\POD}$ should be kept as low as possible not to overload the minimization process. In the Sod test case, we performed a few tests with increasing $N_{\text{few}}^{\POD}$ and $3$ was the first value where the minimization process was giving successful results. While this is an heuristic argument, it still provides a valid starting point for the numerical simulations. A more in depth analysis on the role of $N_{\text{few}}^{\POD}$ and on how to choose it in an automated way is envisioned as a future research direction.}
Again, we solve problem~\eqref{eq:optimal calibration few parameters} with SLSQP using the same constraints defined in the previous section. This extra optimization step can be skipped when other techniques to detect interesting features can be used \cite{taddei_jcp}. One possibility is the use of classical shock detection procedures to find the steepest points of the solutions and calibrating them \cite{mirhoseini2023model}.   We summarize the steps of the whole procedure in Algorithm~\eqref{alg:calibration quasi self similar}.

\begin{algorithm}
    \caption{Calibration for quasi--self--similar solution manifold} \label{alg:calibration quasi self similar}
    \begin{algorithmic}
    \State \textbf{Input:}  the training solution manifold on the physical domain $\lbrace \rho(\mu_\ell, t_j)\text{ for } \ell=1,\dots, N_{\text{train}}, j = 1, \dots,N_{\mu_\ell} \rbrace$ and the number of modes $N^{\POD}_{\text{few}}$.
    \State Select \textit{few} parameters $\tildebbmu_1, \dots, \tildebbmu_{N_{\text{few}}}$.
    \State Find optimal $\thetavec \in \mathbb R^{N_{\text{few}}\times Q}$ by \eqref{eq:optimal calibration few parameters} for the \textit{few} parameters $\tildebbmu_1, \dots, \tildebbmu_{N_{\text{few}}}$.
    \For{$\ell = N_{\text{train}} , \dots, 1$ }
        \For{$j=N_{\mu_{\ell}},\dots, 1$}
            \State Set the initial guess $\bbtheta^{(0)}(\bbmu_{( \ell,j)})$ as in \eqref{eq:initial guess parametric setting}.
            \State Solve Problem \eqref{eq:calibration self similar}.
        \EndFor
    \EndFor
    \State \textbf{Output:} the optimal control points $\bbw^{\text{opt}}(\bbmu_{(\ell,j)})$.
    \end{algorithmic} 
\end{algorithm}
{As an alternative to the whole projection onto $V_{\text{POD}}$ we would also like to remark that, whenever shock detection techniques can be used, the quasi-self-similar setting can be efficiently handled by the convex displacement interpolation (CDI) \cite{taddei_iollo_cucchiara_telib}}.

\section{The geometrical transformation map}\label{sec:geom_trans}
We now present the geometrical transformation map $T[\bbw^{\text{opt}}(\bbmu)]\colon\mathcal{R}\mapsto\Omega$ used to define the calibrated snapshots. We start with the simpler case, namely the 1D setting, and later we consider the 2D case.
\subsection{The 1D setting}\label{sec:geometrical transformation 1D}
Let $\bbw^{\text{opt}}(\bbmu)$ be the control points whose free coordinates are the solution to Problem \eqref{eq:calibration self similar}: to interpolate the values $\{(\overline{\bbw}_{\bbalpha}(\bbmu), \bbw^{\text{opt}}_{\bbalpha}(\bbmu))\}_{\bbalpha}$, we use monotone cubic $C^1$ splines. These are the so-called PCHIPs (Piecewise Cubic Hermite Interpolating Polynomials) interpolators, available in the \texttt{scipy} Python library under the \texttt{interpolate} classes and the built-in function is called \texttt{PchipInterpolator}. 
By employing monotone cubic splines, we obtain a transformation function that preserves the monotonicity of the calibration map $\bbw(\cdot)$, guaranteeing its bijectivity and $C^1$ smoothness if the calibration points are in the ``right order'', as prescribed in Section~\ref{sec:calibration}.
\begin{figure}
    \centering
    \includegraphics[width=\textwidth]{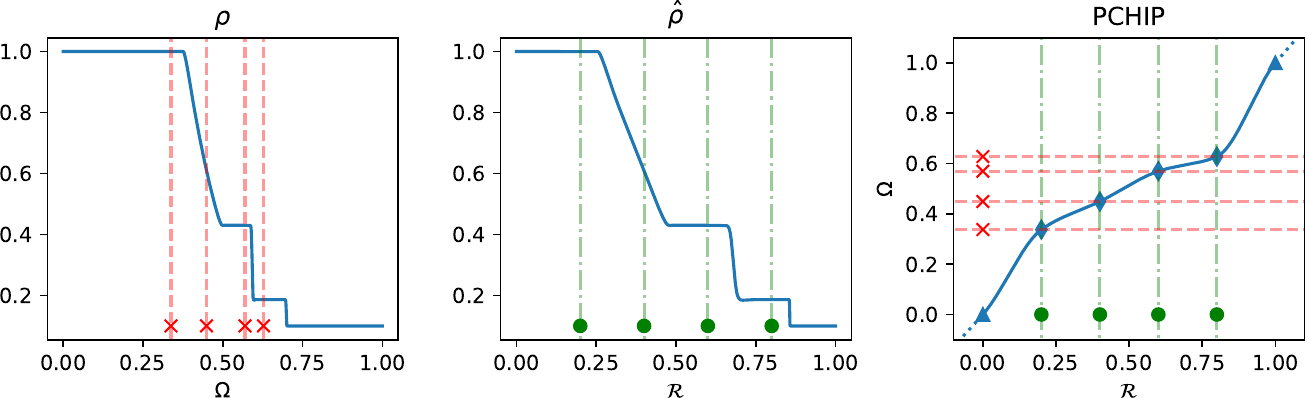}
    \vspace{-3mm}
\caption{Example of PCHIP calibration: solution on physical domain with calibration control points (left), solution on reference domain with reference control points (center) and PCHIP transformation with all control points (right)}
    \label{fig:pchip}
\end{figure}
In Fig.~\ref{fig:pchip}, we can see an example of a PCHIP transformation applied to one of the snapshot calibrated on the detected features. 
On the right of the figure, the transformation is depicted and we can observe that it interpolates the points, it is $\mathcal C^1$ and it is very close to the identity on the boundaries, because we introduce two extra interpolation points outside the boundaries of the domain. 
This helps to keep the regularity of the transformation in the ALE formulation \cite{taddei_sisc}.

In this work, we do not focus on the ALE formulation on $\Omegaref$, nevertheless the PCHIPs allow to easily compute all the necessary ingredients. Indeed, they are polynomials and their derivatives and the inverse of their derivative is easy to compute. 
Moreover, the inverse of the transformation exists and it is unique in each point, hence, with a simple Newton method, we can easily recast the inverse function.

\subsection{The 2D setting}\label{sec:geometrical transformation 2D}
In the 2D setting we use tensor product of one dimensional PCHIPs, in order to exploit their properties for Cartesian geometries.
We refer again to Fig.~\ref{fig:ref_grid} to better understand the transformation map. We need to compute $T[\bbw^{\text{opt}}(\bbmu)]\colon\Omegaref\mapsto\Omega$, such that 
\begin{equation*}
T[\bbw^{\text{opt}}(\bbmu)](\overline{\bbw}^1_{\alpha_1,\alpha_2}, \overline{\bbw}^2_{\alpha_1,\alpha_2})=(\bbw^{1, \text{opt}}_{\alpha_1,\alpha_2}(\bbmu), \bbw^{2, \text{opt}}_{\alpha_1,\alpha_2}(\bbmu))\quad \text{ for } \alpha_i=1,\dots, M_i, \, i=1,2.   
\end{equation*}
Let $\bbxhat\in\Omegaref$ be a point with coordinates $\bbxhat=(\hat{x}, \hat{y})$. We define the map:
\begin{align*}
T[\bbw^{\text{opt}}(\bbmu)](\xhat, \yhat)& := (T^x[\bbw^{\text{opt}}(\bbmu)](\xhat, \yhat), T^y[\bbw^{\text{opt}}(\bbmu)](\xhat, \yhat)), \text{ with }\\
    T^x[\bbw^{\text{opt}}(\bbmu)](\xhat, \yhat) &:= \sum_{\ell=1}^{M_2}{\gamma}^y_{\ell}(\yhat)P^x_{\ell}(\xhat), \qquad T^y[\bbw^{\text{opt}}(\bbmu)](\xhat, \yhat) := \sum_{k=1}^{M_1}{\gamma}^x_k(\xhat)P^y_k(\yhat).
\end{align*}
In the previous equations, we made use of the following quantities:
\begin{enumerate}
    \item $P^x_{\ell}$ is a PCHIP interpolating the points  $\lbrace \overline{\bbw}^1_{\alpha_1, \ell}, \bbw^{1,\text{opt}}_{\alpha_1,\ell}(\bbmu) \rbrace_{\alpha_1=1}^{M_1}$, where the control points $\overline{\bbw}^1_{\alpha_1, \ell} $ for $\alpha_1=1,\dots,M_1$ are on horizontal lines $\hat{y}= \bar{y}_\ell$ in the reference domain, see Fig.~\ref{fig:ref_grid}, namely for all $\alpha_1=1,\dots,M_1$ we have $P^x_{\ell}(\overline{\bbw}^1_{\alpha_1,\ell}) = \bbw^{1,\text{opt}}_{\alpha_1,\ell}(\bbmu);$
    \item $P^y_k$ is a PCHIP interpolating the points  $\lbrace \overline{\bbw}^2_{k,\alpha_2}, \bbw^{2,\text{opt}}_{k,\alpha_2}(\bbmu) \rbrace_{\alpha_2=1}^{M_2}$, where the control points $ \overline{\bbw}^2_{k,\alpha_2} $ for $\alpha_1=1,\dots,M_1$ are on vertical lines $\hat{x}= \bar{x}_k$ in the reference domain, see Fig.~\ref{fig:ref_grid}, namely for all $\alpha_2=1,\dots,M_2$ we have $P^y_k(\overline{\bbw}^2_{k,\alpha_2}(\bbmu)) = \bbw^{2,\text{opt}}_{k,\alpha_2}(\bbmu);$
    \item ${\gamma}^y_{\ell}(\cdot)$ is a PCHIP interpolating the points $\{\overline{\bbw}^2_{\alpha_1,\alpha_2}, \delta_{\alpha_2,\ell}\}_{\alpha_2=1}^{M_2}$, being $\delta_{\alpha_2,\ell}$ the Kronecker delta.
    By doing so, we obtain that $T^x[\bbw^{\text{opt}}(\bbmu)]$ is a convex combination of the $\lbrace P^x_{\ell}(\cdot) \rbrace_{\ell=1}^{M_2}$ such that $T^x[\bbw^{\text{opt}}(\bbmu)](\xhat, \yhat=\overline{\bbw}^2_{\alpha_1, \alpha_2}) = P^x_{\alpha_2}(\xhat)$;
    \item ${\gamma}^x_k(\cdot)$ is a PCHIP interpolating the points $\{\overline{\bbw}^1_{\alpha_1,\alpha_2}, \delta_{\alpha_1,k}\}_{\alpha_1=1}^{M_1}$, as before, leading to the property $T^y[\bbw^{\text{opt}}(\bbmu)](\xhat=\overline{\bbw}^1_{\alpha_1,\alpha_2}, \yhat) = P^y_{\alpha_1}(\yhat)$.
\end{enumerate}
Notice that, ultimately, it holds that 
$$T[\bbw^{\text{opt}}(\bbmu)](\overline{\bbw}^1_{\alpha_1,\alpha_2}, \overline{\bbw}^2_{\alpha_1,\alpha_2})= (P_{\alpha_2}^x(\overline{\bbw}^1_{\alpha_1,\alpha_2}),P_{\alpha_1}^y(\overline{\bbw}^2_{\alpha_1,\alpha_2})) =  \bbw^{\text{opt}}_{\alpha_1,\alpha_2}(\bbmu).$$ 

Also in the 2D case, the Jacobian of the transformation, which is needed in the ALE formulation, is easily accessible, since all the terms are polynomials. Similarly to the 1D case, we can compute the Jacobian of the inverse of the transformation, using the inverse of the Jacobian of the transformation, provided that we can invert the map $T$.

\begin{figure}
    \centering
    \includegraphics[width=\textwidth]{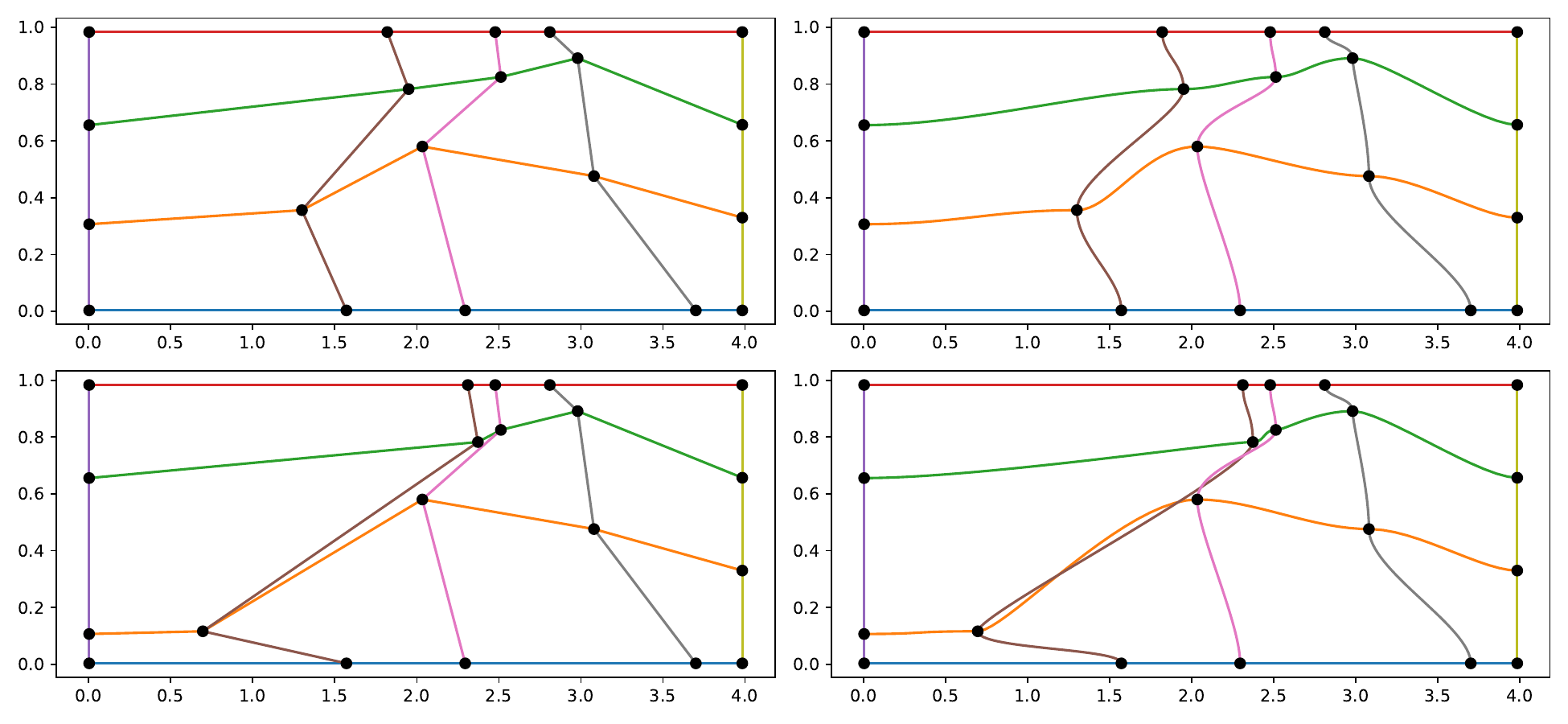}
    \vspace{-5mm}
    \caption{Calibration points ($M_1=5, M_2=4$) on the physical domain $\Omega=[0,4]\times [0,1]$ (left) and the transformation of horizontal and vertical lines on the physical domain (right). On the top calibration points that lead to an invertible map, on the bottom calibration points that respect the monotonicity in each line (left), but whose transformed horizontal and vertical lines cross multiple times (right) and, hence, whose transformation is not invertible}
    \label{fig:crossing_lines}
\end{figure}
\subsubsection{Invertibilty of $T$}
    The computation of the inverse of $T$ is fundamental in many aspect of the algorithm: to display the function on the physical domain, to compute quantities and errors on the physical domain and to compute the inverse of the Jacobian for the ALE formulation.

    Unfortunately, the so defined map $T$ cannot be proven to be invertible, as it might happen that some of the vertical and horizontal lines cross each other multiple times in the physical domain, see Fig.~\ref{fig:crossing_lines} second line. To avoid this, we impose, in the optimization procedure for the calibration, to have positive determinant of the Jacobian of the transformation on the meshpoints of the reference domain, see the constraints in Section~\ref{sec:calibration self similar}. This typically guarantees invertibility. We recall that, being PCHIPs polynomials, the computation of the Jacobian can be performed explicitly in each point. To compute the inverse of the transformation, we perform the following steps. We first apply the transformation map $T$ to the elements of the Cartesian meshgrid of $\Omegaref$: these will be mapped to quadrilateral elements in $\Omega$. Now, given a point $\bbx\in\Omega$, we can easily find to which of these quadrilaterals it belongs and with a Projective transformation (see the python module \texttt{transform} of the  package \texttt{skimage}) we pull it back onto $\Omegaref$. 
    Then, we use the found point as initial guess to find the solution of $T[\bbw^{\text{opt}}(\bbmu)](\hat{\bbx}) = \bbx$ through a Newton type nonlinear solver (\texttt{scipy.root}).

\section{Model Order Reduction with calibration}\label{sec:MOR with calibration}
We are now ready to perform the model order reduction step. In what follows the procedure is similar for the non parametric and the parametric setting: we will therefore present it for the latter case, for the sake of generality.
\subsection{Learning the calibration map}\label{sec:ANN}
Once the calibration procedure presented in Section \ref{sec:calibration} has been carried out, the map $\bbw\colon\mathcal{P}\mapsto\Theta$ is known only through the sample values $\bbw^{\text{opt}}(\bbmu)$, for $\bbmu\in\mathcal{P}^{\text{train}}$. To learn the calibration map $\bbw(\cdot)$ for any parameter $\bbmu$, we employ an Artificial Neural Network (ANN) composed by several layers: an input layer $\bbmu\in \mathcal{Y}^{(0)}=\mathcal{P}$ where we pass the parameters of the problem of interest, L hidden layers $\mathcal{Y}^{(j)}$, $j=1,\dots,L$, and an output layer $\mathcal{Y}^{(L+1)}$, see Fig.~\ref{fig:NN}. As output layer, we would like to obtain $\bbtheta \in \mathbb R^Q$, from which we can extract the calibration points $\bbw$. Keeping in mind the monotonicity constraints applied to the control points in the constrained minimization problems \eqref{eq:calibration self similar} and \eqref{eq:optimization many parameters}, we try to enforce this constraint in the ANN. It is not easy to strongly enforce such constraints, but we can force the output to be positive, using as final activation function a Softplus. Hence, we take as output layer of our ANN not directly $\bbtheta$, but the vector $\bbv \in \mathbb R^Q$ of the differences of the free coordinates with the \textit{previous} ones (in 2D it is referred to the same line). Doing so, the positivity of $\bbv$ is equivalent to the monotonicity in each line of points. In one dimension, it is defined as $v_i = w_i-w_{i-1}$ for $i=2,\dots,M_1-1$, while in two dimensions this operation is done in each horizontal or vertical line.
\tikzset{%
  every neuron/.style={
    circle,
    draw,
    minimum size=0.65cm
  },
  neuron missing/.style={
    draw=none, 
    scale=1,
    text height=0.333cm,
    execute at begin node=\color{black}$\vdots$
  },
}
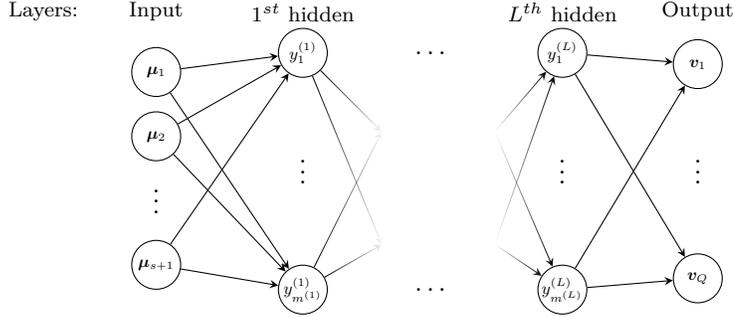
\begin{figure}
\centering
\begin{tikzpicture}[x=1.5cm, y=1.5cm, >=stealth]


\node [align=center, above] at (1.0, 1.6) {\footnotesize Layers:};
\node [align=center, above] at (2.0, 1.6) {\footnotesize Input};
\node [align=center, above] at (3.3, 1.6) {\footnotesize $1^{st}$ hidden};
\node [align=center, above] at (5.6, 1.6) {\footnotesize $L^{th}$ hidden};
\node [align=center, above] at (6.8, 1.6) {\footnotesize Output};

\foreach \m/\l [count=\y] in {1,2,missing,4} 
  \node [every neuron/.try, neuron \m/.try] (input-\m) at (2, 1.8-\y*0.57) {};

\foreach \m [count=\y] in {1,missing,2}
  \node [every neuron/.try, neuron \m/.try ] (hidden1-\m) at (3.3, 2.45-\y*1.05) {};

\foreach \m [count=\y] in {1,missing,2}
  \node [every neuron/.try, neuron \m/.try ] (hidden2-\m) at (5.6, 2.45-\y*1.05) {};

\foreach \m [count=\y] in {1,missing,2}
  \node [every neuron/.try, neuron \m/.try ] (output-\m) at (6.8, 2.25-\y*0.95) {};

\foreach \l [count=\i] in {1,2,4}{
    \ifnum \l=4
        \node[scale=0.7] at (input-\l){\small $\bbmu_{s+1}$};
    \else
        \node[scale=0.7] at (input-\l){\small $\bbmu_\l$};
    \fi
    }
\foreach \l [count=\i] in {1,2} {
    \ifnum \l=1
        \node[scale=0.7] at (hidden1-\l) {\small $y^{(1)}_{\i}$};
    \else 
        \node[scale=0.7] at (hidden1-\l) {\small $y^{(1)}_{m^{(1)}}$};
    \fi
    }
\foreach \l [count=\i] in {1,2} {
    \ifnum \l=1
        \node[scale=0.7] at (hidden2-\l) {\small $y^{(L)}_{\i}$};
    \else 
        \node[scale=0.7] at (hidden2-\l) {\small $y^{(L)}_{m^{(L)}}$};
    \fi
    }
\foreach \l [count=\i] in {1,2}{
        \ifnum \l=2
        \node[scale=0.7] at (output-\l){\small $\bbv_Q$};
    \else
        \node[scale=0.7] at (output-\l){\small $\bbv_\l$};
    \fi
    }
\foreach \i in {1,2,4}
  \foreach \j in {1,...,2}
    \draw [->] (input-\i) -- (hidden1-\j);

\foreach \i in {1,...,2}
    \foreach \j in {1,...,2}
        \draw [->] (hidden2-\i) -- (output-\j);

\node[fill=white,scale=1,inner xsep=0pt,inner ysep=5mm] at ($(hidden1-1)!.5!(hidden2-1)$) {$\dots$};
\node[fill=white,scale=1,inner xsep=0pt,inner ysep=5mm] at ($(hidden1-2)!.5!(hidden2-2)$) {$\dots$};

\foreach \i in {1,...,2}
    \draw [->, scope fading=east] (hidden1-\i) -- (4, 0.7);

\foreach \i in {1,...,2}
    \draw [->, scope fading=east] (hidden1-\i) -- (4, -0.3);

\foreach \i in {1,...,2}
    \draw [->, scope fading=west] (5, 0.7) -- (hidden2-\i);

\foreach \i in {1,...,2}
    \draw [->, scope fading=west] (5, -0.3) -- (hidden2-\i);

\end{tikzpicture}
\vspace{-2mm}
\caption{Example of the architecture of an ANN}\label{fig:NN}
\end{figure}

Each layer $\mathcal{Y}^{(j)}$, $j=1,\dots,L$ is connected to the next and to the previous ones through affine maps $\delta^{(j)}\colon\mathcal{Y}^{(j)}\mapsto\mathcal{Y}^{(j+1)}$ and at every node a nonlinear activation function $\zeta^{(j)}\colon\mathbb{R}\mapsto\mathbb{R}$ is applied component–wise. We used the hyperbolic tangent in all the ANN except in the output layer where the Softplus function is used for the positivity of the outputs. 
On a training set, the learning process changes the weights $\delta^{(k)}$ minimizing the error between the output and the optimal calibration points.

To build and train the ANN, we used a Python based library \texttt{EZyRB} \cite{ezyrb} that uses Adam method \cite{kingma2014adam} to perform the minimization.
More details on the architecture of the ANN for the different test cases will be provided in the numerical section.

\subsection{POD-NN}\label{sec:POD-NN}
In Section~\ref{sec:calibration self similar} and \ref{sec:calibration quasi self similar}, we presented the calibration algorithm in the self-similar and in the quasi self-similar setting. 
Thanks to the calibration of the snapshots, we obtain the \emph{calibrated manifold} $\hat{\mathcal{M}}_\rho = \{\rhohat(\bbmu)(\cdot), \text{ }\bbmu\in\mathcal{P}_{\text{train}}\}$, where $\rhohat(\bbmu)(\cdot)$ is the \emph{calibrated snapshot} defined in Section \ref{sec:calibration}.
We now proceed with the compression of $\hat{\mathcal{M}}_\rho$ by means of the Proper Orthogonal Decomposition (POD): we refer to \cite{stamm_rozza_hesthaven} for more details.
The dimension of the reduced basis can be chosen either setting a maximum number of basis $N^{\text{max}}_{\text{POD}}$ or choosing the most relevant modes such that the discarded energy is smaller than a certain tolerance $\tau_{\text{POD}}$. 
Once the POD has been carried out, we obtain a linear space $V^{\rho}_N$ spanned by the $N$ orthonormal reduced basis functions $\{\phiref_i\}_{i=1}^N$ on the reference domain $\Omegaref$. $V^{\rho}_N$ should now represent with a good accuracy any element of the calibrated solution manifold $\hat{\mathcal{M}}_{\rho}$.
\subsection{Online solution by means of ANN}\label{sec:ROM-NN} 
In this work, we mainly focus on the calibration procedure and the offline phase. 
Hence, we will use a non--intrusive approach for the reconstruction of the online solutions.
Let $\bbmu\in\mathcal{P}$ be a parameter of choice: the goal is to construct a linear approximation $\rho_N(\bbmu)$ of the snapshot $\rho(\bbmu)$. It should be clear by now that, since we are dealing with advection dominated problems, this is not a simple task within the standard MOR setting. However, in Sections~\ref{sec:calibration} and \ref{sec:geom_trans}, we presented a calibration technique that allows us to obtain a linear space $V^{\rho}_N$ that approximates with good accuracy the calibrated manifold $\hat{\mathcal{M}}_\rho$: we are therefore able to construct a linear approximation of the calibrated solution of interest $\rhohat(\bbmu)$ \emph{in the reference domain} $\Omegaref$. This means that we can approximate $\rhohat(\bbmu)$ with
\begin{equation*}
    \rhohat(\bbmu)\approx \rhohat_N(\bbmu):=\sum_{i=1}^N\underline{\rhohat}_N^i(\bbmu)\phiref_i,
\end{equation*}
where $\underline{\rhohat}_N^i(\bbmu) = (\rhohat(\bbmu), \phiref_i)_{L^2(\Omegaref)}$ for $i=1,\dots, N$, being $(\cdot,\cdot)_{L^2(\Omegaref)}$ the $L^2$ scalar product on the reference domain.
In order to find the vector $\underline{\rhohat}_N(\bbmu)$ one can adopt two alternative ways: an intrusive approach, by means of a Galerkin projection of the high order algebraic system onto the reduced space, or a non-intrusive approach by means of an ANN. If the standard Galerkin projection setting is adopted, the online system and the reconstruction of the online solution is carried out within an ALE formalism \cite{richter, torlo,taddei_sisc}: the original problem of interest, formulated over $\Omega$, has to be re-written into a problem formulated over the reference domain $\Omegaref$. 
In this approach, a hyper-reduction procedure \cite{barrault2004empirical,yano2019discontinuous} will be necessary to tackle the nonlinearities of the problem and of the transformation map. This approach is currently under investigation, and it will be part of a future extensions.

An alternative to the intrusive approach is represented by the use of ANN, in the spirit of Section~\ref{sec:ANN}. 
We consider the map $\Pi_N\colon\mathcal{P} \mapsto\mathbb{R}^N$
\begin{equation*}
\begin{split}
    \Pi_N(\bbmu) := \underline{\rhohat}_N(\bbmu) = [\underline{\rhohat}_N^1(\bbmu),\dots,\underline{\rhohat}_N^N(\bbmu)]^T.
\end{split}
\end{equation*}
We make use of an ANN to learn the $L^2$ projection map $\Pi_N$: to train the map, we employ the set of input samples $\bbmu$ for $\bbmu\in\mathcal{P}^{\text{train}}$, and the set of output samples $\{\underline{\rhohat}_N(\bbmu)\text{ for }\bbmu\in\mathcal{P}^{\text{train}}\}$. 
In this algorithm, we do not use the optimal calibrated points obtained with the optimization process, but we use the ones predicted by the ANN of Section~\ref{sec:ANN}. Doing so, any systematic error in the online calibration should be already taken into account while performing the $L^2$ projection and automatically corrected by this approach. Moreover, it is also possible to use different training sets for the calibration optimization procedure and the model order reduction ones. 
More details on the architecture of the ANN employed will be provided in the numerical section.

\subsection{Offline--online splitting}
{We summarize in this section the operations that should be performed in the offline part (i.e. only once for every new parametric problem one wishes to consider) and the online operations that should be performed to predict the solution for a new parameter.\\
\textbf{Offline phase}: FOM simulations for all parameters in the training set, choice of the reference calibration points, optimization procedure to find the calibration points for a training set, learning of the regression map for the calibration points for all parameters, transformation of the FOM solutions onto the reference domain, computation of the RB, learning of the regression for the coefficients-to-(reduced coefficient) map.\\
\textbf{Online phase}: given a new parameter evaluate the regression of the calibration points to know the geometrical transformation, evaluate the regression for the RB coefficients to get the solution on the reference domain, combine the two maps to get the solution on the physical domain.
}

\section{Numerical results}\label{sec:numerical}
We now present some numerical results to validate the proposed methodology. We will consider different time dependent test cases: the Sod shock tube problem in 1D, in the non parametric and in the parametric setting, already introduced in Section \ref{sec:calibration self similar} and \ref{sec:calibration quasi self similar}, respectively. To further test the performance of our methodology, we subsequently consider a 2D problem, namely the double Mach reflection (DMR) problem, again in the non parametric and in the parametric setting. {To conclude, some numerical results for the non parametric triple point problem will be shown.}
For all presented tests, the FOM consists of high order solutions obtained with an explicit Finite Volume discretization with WENO reconstruction of order 5, with explicit time discretization given by the optimal SSPRK54, with CFL coefficient 0.8 and Rusanov numerical flux.
\subsection{Non-parametric Sod shock tube problem in 1D}\label{sec:numerical results Sod non-parametric}
We consider Problem~\eqref{eq:euler conservative} introduced in Section~\ref{sec:motivation}, where the physical domain is $\Omega=\Omegaref=[0,1]$.  The number of spatial degrees of freedom is 1500 and this leads to computational costs of around 2 minutes using a Fortran code~\cite{Nunezrepo} on a Intel(R) Xeon(R) CPU E3-1245 v5 @ 3.50GHz. In both cases, we use Dirichlet BC as the waves do not exit the domain before the final time. The details of all the relevant quantities are presented in Table~\ref{tab:Sod non-parametric}.
\begin{table}
\caption{Test details for the 1D Sod shock tube problem, non parametric (Nonparam) and parametric (Param) setting}
\label{tab:Sod non-parametric}
\centering
\begin{tabular}{|c| c| c||c |c| c|} 
\hline
 Quantity & Nonparam & Param & Quantity & Nonparam & Param  \\[0.5ex]
 \hline
 $\rho_L$ & $1$ & $[0.7,1.3]$ & $\rho_R$ & $0.1$ &$[0.1,0.15]$ \\
 $u_L$ & $0$ & 0& $u_R$&0 & $0$ \\
 $p_L$ & $1$ &$[0.7, 1.3]$& $p_R$& $0.125$  & $[0.05, 0.15]$\\ \hline
$t_f$ & $0.2$s& $0.2$s & $N_h$ & 1500&1500\\ 
 \hline
\end{tabular}
\end{table}
\begin{figure}
    \centering
    \includegraphics[width=\textwidth, trim={0 10 0 10},clip]{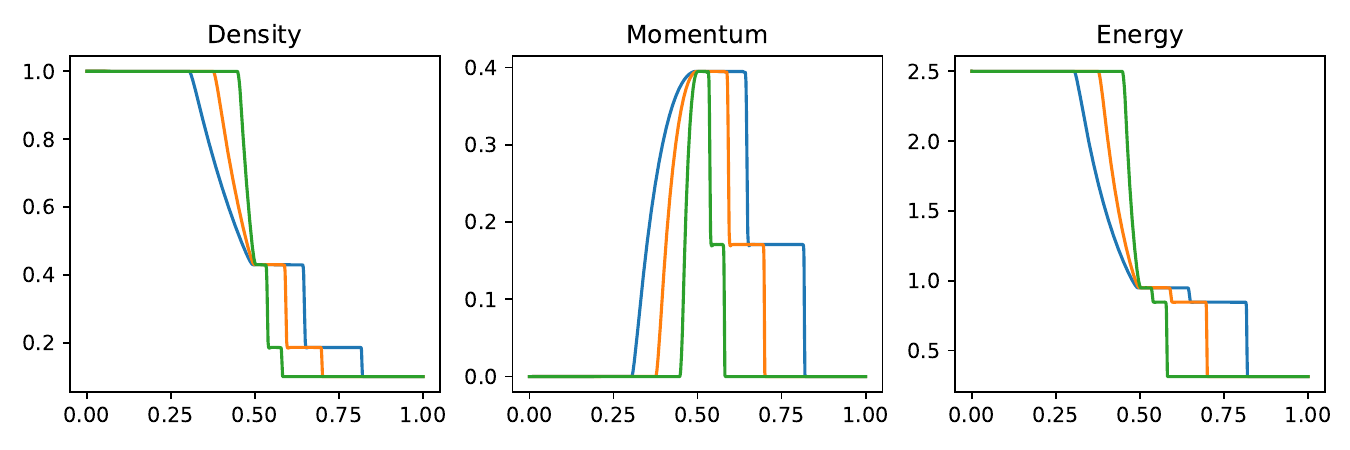}
    \includegraphics[width=\textwidth, trim={0 10 0 10},clip]{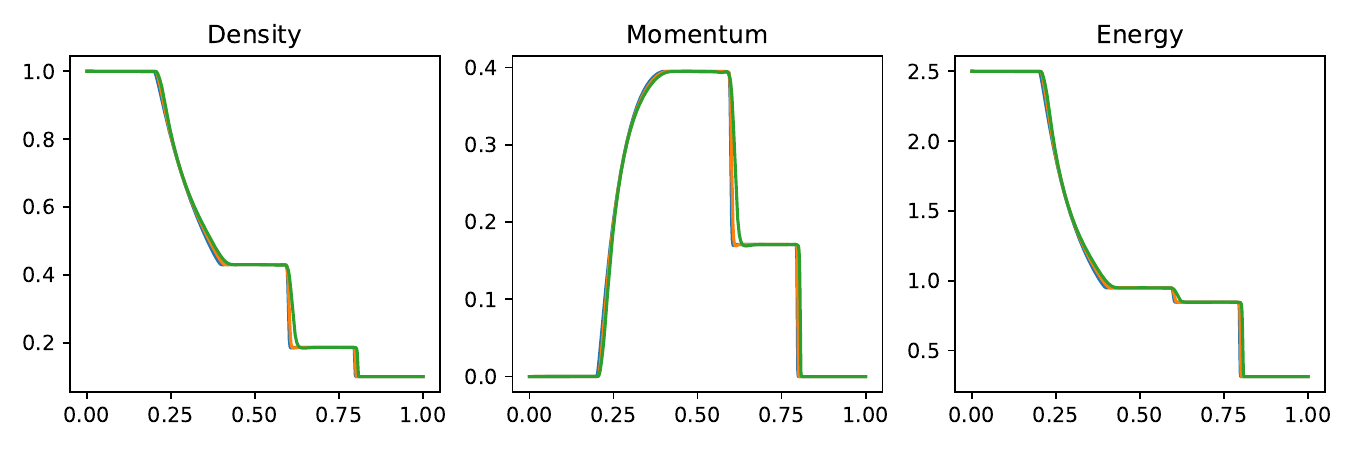}
    \vspace{-3mm}
    \caption{FOM simulation of Sod 1D problem non parameteric at times 0.04 (green), 0.1 (orange) and 0.16 (blue), on the original domain (top) and calibrated on the reference domain (bottom)}
    \label{fig:Sod_FOM}
\end{figure}

Fig.~\ref{fig:Sod_FOM} shows some snapshots for the density $\rho$, the moment $m$ and the energy $E$ at three different times of the simulation: it is clear that the structures of all three the components of the solution present discontinuities that travel in the domain at the same locations.

We carry out the calibration technique proposed in Algorithm~\ref{alg:calibration self similar}: we choose as reference solution $\rhoref$ the density $\rho(t=0.16)$. We then choose $M=4$ control points $\bbwref_1=0.2,\,\bbwref_2=0.4,\,\bbwref_3=0.6,\,\bbwref_4=0.8$ equispaced in the reference domain $\Omegaref=[0,1]$. 
The calibrated solutions (using the ANN to forecast the calibration points) are shown in Fig.~\ref{fig:Sod_FOM}: the main features of the solutions, namely the shock, the contact discontinuity and the rarefaction wave are well aligned with the reference solution. 
The details of the quantities required to carry out the calibration step are listed in Table~\ref{tab: calibration Sod non-parametric}. 
We point out that the calibration step and the training of the ANN have been carried out on the training set $\mathcal{P}^{\text{train}}\subset[0.01, 0.16]\subset [0, 0.2]$ sampled with 25 equispaced parameters. 
The reason for excluding the first timesteps from the training is that the minimization is tricky during the first timesteps: indeed we have a transition phase, during which all the features are in the same point, leading to non invertible maps. An alternative way to overcome this difficulty could be to restore to local reduced basis spaces \cite{amsallem2012nonlinear}, or to use a FOM approach for the first timesteps.
The final times are excluded to test the extrapolatory performances of the ROM.
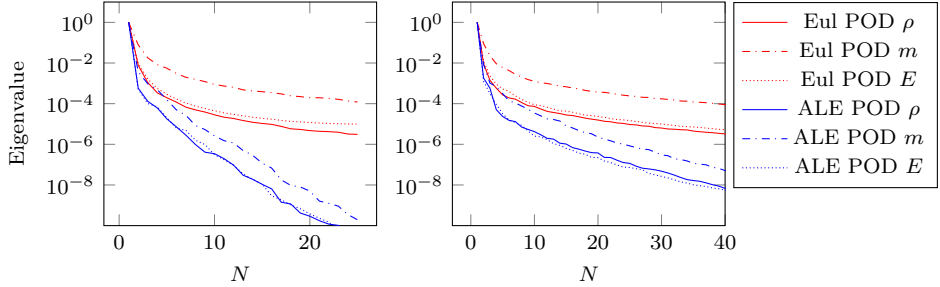
\begin{figure}
    \centering
   \begin{tikzpicture}
        \begin{axis}[ymode=log,
                    xlabel={$N$},
                    xmax=27,
                    ymin=1e-10,
                    ylabel={Eigenvalue},
			      ytick={1,0.01,0.0001,0.000001, 1e-8},
                    width=0.4\textwidth,
                    height=0.35\textwidth,
                    style={font=\footnotesize}]  
            \addplot[mark size=1.3pt, red] table [mark=triangle*,x=Number, y=Eul_eig_norm, col sep=comma] {Dat001.csv};
            \addplot[dashdotted, mark size=1.3pt, red] table [x=Number, y=Eul_eig_norm, col sep=comma] {Dat002.csv};
            \addplot[densely dotted,mark size=1.3pt, red] table [x=Number, y=Eul_eig_norm, col sep=comma] {Dat003.csv};
            
            \addplot[mark size=3pt, blue] table [mark=square*,x=Number, y=ALE_eig_norm, col sep=comma] {Dat004.csv};
            \addplot[dashdotted, mark size=3pt, blue] table [x=Number, y=ALE_eig_norm, col sep=comma] {Dat005.csv};
            \addplot[densely dotted, mark size=3pt, blue] table [x=Number, y=ALE_eig_norm, col sep=comma] {Dat006.csv};
        \end{axis}
    \end{tikzpicture}
   \begin{tikzpicture}  
        \begin{axis}[ymode=log,
                    xlabel={$N$},
                    xmax=40,
                    ymin=1e-10,
                    ytick={1,0.01,0.0001,0.000001, 1e-8},
                    legend pos=outer north east,
                    width=0.4\textwidth,
                    height=0.35\textwidth,
                    style={font=\footnotesize}]
            
            \addplot[mark size=1.3pt, red] table [mark=triangle*,x=Number, y=Eul_eig_norm, col sep=comma] {Dat007.csv};
            \addlegendentry{Eul POD $\rho$}
            \addplot[dashdotted, mark size=1.3pt, red] table [x=Number, y=Eul_eig_norm, col sep=comma] {Dat008.csv};
            \addlegendentry{Eul POD $m$}
            \addplot[densely dotted,mark size=1.3pt, red] table [x=Number, y=Eul_eig_norm, col sep=comma] {Dat009.csv};
            \addlegendentry{Eul POD $E$}
            
            \addplot[mark size=3pt, blue] table [mark=square*,x=Number, y=ALE_eig_norm, col sep=comma] {Dat010.csv};
            \addlegendentry{ALE POD $\rho$}
            \addplot[dashdotted, mark size=3pt, blue] table [x=Number, y=ALE_eig_norm, col sep=comma] {Dat011.csv};
            \addlegendentry{ALE POD $m$}
            \addplot[densely dotted, mark size=3pt, blue] table [x=Number, y=ALE_eig_norm, col sep=comma] {Dat012.csv};
            \addlegendentry{ALE POD $E$}
        \end{axis}
    \end{tikzpicture}
    \vspace{-2mm}
    \caption{Sod 1D: Eigenvalue decay of the PODs (normalized to have $\lambda_1=1$) non parametric (left) and parametric (right)}
    \label{fig:eig_Sod}
\end{figure}

\begin{table}
\caption{Calibration of the 1D Sod shock tube problem, non parametric and parametric setting}
\centering
\begin{tabular}{|c|c| c||c|c| c|} 
 \hline
 Quantity & Nonparam & Param & Quantity & Nonparam & Param\\[0.5ex]
 \hline
 $\delta$ & $10^{-6}$ & $10^{-6}$& $\alpha$ & $0.$& $0.$\\
 $\rhoref$ & $\rho(t=0.16)$ &- & $M_1$ & 6 &6\\
 max. iter. & $100$ & $100$ &  minim. alg. & \texttt{SLSPQ}& \texttt{SLSPQ}\\
 $N_{\text{train}}$ & 1 & $16$&$N_{\mu}$&25&25\\
 $N_{\text{few}}$& - & $10$ & $N^{\text{POD}}_{\text{few}}$ &-& $3$\\
 \hline
\end{tabular}
\label{tab: calibration Sod non-parametric}
\end{table}

\begin{table}
\caption{Architecture and details of the calibration-NN and the POD-NN in both non parametric and parametric setting}
\centering
\begin{tabular}{|c ||c |c||c|c|} 
 \hline
 & \multicolumn{2}{|c||}{Non parametric case}& \multicolumn{2}{|c|}{Parametric case}\\
 \hline
 Parameter & Calibration-NN & POD-NN& Calibration-NN & POD-NN\\
 \hline
 $L$ & $4$ & 4 & $4$ & 4 \\
 neurons per layer & $16$ & $16$ & 16 & 16 \\
 max. epochs & $20000$& $10000$& $20000$ & $10000$ \\
 loss fun. tol. & $10^{-6}$ & $10^{-5}$& $10^{-6}$ & $10^{-5}$ \\
 $\Tilde{\zeta}$ & \texttt{Tanh}/\texttt{Softplus} & \texttt{Tanh} & \texttt{Tanh}/\texttt{Softplus} & \texttt{Tanh} \\
 $\tau_{\text{POD}}$ & -&$10^{-4}$&$10^{-4}$& $10^{-4}$\\
 $N^{\text{max}}_{\text{POD}}$ & -&$7$&3&7\\
 \hline
 \end{tabular}
\label{tab: NN_architecture_details}
\end{table}

\begin{figure}
    \centering
    {Eulerian approach}\\
    \includegraphics[width=0.32\textwidth,trim={22 20 35 35},clip]{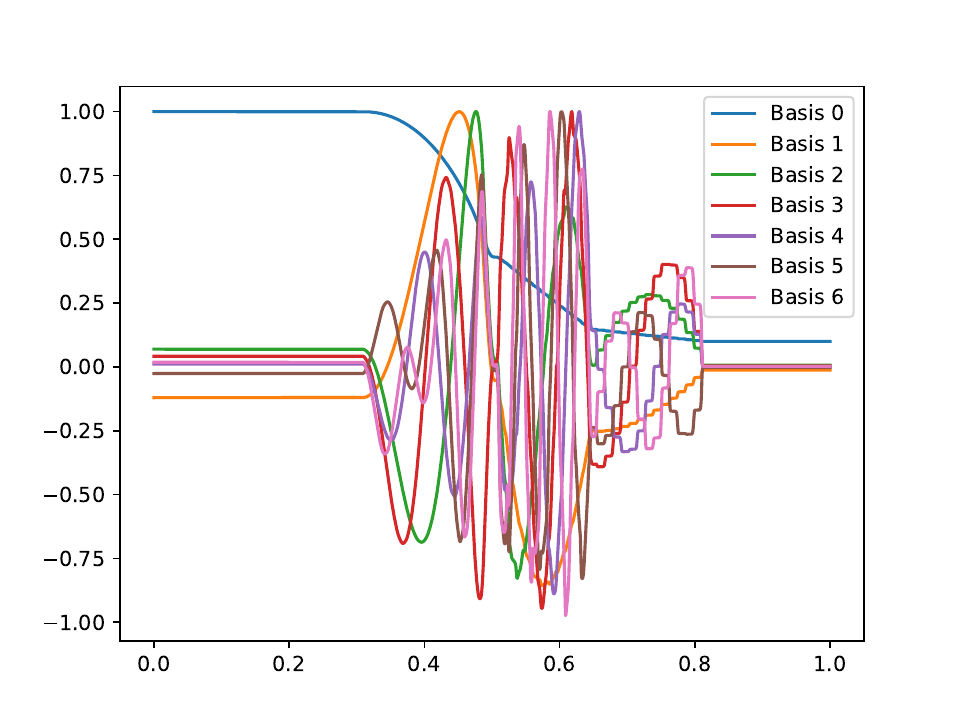}\,
    \includegraphics[width=0.32\textwidth,trim={22 20 35 35},clip]{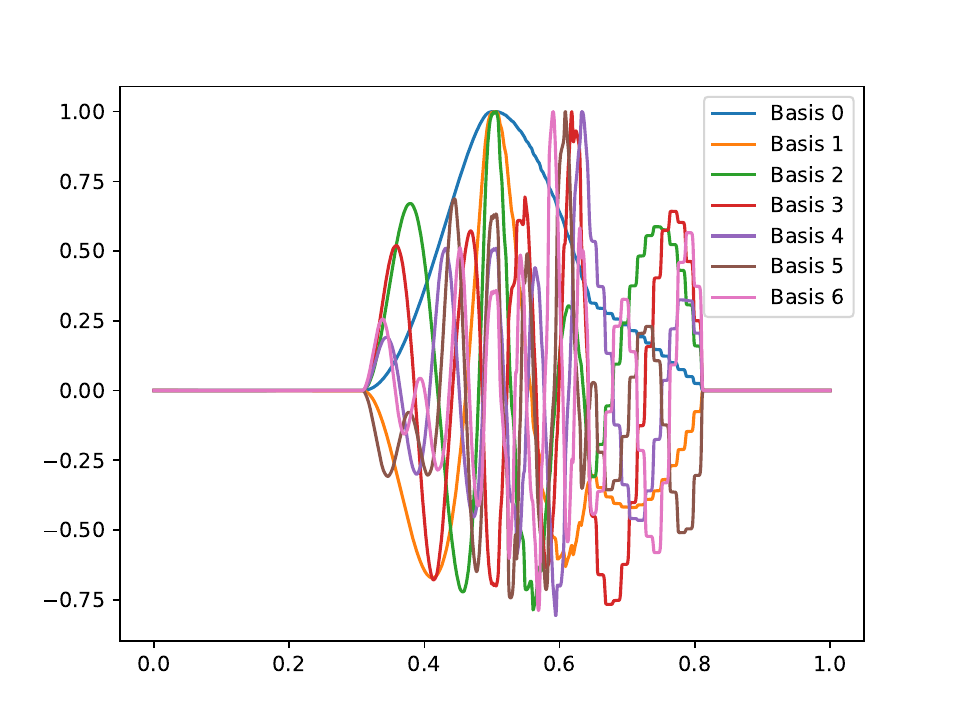}\,
    \includegraphics[width=0.32\textwidth,trim={22 20 35 35},clip]{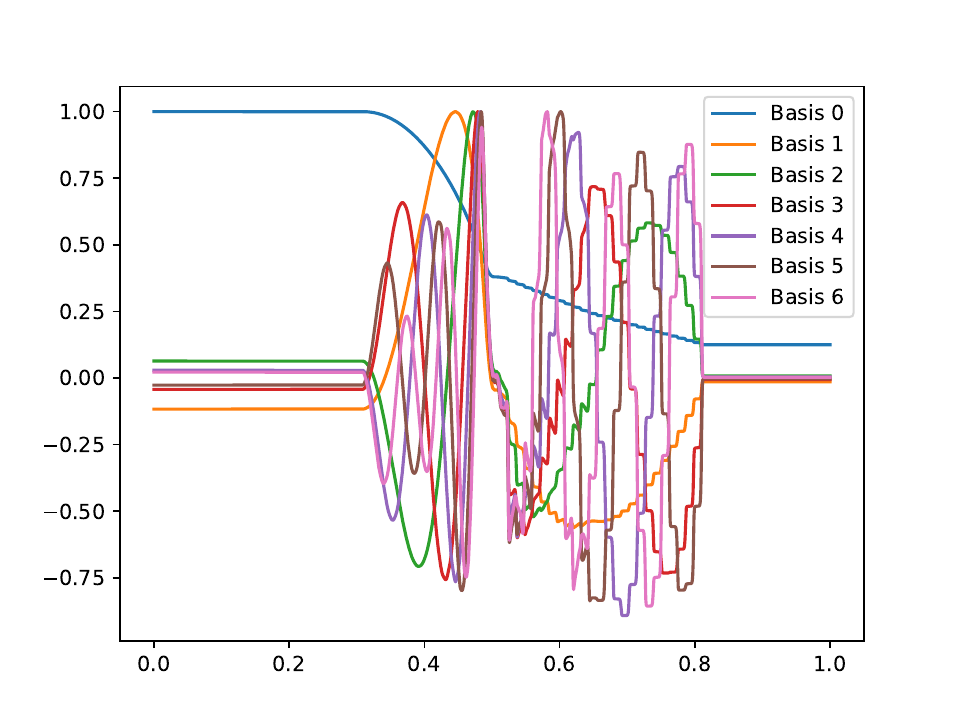}\\[2mm]
    {ALE approach}\\
    \includegraphics[width=0.32\textwidth,trim={22 20 35 35},clip]{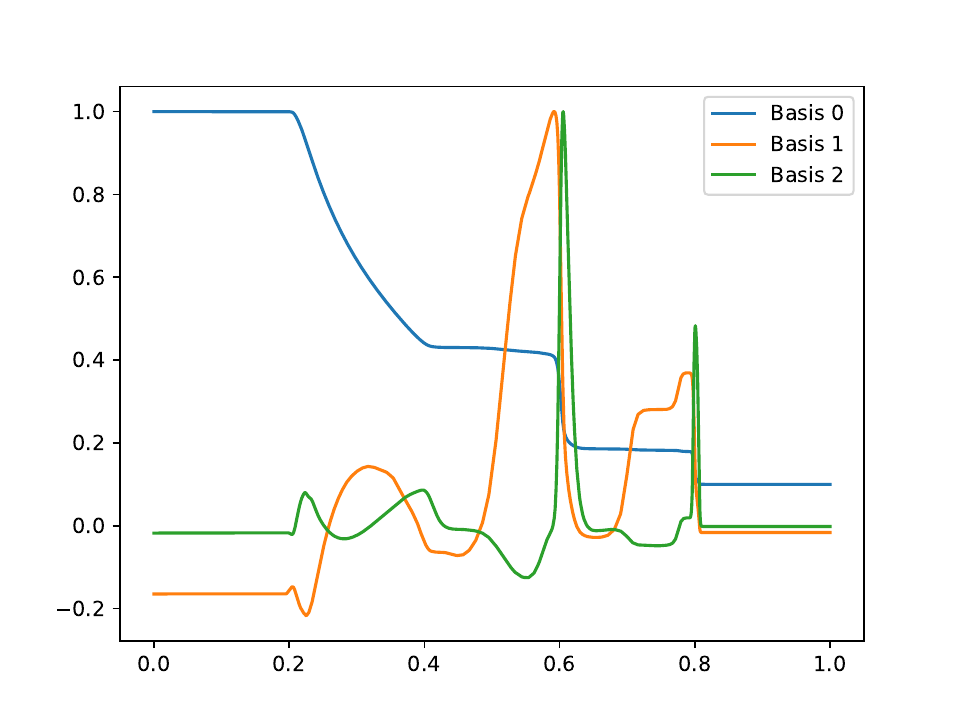}\,
    \includegraphics[width=0.32\textwidth,trim={22 20 35 35},clip]{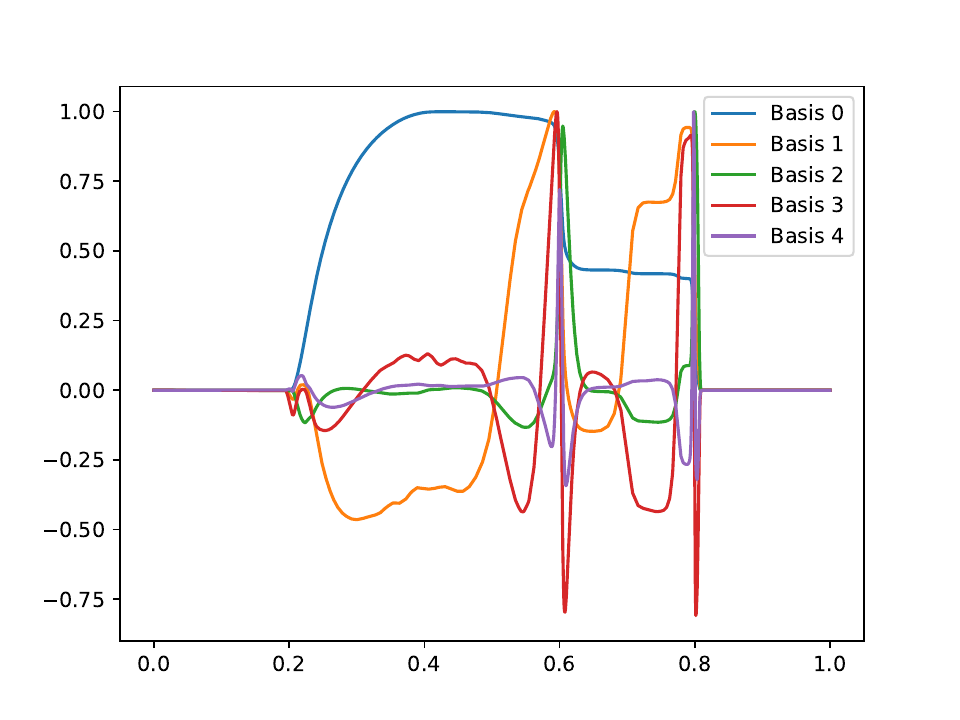}\,
    \includegraphics[width=0.32\textwidth,trim={22 20 35 35},clip]{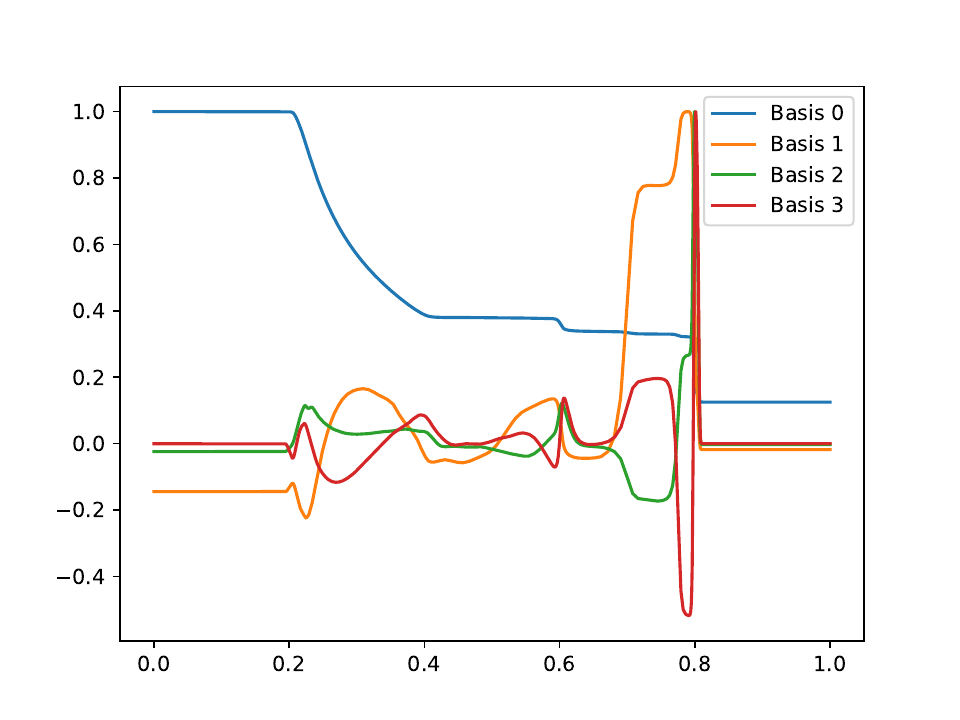}
    \caption{Sod 1D non parametric: The first modes obtained by compressing the original manifolds (Eulerian on the top) and the calibrated manifolds (ALE on the bottom) for $\rho$ (left), $m$ (center) and $E$ (right), with POD with $\tau_{\text{POD}}=10^{-4}$ and $N^{\text{max}}_{\text{POD}}=7$}
    \label{fig:sod_no_param_modes}
\end{figure}

\begin{figure}
    \centering
    {Eulerian ROM}\\
    \includegraphics[width=0.32\textwidth,trim={22 15 35 15},clip]{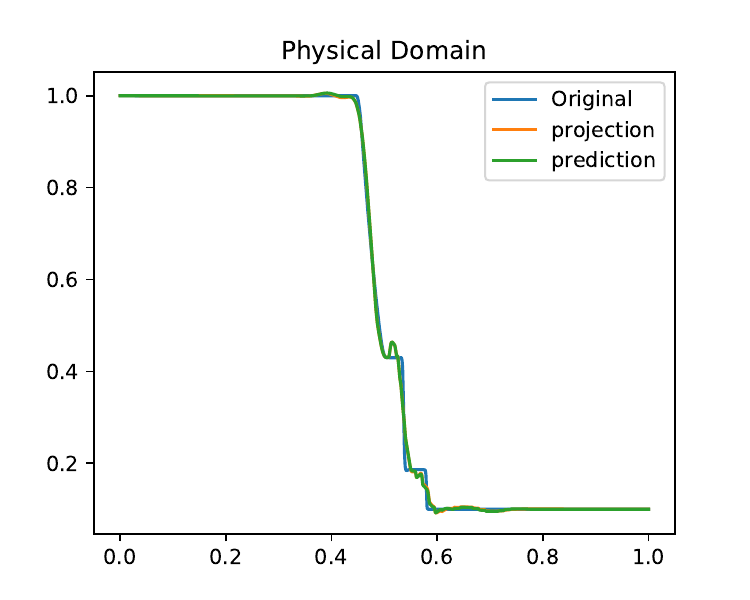}\,
    \includegraphics[width=0.32\textwidth,trim={22 15 35 15},clip]{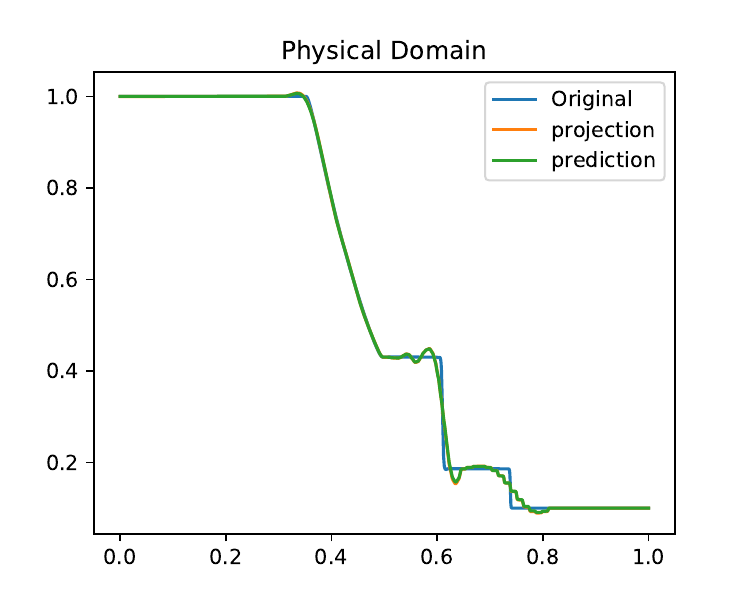}\,
    \includegraphics[width=0.32\textwidth,trim={22 15 35 15},clip]{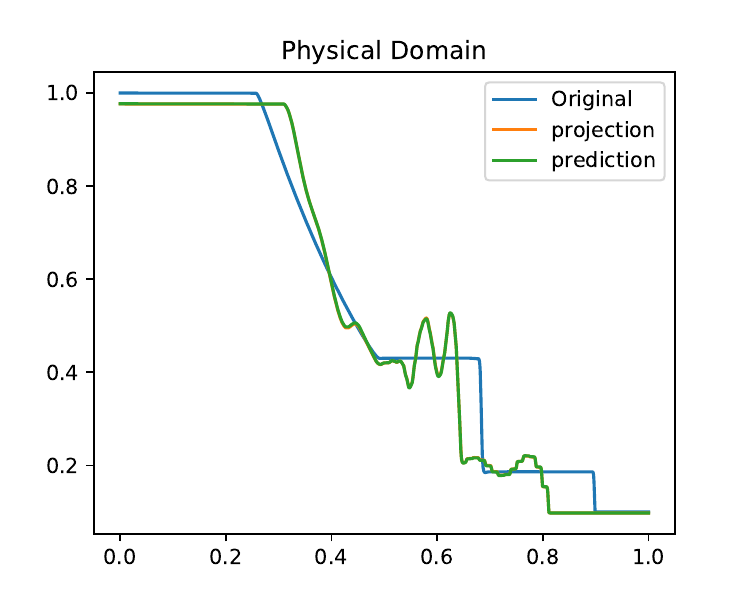}\\[2mm]
    {ALE ROM}\\
    \includegraphics[width=0.32\textwidth,trim={22 20 35 20},clip]{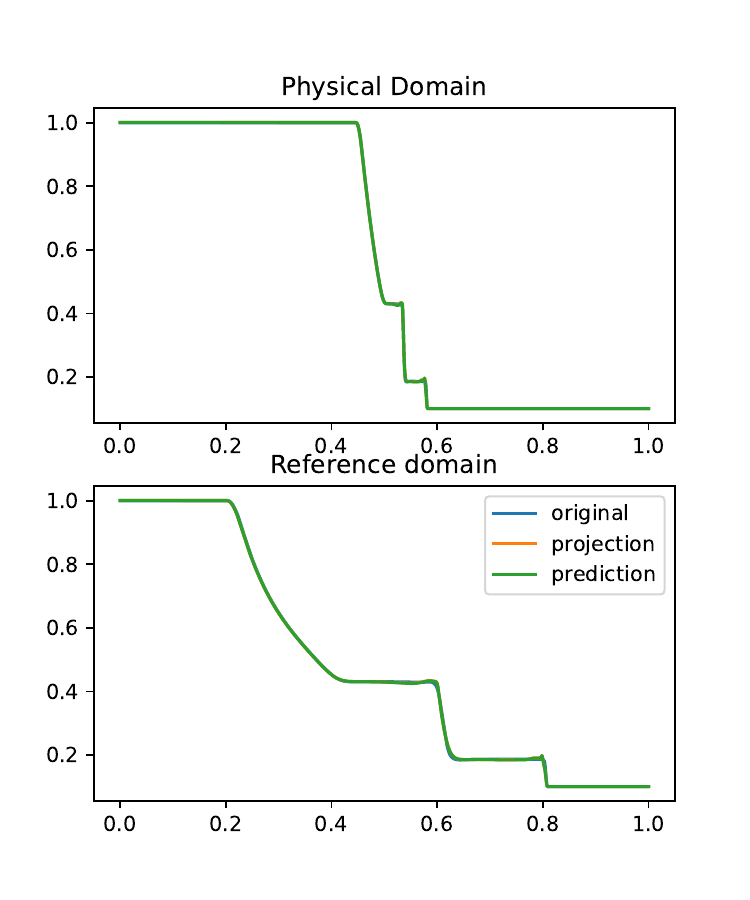}\,
    \includegraphics[width=0.32\textwidth,trim={22 20 35 20},clip]{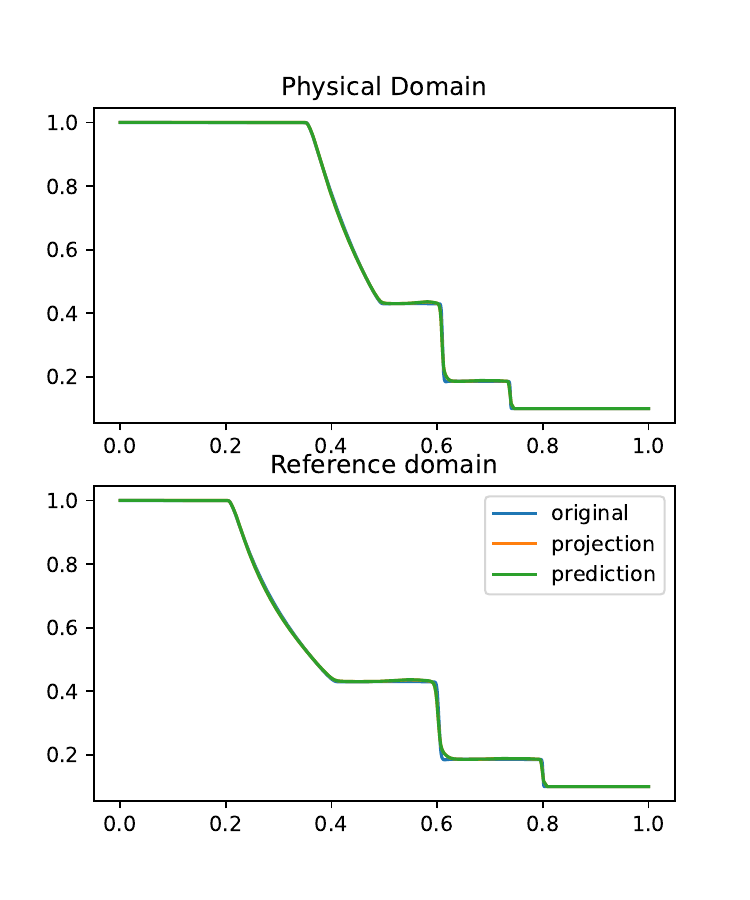}\,
    \includegraphics[width=0.32\textwidth,trim={22 20 35 20},clip]{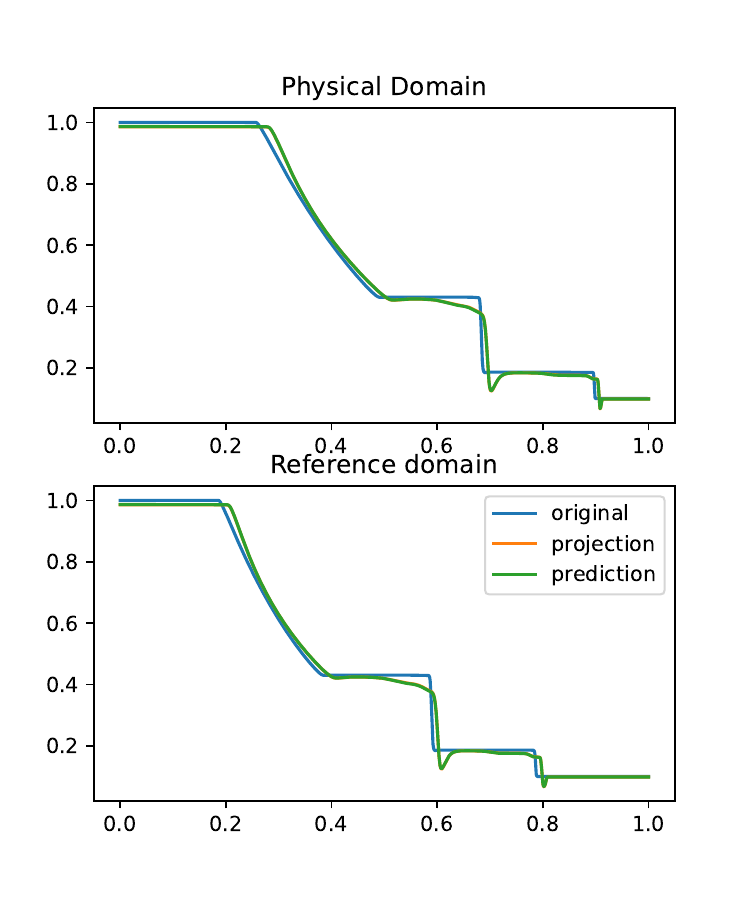}
    \vspace{-2mm}
    \caption{Sod 1D non parametric: Online approximation of the density $\rho$ at times 0.04, 0.12 and 0.2 (left to right). Top row: Eulerian ROM simulations with $N=7$ modes on $\Omega$. Central row: ALE ROM simulations on $\Omega$ with $N=3$ modes. Bottom row: ALE ROM simulations on $\Omegaref$ with $N=3$ modes}
    \label{fig:sod_no_param_ROM}
\end{figure}
Fig.~\ref{fig:eig_Sod} (left) shows the eigenvalue decay of the POD for both calibrated (ALE, in blue) and original (Eulerian, in red) approaches. We can see that, differently from the Eulerian approach, for the calibrated approach the first eigenvalue is much more relevant than all the others and the Kolmogorov $n$--width decay is much faster. 
Fig.~\ref{fig:sod_no_param_modes} shows the behavior of the first modes obtained by compressing with a POD the non-calibrated manifolds, for the three conservative variables $\rho$, $m$ and $E$.
We remark that also here we consider the FOM solutions for $t\in\mathcal{P}^{\text{train}}$, thus excluding the initial and the final times from the compression. The modes are highly oscillatory, because they struggle to correctly represent the positions of the three discontinuities in the domain.
Fig.~\ref{fig:sod_no_param_modes} also shows the POD modes obtained by running a POD on the calibrated manifolds: after the calibration the oscillations in the modes are much milder and focused on the discontinuity locations.
Fig.~\ref{fig:sod_no_param_ROM} shows some POD-NN results with $N=7$, without calibration, for the density $\rho$ at different times $t$ (including the extrapolatory regime at $t=0.2$). 
The comparison is made between the FOM solution $\rho$, its $L^2$ projection onto the reduced space generated by the first $N=7$ modes of the non-calibrated manifold (Eulerian modes) and the online reconstruction obtained using a linear combination of said modes, with coefficients that are predicted by an ANN. 
We can see that the Eulerian approach fails to correctly reproduce the calibrated FOM solutions: in particular, the standard MOR struggles to capture the correct position of the discontinuities, and it shows some oscillations in the online approximation that are most likely due the highly oscillating nature of the Eulerian modes themselves. 
Fig.~\ref{fig:sod_no_param_ROM} also shows some POD-NN simulations obtained after the calibration procedure (computational cost of prediction of both calibration points and ROM coefficients below 0.001s); here we use $N=3$ modes as the POD algorithm stops earlier for the imposed tolerance. 
We can see that we obtain a very good approximation of the calibrated snapshots in the reference configuration $\Omegaref$, i.e., $\rhohat$, and in the physical domain $\Omega$, i.e. $\rho$, with the main features correctly reproduced by the online solution. 
We stress the fact that $t=0.2$ is \emph{outside} the training interval $\mathcal{P}^{\text{train}}$: in this case, the positions of the shock, the contact and the rarefaction wave have been slightly misplaced by the online model, hence, the approximation is not as great as in the interpolatory regime. 
All the details on the architecture of the ANN used to learn the calibration map and to predict the online solution are summarized in Table~\ref{tab: NN_architecture_details}.

\subsubsection{Validation of the calibration strategy}\label{sec:validation_calibration}

\begin{figure}
    \centering
    \includegraphics[width=0.49\linewidth]{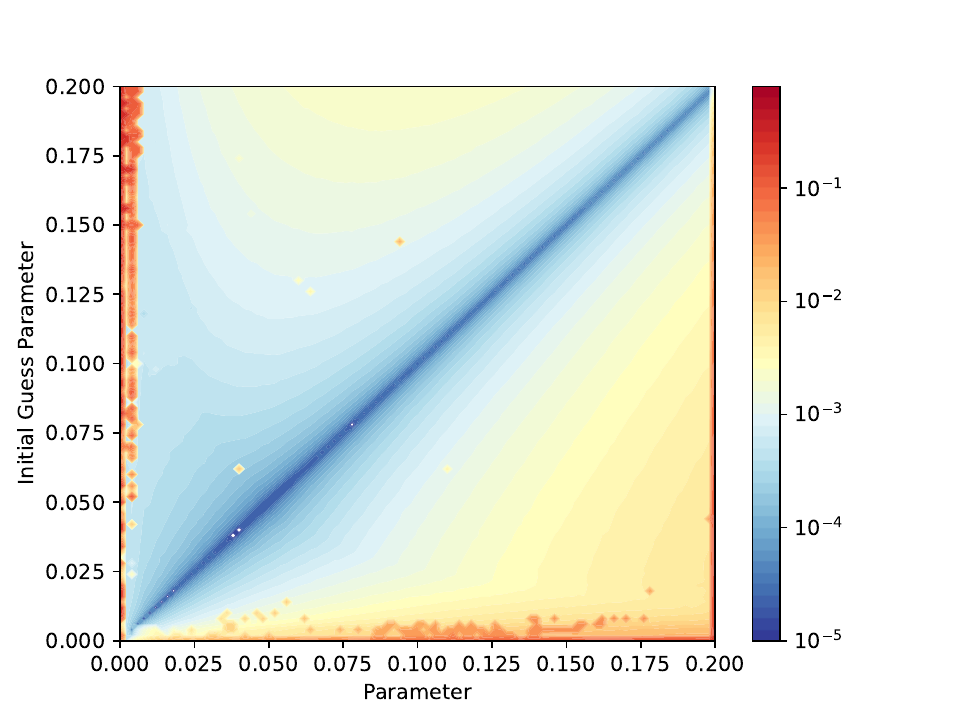}
    \includegraphics[width=0.49\linewidth]{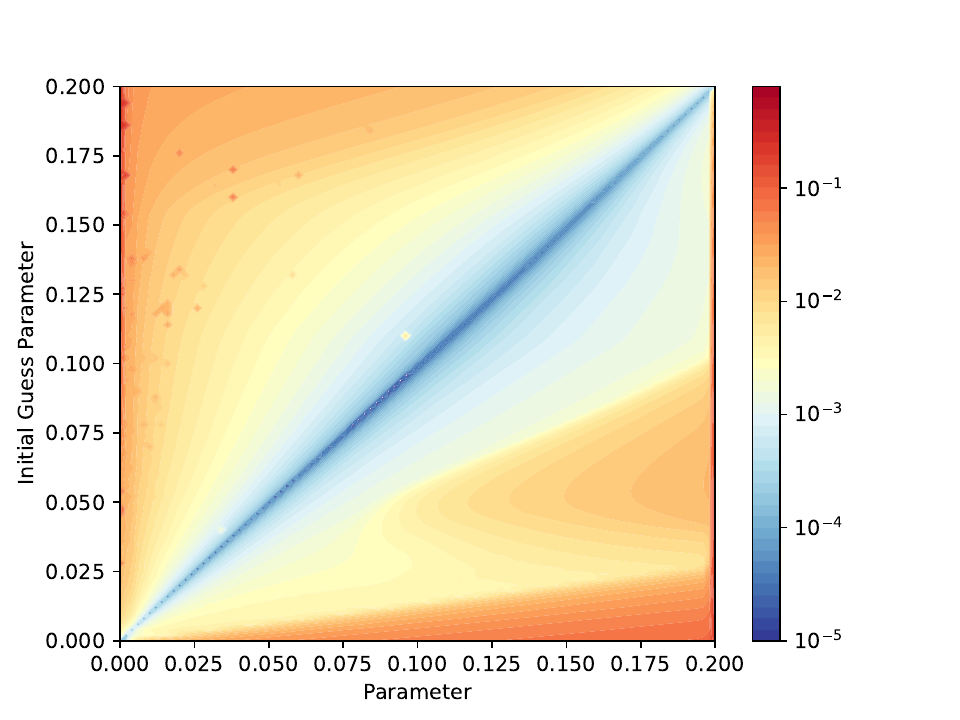}
\caption{Calibration error for Sod 1D at different times (parameter $\bbmu$) using different initial guess times (parameter $\bar\bbmu$), for the FV solutions. Left: $\rhoref=\rho(\bar\bbmu)$, $\bbw^{(0)}(\bbmu)=\bbw^{ex}(\bar\bbmu)$, error measure \eqref{eq:error_characteristics}. Right: $\rhoref=\rho(\bar\bbmu)$ and $\bar\bbw=\bbw^{(0)}(\bbmu) = \lbrace 0.2,0.4,0.6,0.8\rbrace$ for all $\bbmu$, error measure \eqref{eq:error_projection}}
    \label{fig:calibration_all_guesses}
\end{figure}

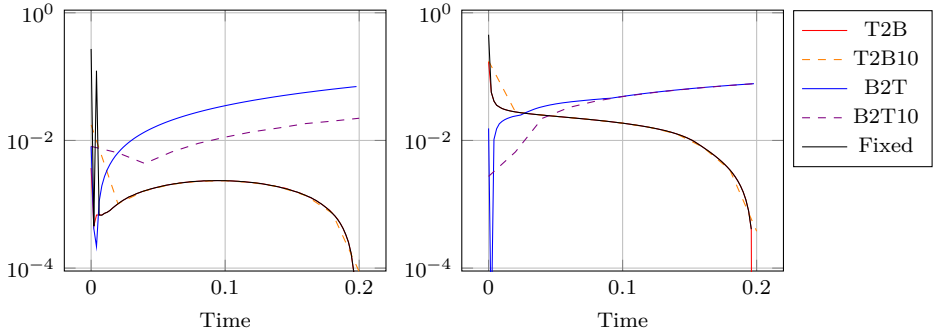
\begin{figure}
    \centering
    \begin{tikzpicture}
        \begin{axis}[ymode=log,
                    xlabel={Time},
                    grid=both,
                    ymin=0.9e-4,
                    ymax=1.1,
                    width=0.45\textwidth,
                    style={font=\footnotesize}]
            \addplot[red] table [x=time, y=calibration_error, col sep=comma] {Dat013.txt};
            \addplot[orange,dashed] table [x=time, y=calibration_error, col sep=comma] {Dat014.txt};
            \addplot[blue] table [x=time, y=calibration_error, col sep=comma] {Dat015.txt};
            \addplot[violet,dashed] table [x=time, y=calibration_error, col sep=comma] {Dat016.txt};
            \addplot[black] table [x=time, y=calibration_error, col sep=comma] {Dat017.txt};
        \end{axis}
    \end{tikzpicture}
    \begin{tikzpicture}
        \begin{axis}[ymode=log,
                    xlabel={Time},
                    grid=both,
                    ymin=0.9e-4,
                    ymax=1.1,
                    width=0.45\textwidth,
                    legend pos=outer north east,
                    style={font=\footnotesize}]
            \addplot[red] table [x=time, y=calibration_error, col sep=comma] {Dat018.txt};
            \addlegendentry{T2B}
            \addplot[orange,dashed] table [x=time, y=calibration_error, col sep=comma] {Dat019.txt};
            \addlegendentry{T2B10}
            \addplot[blue] table [x=time, y=calibration_error, col sep=comma] {Dat020.txt};
            \addlegendentry{B2T}
            \addplot[violet,dashed] table [x=time, y=calibration_error, col sep=comma] {Dat021.txt};
            \addlegendentry{B2T10}
            \addplot[black] table [x=time, y=calibration_error, col sep=comma] {Dat022.txt};
            \addlegendentry{Fixed}
        \end{axis}
    \end{tikzpicture}
    \caption{Calibration error for Sod 1D at different times (parameters) using different calibration strategies, for the FV solutions. Left: exact waves calibration as initial guess and right for equispaced initial guess calibration points. Error measure: \eqref{eq:error_characteristics} (left) and \eqref{eq:error_projection} (right)}
    \label{fig:calibration_strategies_comparison}
\end{figure}

{
In this section, we focus on the two following aspects: the initial guess considered in the calibration, for each parameter, and the order of the parameters for the calibration. To validate our methodology, we perform two types of tests, described in the following paragraphs. We recall that, since we are in the non-parametric case here, the parameter $\bbmu$ represents only the time $t$. Again, we also remark that we focus here on the density $\rho$.}
\paragraph{Reference solution choice in the calibration}
{We are interested in performing the calibration for every parameter, by changing the reference calibration points $\bar\bbw$ through the free coordinates $\bar\bbtheta$, the reference solution $\bar\rho = \rho(\bar\bbmu)$ and the initial guess $\bbtheta^{(0)}(\bbmu)$ in the minimization problem~\eqref{eq:calibration self similar}. We then evaluate the calibration error: the chosen measures for estimating this error will be clarified later on in this paragraph.}
\subparagraph{Exact characteristics calibration reference points}
{For this test, we know the exact solution of the problem~\eqref{eq:euler conservative}, and therefore the exact location of the rarefaction and the discontinuities of the density function. To perform the calibration, we 
choose a reference parameter $\bar\bbmu\in\mathcal{P}$, and we consider $\rhoref=\rho(\bar\bbmu)$ and $\bar\bbtheta=\bbtheta^{ex}(\bar\bbmu)$, where $\bbtheta^{ex}(\bar\bbmu)\in \mathbb R^4$ is the vector of the exact locations of the beginning and end of the rarefaction wave, the contact and the shock discontinuities for $\rho(\bar\bbmu)$. We then carry on the calibration, for every $\bbmu\in\mathcal{P}^{train}$: the initial guess $\bbtheta^{(0)}(\bbmu)$ will be $\bbtheta^{ex}(\bar\bbmu)$. We study what happens if we change $\bar\bbmu$: the measure that we consider to estimate the calibration error in this case is the $1$-norm of the difference between $\bbtheta^{ex}(\bbmu)$ (the vector of exact locations of beginning and end of the rarefaction wave, the contact and the shock discontinuities for the $\rho(\bbmu)$) and $\bbtheta^{opt}(\bbmu)$ (the output of the calibration for the parameter $\bbmu$), namely 
\begin{equation}
    \label{eq:error_characteristics}|| \bbtheta^{opt}(\bbmu)-\bbtheta^{ex}(\bbmu)||_1.
\end{equation}
This choice allows us to measure how far our output is from the actual location of the discontinuities and rarefaction of the solution.
The result is reported in Fig.~\ref{fig:calibration_all_guesses} (left). Using the initial guess close to the final time results in the largest errors when calibrating early time solutions, while the opposite seem less problematic and optimally one could choose a reference solution close to the initial time to have the least amount of error for all the considered optimized parameters. It is also clear from this test that the best choice would be to choose the initial guesses close to the parameter that we are calibrating minimizes the error, i.e. staying close to the diagonal. And this is what we are actually doing in presented strategy in Section~\ref{sec:calibration self similar} in \eqref{eq:initial guess parametric setting}.}
\subparagraph{Equispaced calibration reference points}
{For this second test, we assume we do not know the exact solution, nor the wave locations. In this case, we choose as reference calibration points $\bar{\bbtheta} = \lbrace 0.2,0.4,0.6,0.8 \rbrace$ equispaced points in $[0,1]$: these points will be fixed for all the reference parameters $\bar\bbmu$ we will use and they will also be used as initial guess of the calibration process. We then choose $\rhoref=\rho(\bar\bbmu)$, and we analyse the results of the calibration when we change $\bar{\bbmu}$. 
The measure chosen to study the error of the calibration (since we do not know $\bbtheta^{ex}(\bbmu)$) is the following: 
\begin{equation}\label{eq:error_projection}\lvert\lvert\rhohat(\bbmu, \bbw(\bbtheta^{opt}(\bbmu))) - \Pi_{\rhoref}\rhohat(\bbmu, \bbw(\bbtheta^{opt}(\bbmu)))\rvert\rvert_2^2,
\end{equation}
where $\Pi_{\rhoref}$ is the projection onto the linear space spanned by the chosen reference solution.
The results are reported in Fig.~\ref{fig:calibration_all_guesses}: again, we can see that choosing a reference solution close to the parameter of interest leads to smaller errors and that using one reference solution the optimal choice would be something close to $\bar\bbmu \approx 0.1$. 
}
\paragraph{Order in which we perform the calibration}
{
\begin{table}
\centering
    \begin{tabular}{|c|c|c|c|}
    \hline
         Strategy & Reference sol & Order & Initial guess  \\ \hline
         T2B &  Final time & {Last time to first} & Last computed\\
         T2B10 &  Final time &{Last time to first, every 10 steps}  & Last computed\\
         B2T &  Initial time & {First time to last}  & Last computed\\
         B2T10 &  Initial time & {First time to last, every 10 steps}  & Last computed\\
         Fixed &  Final time &{Last time to first} & Final time\\
         \hline
    \end{tabular}
    \caption{Details for the order calibration test, where we analyze the calibration output according to the order (in time) in which we perform it}
    \label{tab:details_calibration_order}
\end{table}
Motivated by the results shown in Fig.~\ref{fig:calibration_all_guesses}, we now use as initial guess the calibration points found for the closest parameter already computed. We will compare the different orders described in Table~\ref{tab:details_calibration_order}. In particular, we will try to start from the final time or the initial one and we will use either all 100 timesteps we have in the training set or only one every 10 of these.
Finally, we also compare it with a fixed initial guess strategy.
As we can observe in Fig.~\ref{fig:calibration_strategies_comparison}, starting with the final solution as reference solution, we obtain considerably lower errors all along the time domain, almost independently on how many parameters we sample (this will not be true for more complicated tests in 2D), with respect to using the initial solution as a reference one.
The authors are anyway surprised by how the error using the initial solution as reference is not too large. 
Also the fixed initial guess at the final time produces interesting results.}

{In 2D tests, we cannot perform similar comparison, as no exact calibration could be performed. 
We can undoubtedly say that in 2D being close to the initial guess is of more importance than in 1D, hence, in general, it is better to have a large enough training set.
}
\subsection{Parametric Sod shock tube problem in 1D}\label{sec:numerical results Sod parametric}

\begin{figure}
    \centering
    {Eulerian approach}\\
    \includegraphics[width=0.32\textwidth,trim={22 20 35 25},clip]{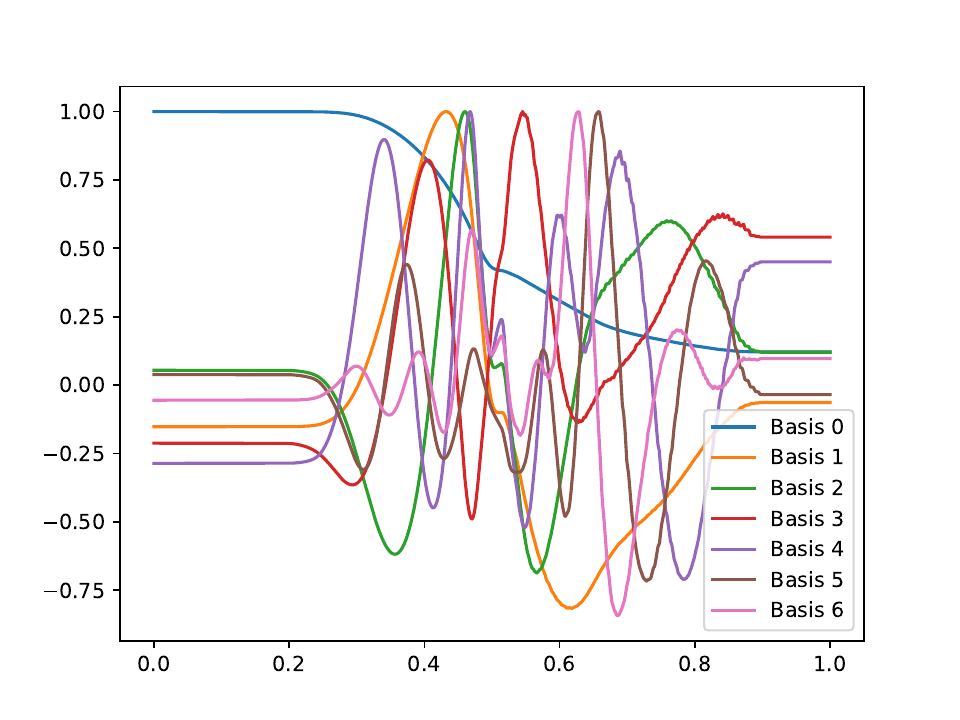}\,
    \includegraphics[width=0.32\textwidth,trim={22 20 35 25},clip]{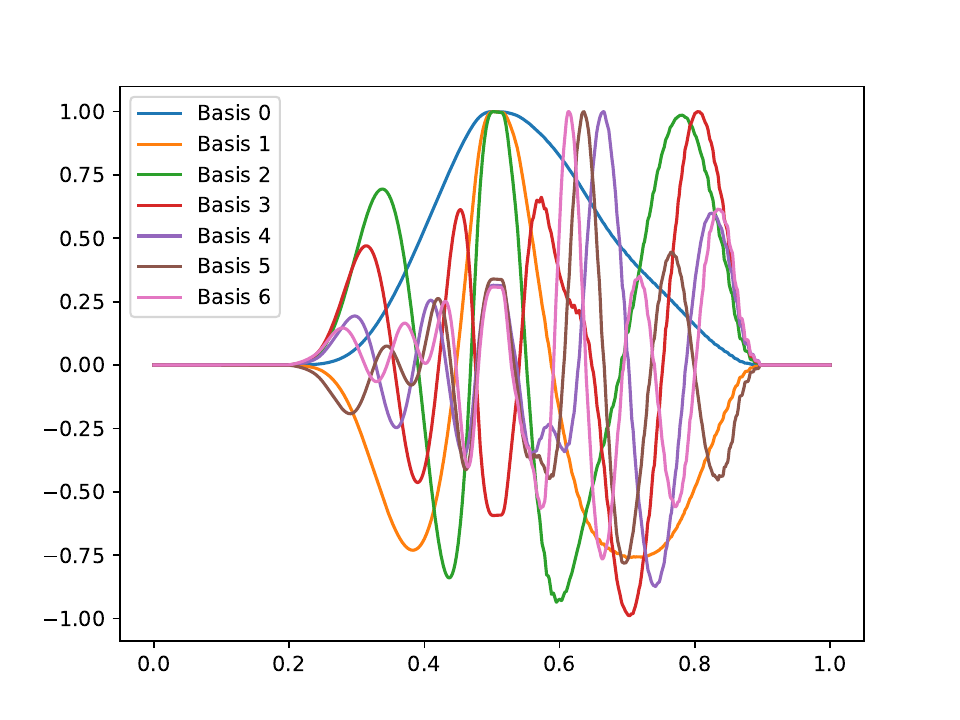}\,
    \includegraphics[width=0.32\textwidth,trim={22 20 35 25},clip]{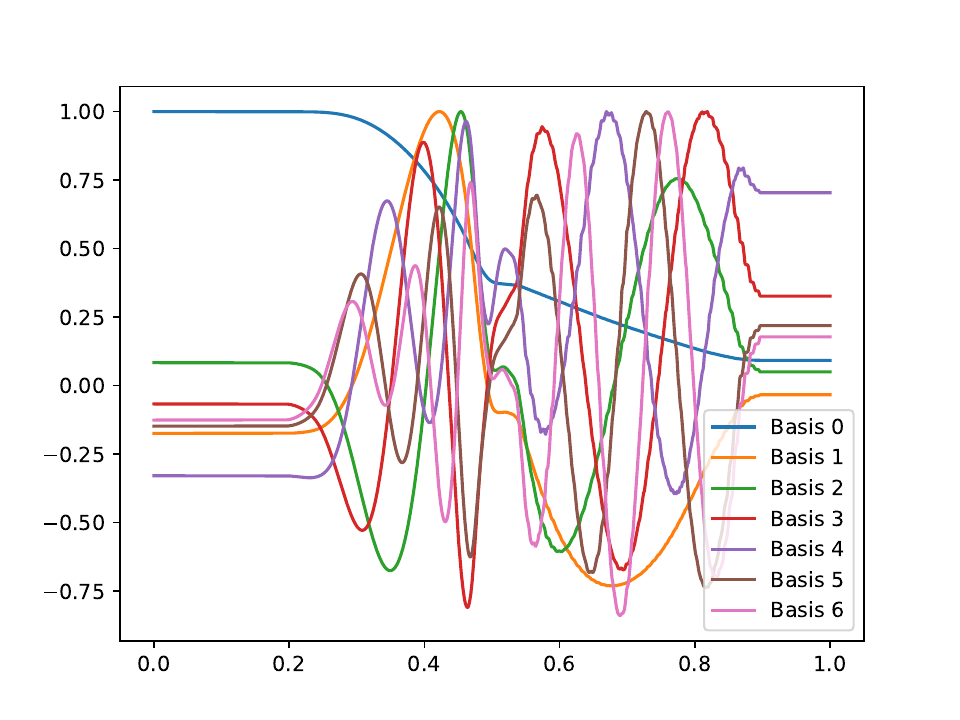}\\
    {ALE approach}\\
    \includegraphics[width=0.32\textwidth,trim={22 20 35 25},clip]{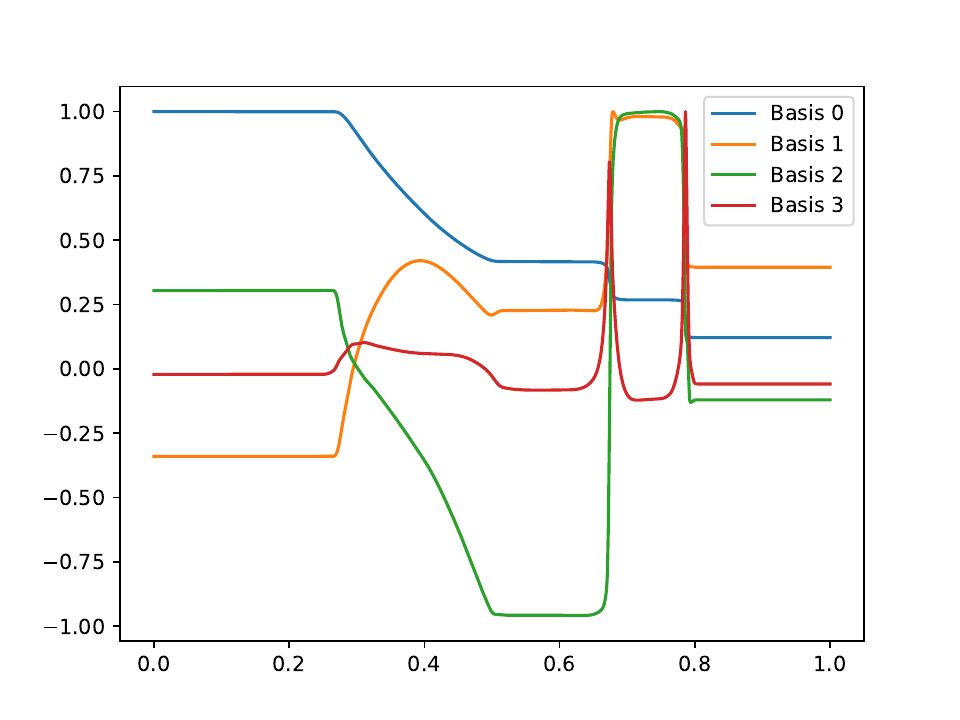}
    \includegraphics[width=0.32\textwidth,trim={22 20 35 25},clip]{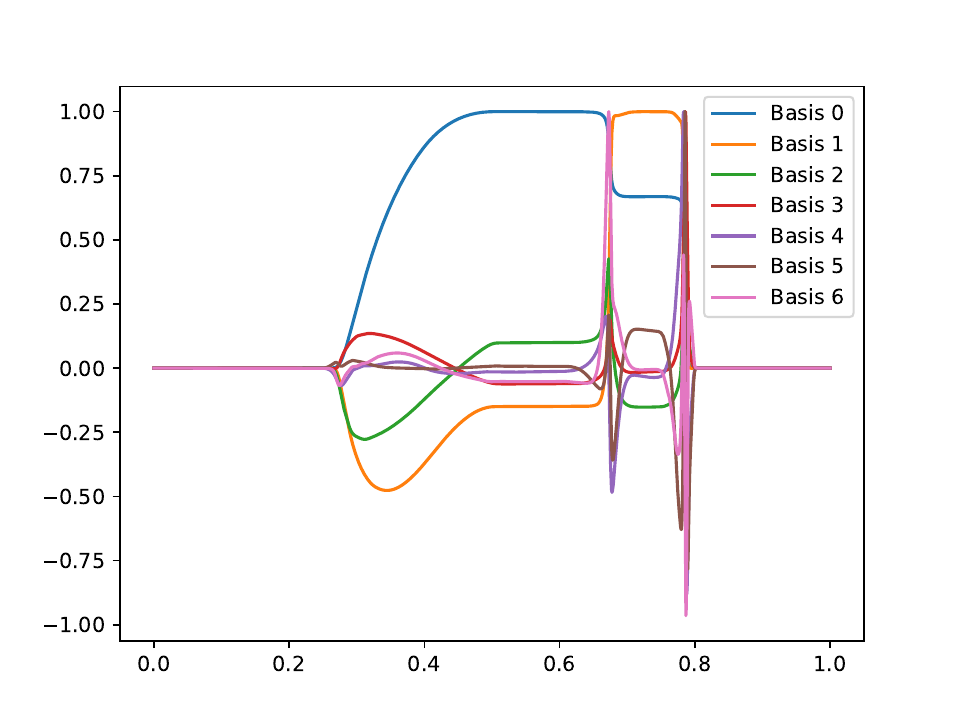}
    \includegraphics[width=0.32\textwidth,trim={22 20 35 25},clip]{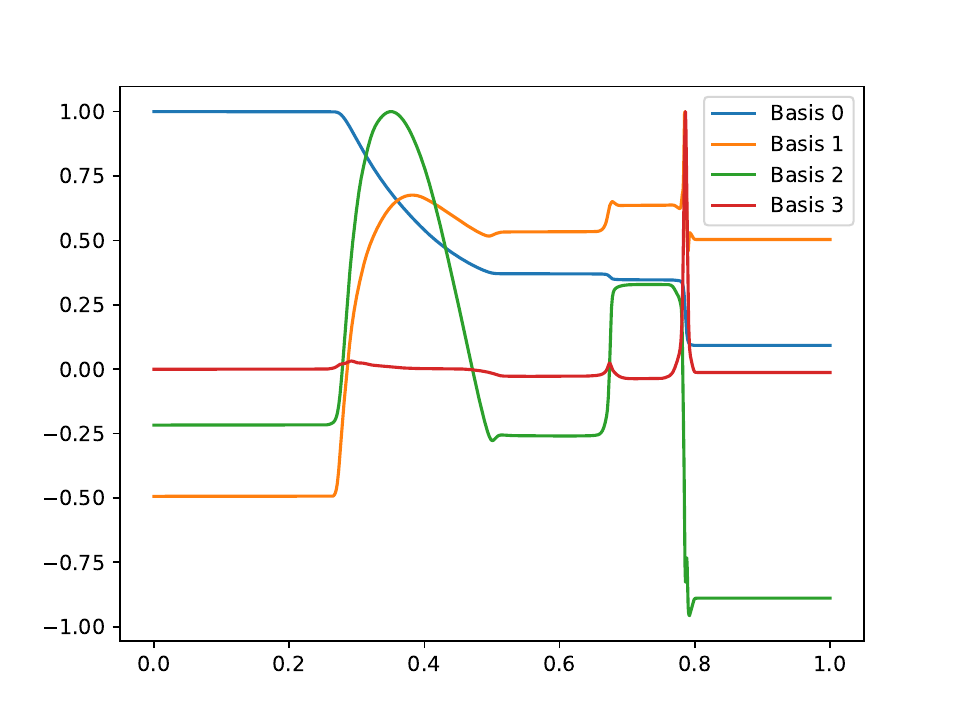}
    \vspace{-1mm}
    \caption{Sod 1D parametric: The first modes obtained by compressing the original manifolds (Eulerian on the top) and the calibrated manifolds (ALE on the bottom) for $\rho$ (left), $m$ (center) and $E$ (right), with POD with $\tau_{\text{POD}}=10^{-4}$ and $N^{\text{max}}_{\text{POD}}=7$}
    \label{fig:sod_param_POD modes}
\end{figure}

\begin{figure}
    \centering
    {Eulerian ROM}\\
    \includegraphics[width=0.32\textwidth,trim={22 15 35 15},clip]{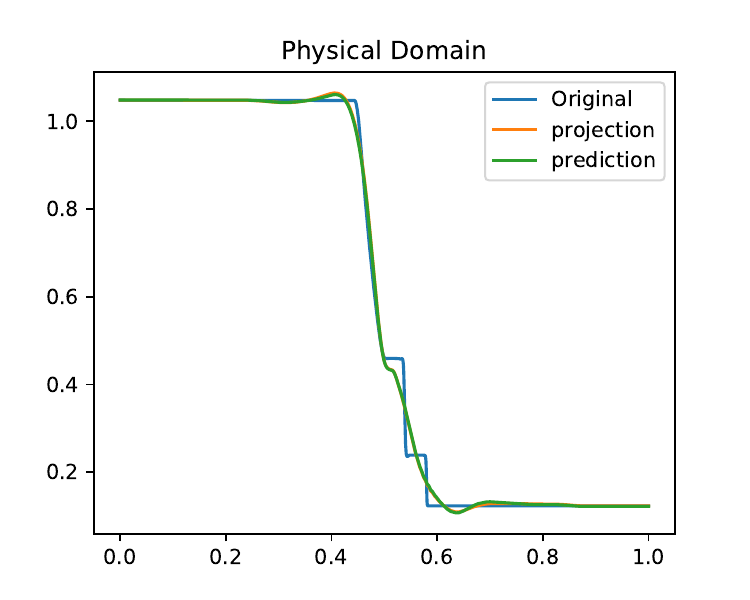}\,
    \includegraphics[width=0.32\textwidth,trim={22 15 35 15},clip]{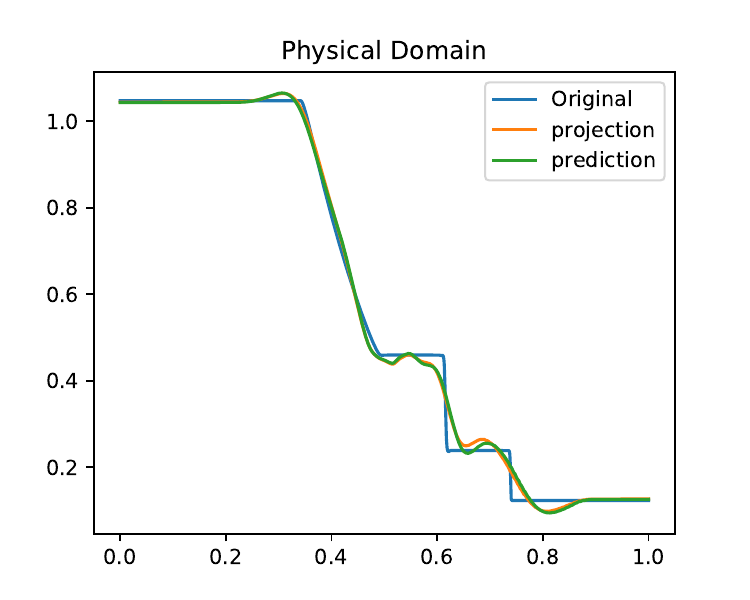}\,
    \includegraphics[width=0.32\textwidth,trim={22 15 35 15},clip]{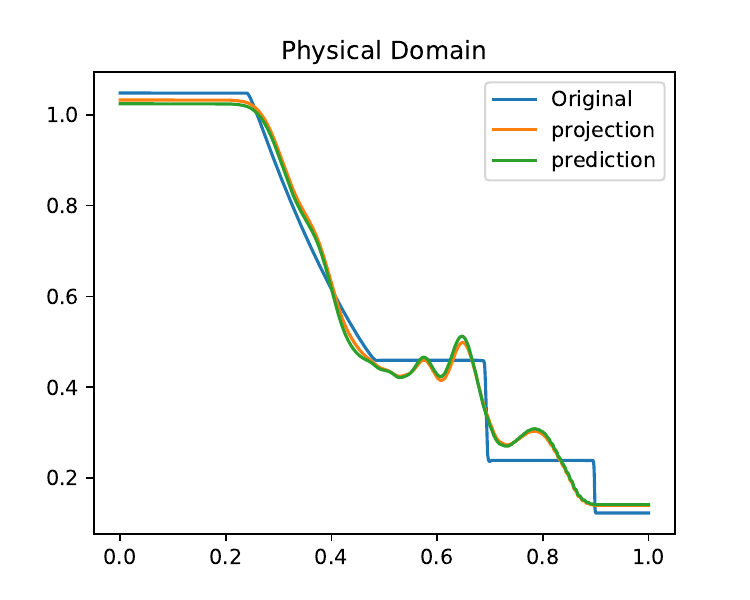}\\
    {ALE ROM}\\
    \includegraphics[width=0.32\textwidth,trim={22 30 35 30},clip]{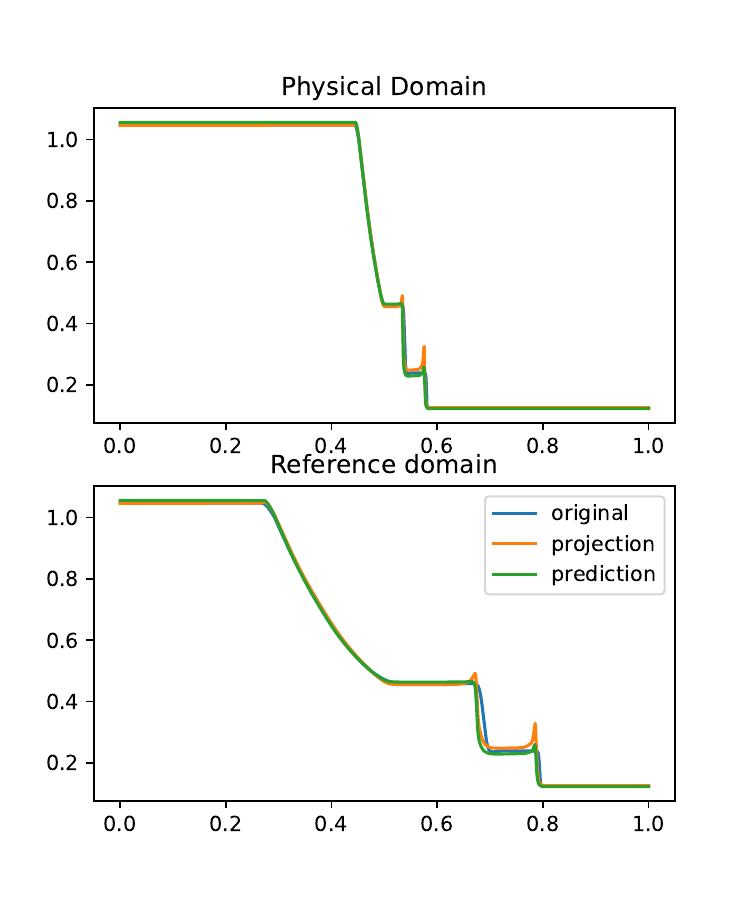}\,
    \includegraphics[width=0.32\textwidth,trim={22 30 35 30},clip]{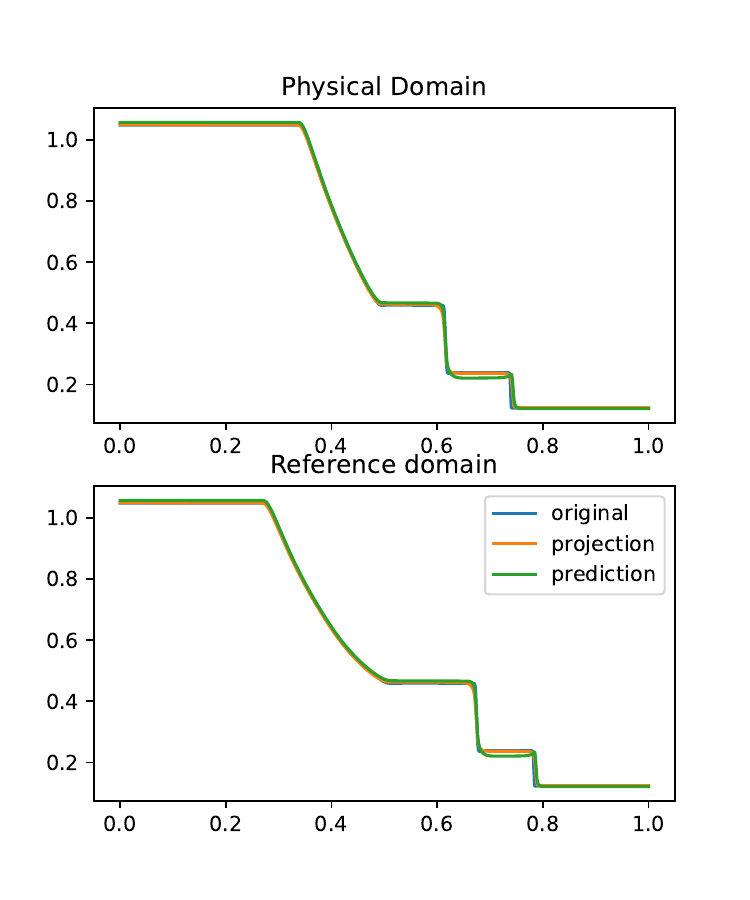}\,
    \includegraphics[width=0.32\textwidth,trim={22 30 35 30},clip]{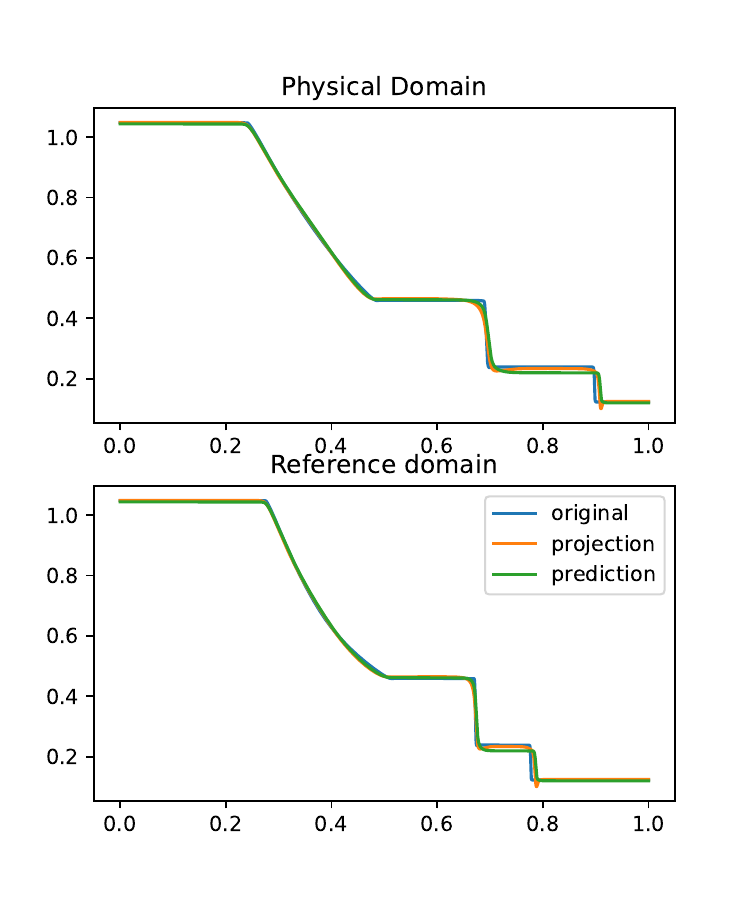}
    \vspace{-1mm}
    \caption{Sod 1D parametric: Online approximation of the density $\rhoref$ at times 0.04, 0.12 and 0.2 (left to right), for $\rho_L=1.047937,\, \rho_R=0.122810,\, p_L=1.203980,\, p_R=0.144468$. Top row: Eulerian ROM simulations on $\Omega$ with $N=7$ modes. Central row: ALE ROM simulations on $\Omega$ with $N=4$ modes. Bottom row: ALE ROM simulations on $\Omegaref$ with $N=4$ modes}
    \label{fig:sod_param_mu0_ROM}
\end{figure}

\begin{figure}
    \centering
    {Eulerian ROM}\\
    \includegraphics[width=0.32\textwidth,trim={22 15 35 15},clip]{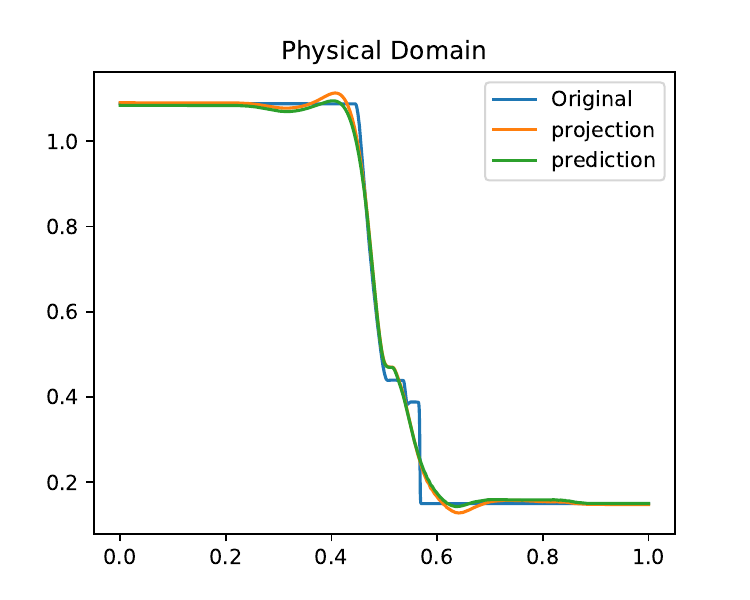}\,
    \includegraphics[width=0.32\textwidth,trim={22 15 35 15},clip]{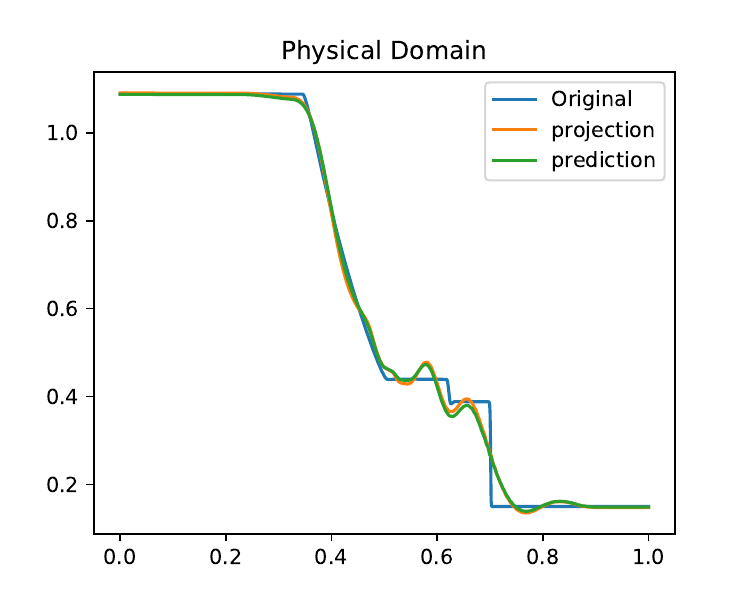}\,
    \includegraphics[width=0.32\textwidth,trim={22 15 35 15},clip]{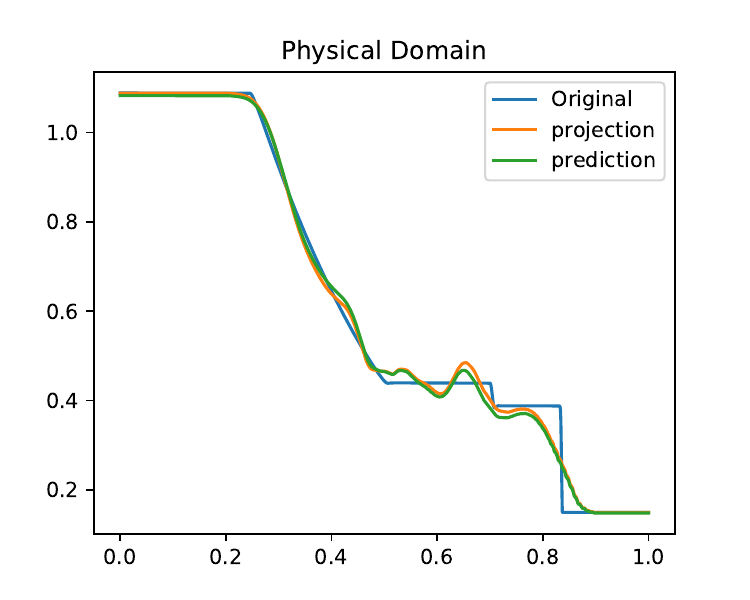}\\
    {ALE ROM}\\
    \includegraphics[width=0.32\textwidth,trim={22 30 35 30},clip]{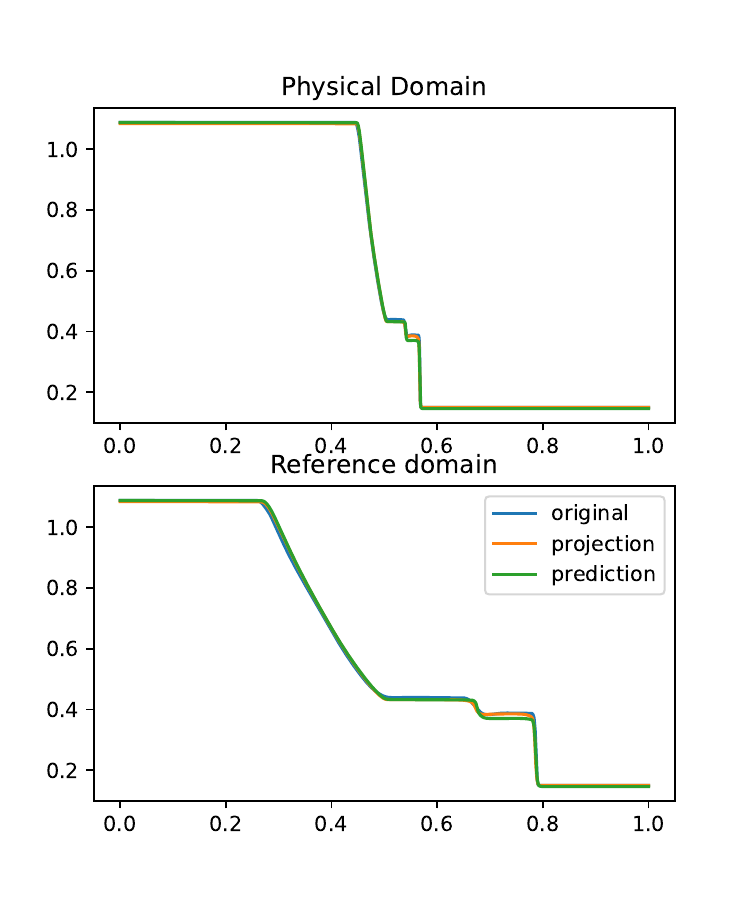}\,
    \includegraphics[width=0.32\textwidth,trim={22 30 35 30},clip]{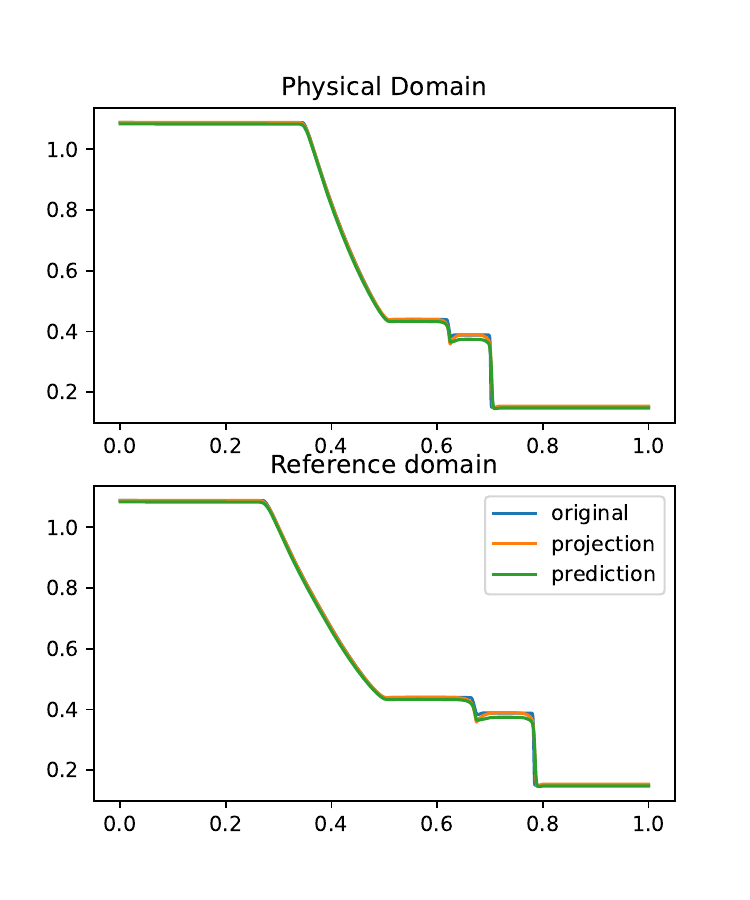}\,
    \includegraphics[width=0.32\textwidth,trim={22 30 35 30},clip]{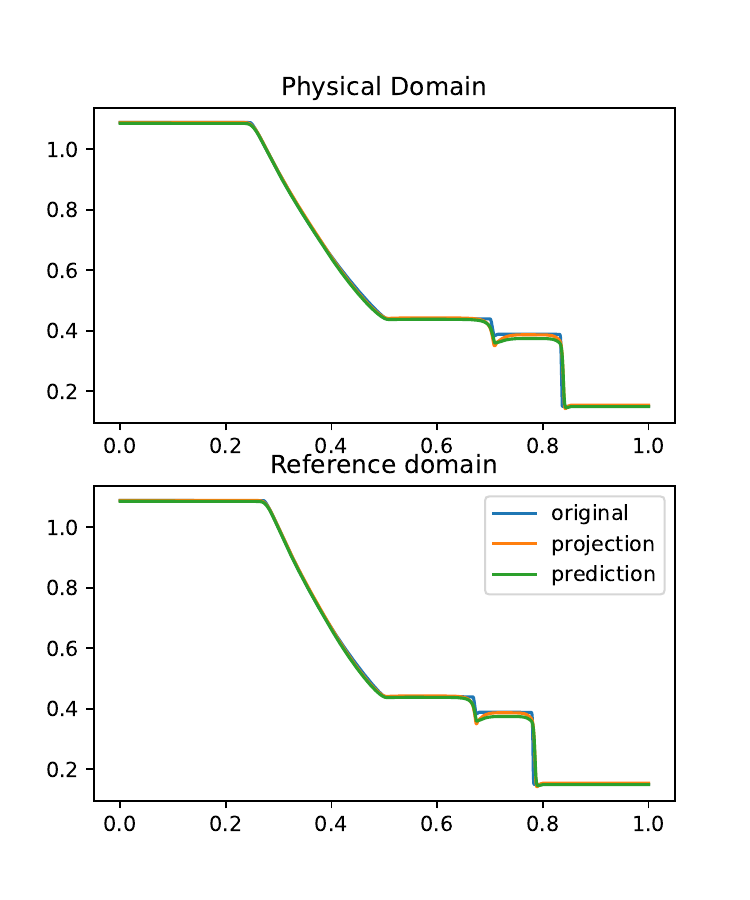}
    \caption{Sod 1D parametric: Online approximation of the density $\rhoref$ at times 0.04, 0.12 and 0.2 (left to right), for $\rho_L=1.08827,\, \rho_R=0.149654,\, p_L=1.193154,\, p_R=0.078459$. Top row: Eulerian ROM simulations on $\Omega$ with $N=7$ modes. Central row: ALE ROM simulations on $\Omega$ with $N=4$ modes. Bottom row: ALE ROM simulations on $\Omegaref$ with $N=4$ modes}
    \label{fig:sod_param_mu2_ROM}
\end{figure}
We now consider the parametric version of the Sod problem, already introduced in Section \ref{sec:calibration quasi self similar}. We recall that we consider $\mu=(\rho_L, \rho_R, p_L, p_R)\in\mathcal{P}_{\text{phys}}\subset\mathbb{R}^4$. All the details for the numerical simulations are provided in Table \ref{tab:Sod non-parametric}. 
Also in this case, the FOM solutions have been obtained employing the same FV discretization. 
In the parametric setting, we generate the training space $\mathcal{P}_{\text{phys}}^{\text{train}}$ using $N_{\text{train}} = 16$ randomly selected parameters $\mu$ from $\mathcal{P}_{\text{phys}}$. Again, we consider the training time interval $[0.01,0.16]$ discretized with around 45 times for each physical parameter.

Fig.~\ref{fig:sod_param_POD modes} shows the first modes for the three components $\rho$, $m$ and $E$, without calibration (Eulerian approach) and with calibration (ALE approach), respectively. 
Also in the parametric case, we can notice that the Eulerian modes are highly oscillating, similarly to the non-parametric test case: the calibration helps to significantly mitigate this phenomenon. 
To further validate this, we show in Fig.~\ref{fig:eig_Sod} (right) a comparison between the rate of decay of the eigenvalues obtained with a POD on the non-calibrated (red) and calibrated (blue) manifolds. The calibration results in an improvement in the rate of decay and we clearly observe that, in comparison to the non parametric case, the decay is slower and we need more basis functions to represent our solution manifold.
All the details for the numerical implementation of the calibration procedure are summarized in Table~\ref{tab: calibration Sod non-parametric}. 

Fig.~\ref{fig:sod_param_mu0_ROM} and Fig.~\ref{fig:sod_param_mu2_ROM} represent the FOM, the $L^2$ projection on the reduced space and the POD--NN online approximation for $\rho$, for two parameters in the test set. 
We plot both the Eulerian ROM and the ALE one, for the latter both in the physical $\Omega$ and in the reference domain $\Omegaref$.
The online approximations are obtained with $N=4$ modes.
As we can see, in both cases the Eulerian ROM is struggling to correctly capture the positions of the discontinuities, and it provides an approximated solution that exhibits some non-negligible oscillations, most likely due to the oscillating nature of the Eulerian modes themselves. 
The results provided with the calibration are much more accurate, since the MOR is now able to correctly represent the positions of the discontinuities, and it does not present any oscillations in the approximations.
There are minor flaws in the extrapolatory regime and in the early times, still keeping the quality of the reduced solution very high. 

All the details on the architecture of the ANN used to learn the calibration map and to predict the online solution are summarized in Table~\ref{tab: NN_architecture_details}.

\subsection{Non-parametric DMR problem in 2D}\label{sec:numerical results DMR non-parametric}
We now consider a 2D test case, namely the Double Mach Reflection (DMR) problem \cite{woodward}. 
Let $\Omega = [0, 4]\times[0,1]$: we consider the Euler equations \eqref{eq:euler conservative}, in the time interval $[0, 0.25]$, with the following IC
\begin{align}
    &\begin{cases}
        (\rho_L, u_L, v_L, p_L) = ( 8, 8.25 \cos(\beta) , - 8.25 \sin(\beta), 116.5  ) & \bbx \in \Omega_L(\beta,t=0)\\
        (\rho_R, u_R, v_R, p_R) = ( 1.4, 0 , 0, 1  ) & \bbx \in \Omega_R(\beta,t=0) 
    \end{cases}\\
    &\Omega_L(\beta,t) = \left\lbrace \bbx \in \Omega : x<\frac16 + \tan(\beta) y + \frac{10}{\cos(\beta)}t \right\rbrace, \qquad \Omega_R(\beta,t)=\Omega \setminus \Omega_L(\beta,t),
\end{align}
with $\beta=\frac{\pi}{6}$. The BCs are assigned through ghost cells as 
\begin{equation}
    (\rho, u,v,p)=\begin{cases}
        (\rho_L, u_L,v_L,p_L), & \text{if }x=0 \text{ or }(y=1 \text{ and } \bbx \in \Omega_L(\beta,t)),\\
        (\rho_R, u_R,v_R,p_R), & \text{if } y=1 \text{ and } \bbx \in \Omega_R(\beta,t),\\
        (\rho_{in}, u_{in},-v_{in},p_{in}) ,&\text{if }y=0 \text{ (wall BC)},\\
        (\rho_{in}, u_{in},v_{in},p_{in}),&\text{if }x=4 \text{ (outflow)},
    \end{cases}
\end{equation}
where $\cdot_{in}$ denotes the value inside the domain at the corresponding boundary cell.

\begin{table}
\caption{Calibration of the 2D DMR problem, non parametric and parametric setting.}
\centering
\begin{tabular}{|c|c| c||c|c| c|} 
 \hline
 Quantity & Nonparam & Param & Quantity & Nonparam & Param\\[0.5ex]
 \hline
 $\rhoref$ & $\rho(t=0.2)$ &-& $\alpha$ & $10^{-4}$& $10^{-4}$\\
 $M_1$ & 7 &7 & $\delta$ & $10^{-2}$ & $10^{-1}$  \\
 $M_2$ & 6 &6 & minim. alg. & \texttt{SLSPQ}& \texttt{SLSPQ} \\
 $N_{\text{train}}$ & 1 & $16$ & max. iter. & $100$ & $100$ \\
 $N_{\mu}$ & 180 & 23 &&&\\
 \hline
\end{tabular}
\label{tab: calibration DMR non-parametric}
\end{table}

\begin{figure}
    \centering
    \begin{tikzpicture}
        \begin{axis}[ymode=log,
                    xlabel={$N$},
                    xmax=20,
                    ymin=1e-6,
                    grid=both,
				ylabel={Eigenvalue},
                    ytick={1,0.01,0.0001,0.000001},
                    width=0.45\textwidth,
                    style={font=\footnotesize}]
            \addplot[red] table [x=Number, y=Eul_eig_norm, col sep=comma] {Dat023.csv};
            \addlegendentry{Eul POD}
            \addplot[blue] table [x=Number, y=ALE_eig_norm, col sep=comma] {Dat024.csv};
            \addlegendentry{ALE POD}
        \end{axis}
    \end{tikzpicture}\hfill
    \begin{tikzpicture}
        \begin{axis}[ymode=log,
                    xlabel={$N$},
                    xmax=20,
                    ymin=1e-6,
                    grid=both,
				ylabel={Eigenvalue},
                    ytick={1,0.01,0.0001,0.000001},
                    width=0.45\textwidth,
                    style={font=\footnotesize}]
            \addplot[red] table [x=Number, y=Eul_eig_norm, col sep=comma] {Dat025.csv};
            \addlegendentry{Eul POD}
            \addplot[blue] table [x=Number, y=ALE_eig_norm, col sep=comma] {Dat026.csv};
            \addlegendentry{ALE POD}
        \end{axis}
    \end{tikzpicture}   
    \vspace{-3mm}
    \caption{DMR: Eigenvalue decay of the PODs (normalized to have $\lambda_1=1$), non parametric on the left, parametric on the right}
    \label{fig:eig_DMR_param}
\end{figure}
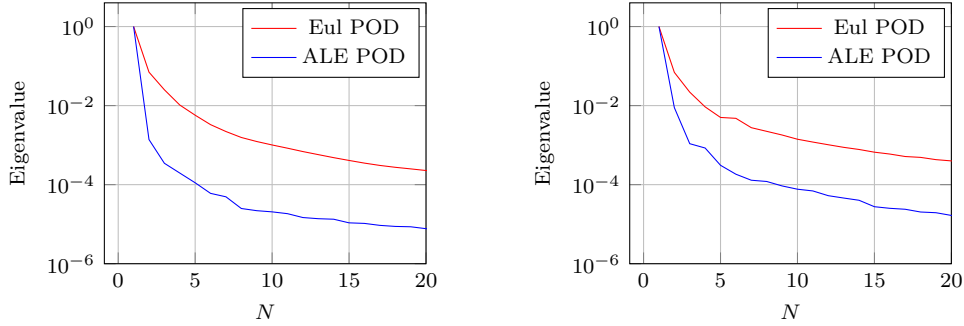
\begin{figure}
    \centering
    \begin{minipage}{0.045\textwidth}
        \includegraphics[width=\textwidth,trim={850 0 115 0},clip]{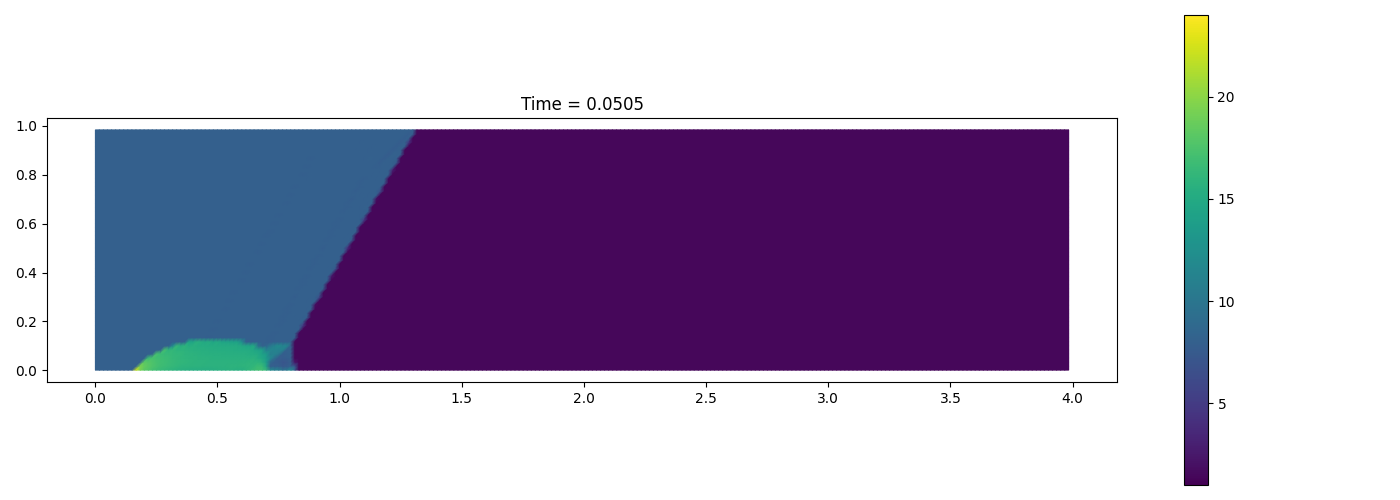}
    \end{minipage}
    \begin{minipage}{0.465\textwidth}
        \includegraphics[width=\textwidth,trim={0 0 280 25},clip]{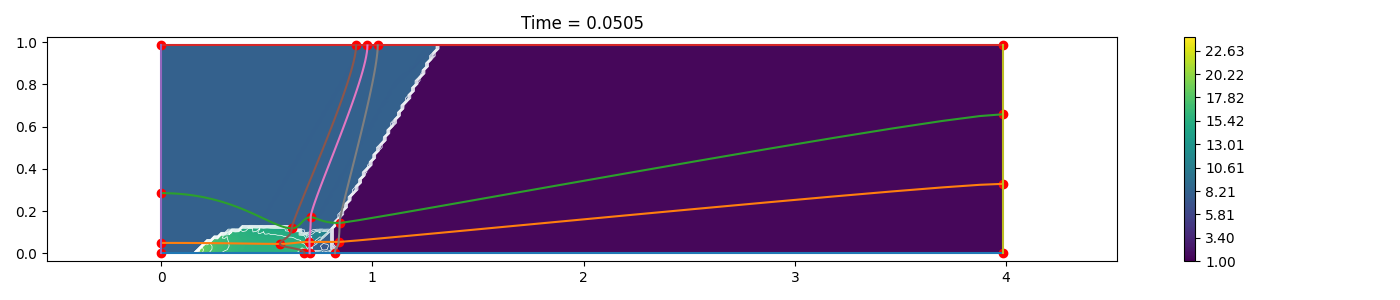}\\        
        \includegraphics[width=\textwidth,trim={0 0 280 25},clip]{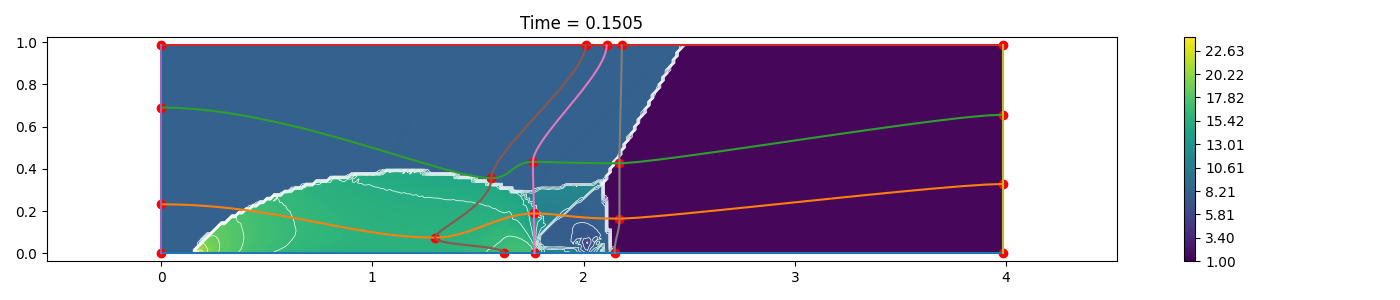}\\        
        \includegraphics[width=\textwidth,trim={0 0 280 25},clip]{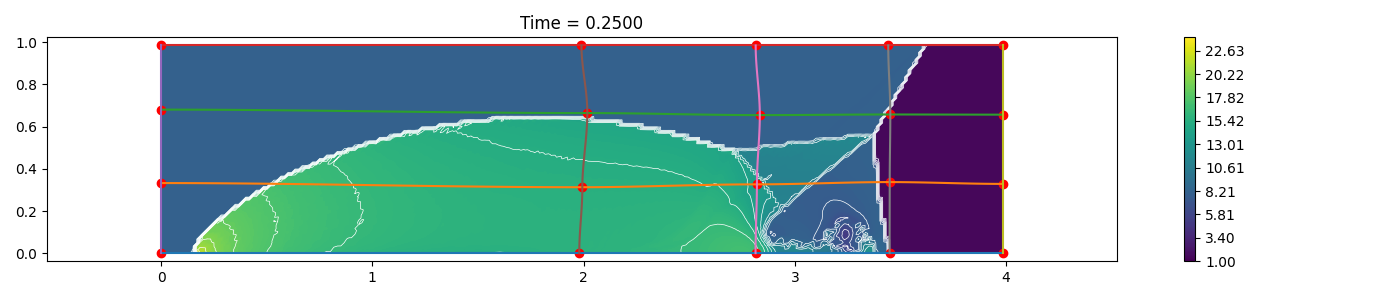}
    \end{minipage}
    \begin{minipage}{0.465\textwidth}
        \includegraphics[width=\textwidth,trim={0 0 280 25},clip]{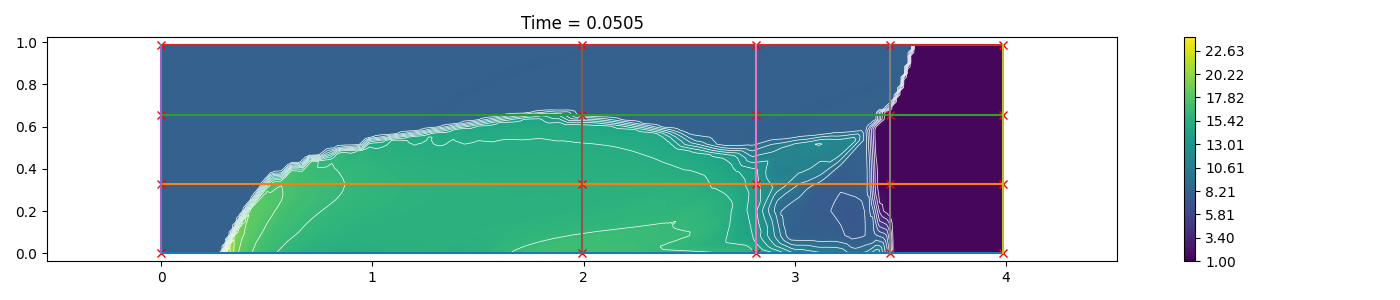}\\        
        \includegraphics[width=\textwidth,trim={0 0 280 25},clip]{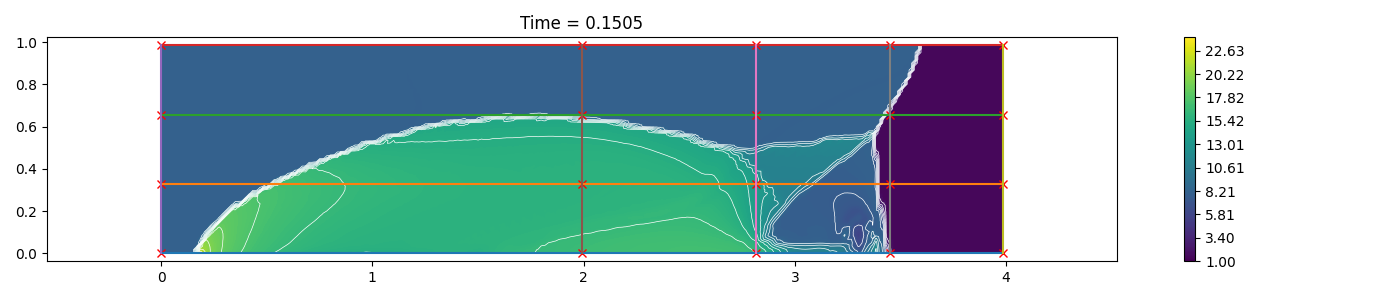}\\        
        \includegraphics[width=\textwidth,trim={0 0 280 25},clip]{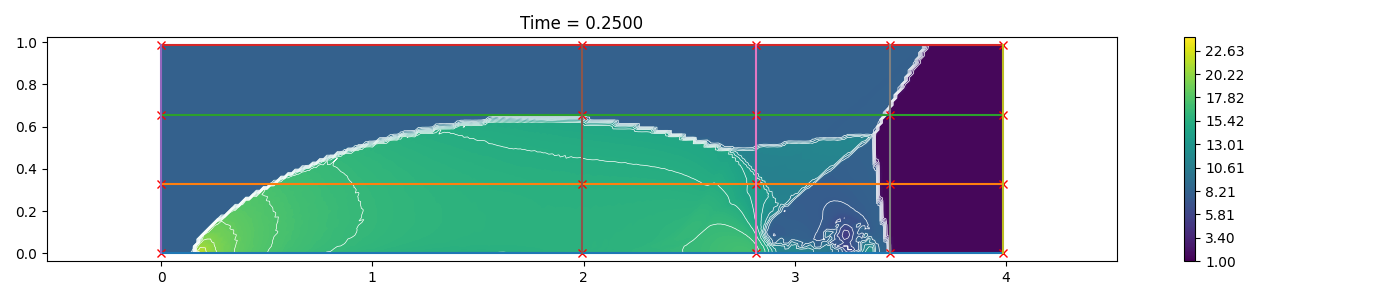}
    \end{minipage}
    \caption{FOM solution for $\rho$ for DMR non parametric at times 0.05 (top), 0.15 (center) and 0.25 (bottom) in the physical domain $\Omega$ (left) and, after calibration, in the reference configuration $\Omegaref$ (right). We mark on the plots the control points and the Cartesian grid that links them in the reference domain and its image through $T$ on the physical one. {We plot in white 20 contour lines at equispaced values between 1 and 25}}
    \label{fig:DMR_no_param_FOM_solutions}
\end{figure}

Fig.~\ref{fig:DMR_no_param_FOM_solutions} (left column) represents the FOM snapshots for the density $\rho$ at three different times of the numerical simulation; here, the same FV scheme has been employed at the FOM level on a mesh of $2400\times 600$ cells (computational time of 5 days), then downsampled to a mesh of $240\times 60$ cells to perform the offline phase (including calibration) in reasonable computational times. We retain 500 time samples in $[0,0.25]$ of which we include {100 in the training set of the calibration procedure (every $\Delta t =0.0025$) and 45 in the training set of the reduced algorithms all in $[0.02,0.2]$}.

\begin{figure}
    \centering
    \newcommand{\whicherror}{L2_errors}
    \begin{tikzpicture}
        \begin{axis}[ymode=log,
                    xlabel={Time},
				ylabel={Error},
				width=.68\textwidth,
				height=.45\textwidth,
                    legend pos=outer north east,
                    grid=both,
                    xmin = -0.01,
                    xmax= 0.26,
                    ymax=1,
                    style={font=\footnotesize}]
            \addplot[mark=square,mark size=1.3pt,red] table [x=Time, y=L2Error_eul_FOM_phy, col sep=comma] {Dat035.csv};
            \addlegendentry{Eul POD $N=2$}
            \addplot[mark=triangle*,mark size=1.3pt,red] table [x=Time, y=L2Error_eul_FOM_phy, col sep=comma] {Dat036.csv};
            \addlegendentry{Eul POD $N=4$}            
            \addplot[mark=diamond,mark size=1.3pt,red] table [x=Time, y=L2Error_eul_FOM_phy, col sep=comma] {Dat037.csv};
            \addlegendentry{Eul POD $N=7$}           
            \addplot[mark=otimes*,mark size=1.3pt,red] table [x=Time, y=L2Error_eul_FOM_phy, col sep=comma] {Dat038.csv};
            \addlegendentry{Eul POD $N=12$}
            \addplot[mark=pentagon,mark size=1.3pt,red] table [x=Time, y=L2Error_eul_FOM_phy, col sep=comma] {Dat039.csv};
            \addlegendentry{Eul POD $N=30$}

            \addplot[mark=square,mark size=1.3pt,blue] table [x=Time, y=L2Error_ALE_FOM_phy, col sep=comma] {Dat040.csv};
            \addlegendentry{ALE POD $N=2$}
            \addplot[mark=triangle*,mark size=1.3pt,blue] table [x=Time, y=L2Error_ALE_FOM_phy, col sep=comma] {Dat041.csv};
            \addlegendentry{ALE POD $N=4$}
            \addplot[mark=diamond,mark size=1.3pt,blue] table [x=Time, y=L2Error_ALE_FOM_phy, col sep=comma] {Dat042.csv};
            \addlegendentry{ALE POD $N=7$}
            \addplot[mark=otimes*,mark size=1.3pt,blue] table [x=Time, y=L2Error_ALE_FOM_phy, col sep=comma] {Dat043.csv};
            \addlegendentry{ALE POD $N=12$}
            \addplot[mark=pentagon,mark size=1.3pt,blue] table [x=Time, y=L2Error_ALE_FOM_phy, col sep=comma] {Dat044.csv};
            \addlegendentry{ALE POD $N=30$}
        \end{axis}
    \end{tikzpicture}    
    \vspace{-1mm}
    \caption{DMR non parametric: Error in time of reduced methods with different number $N$ of modes}
    \label{fig:error_DMR2D}
\end{figure}
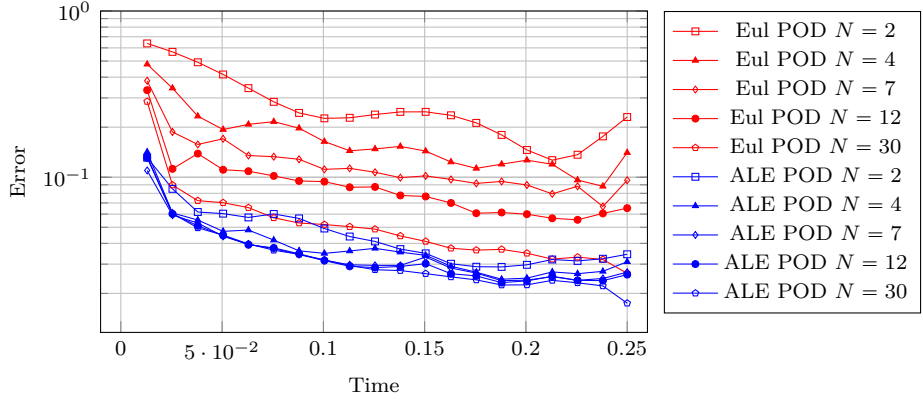

Fig.~\ref{fig:eig_DMR_param} (left) shows the rate of decay of the eigenvalues returned by the POD on $\mathcal{M}_{\rho}$ (red): also for this test case we have a solution manifold with a slowly decaying Kolmogorov $n$--width, due to the fact that the shock moves inside $\Omega$. 
We therefore perform a calibration procedure, using Algorithm~\ref{alg:calibration self similar} and the 2D geometrical transformation map $T$ introduced in Section \ref{sec:geometrical transformation 2D}: all the details for the calibration step are summarized in Table~\ref{tab: calibration DMR non-parametric}. 
Fig.~\ref{fig:DMR_no_param_FOM_solutions} (right column) shows the outcome of the calibration for the density $\rho$, at different times: here the snapshots are represented in the reference configuration $\Omegaref$ (computational time for forecasting the calibration points around 0.05s). With the calibration procedure, we obtain an improvement in the rate of decay of the eigenvalues, as it is shown in Figure \ref{fig:eig_DMR_param} (blue line). 
To conclude, in Fig.~\ref{fig:error_DMR2D} we show the behavior in time of the approximation error, between the FOM solution and the online solution (computational time to evaluate the NN for the ROM coefficients 0.001s), with or without calibration, according to the number $N$ of modes used. 
Both errors have been computed in the physical domain $\Omega$ and are defined as:
\begin{align}\label{eq:error_definition}
\frac{\lvert\lvert \rho(t) - \rho_N(t)\rvert\rvert_{L^2(\Omega)}}{\lvert\lvert \rho(t)\rvert\rvert_{L^2(\Omega)}}\quad \text{ and }  \quad \frac{\lvert\lvert \rho(t) - \rhohat_N(t)\circ T^{-1}[\bbw^{\text{opt}}(t)]\rvert\rvert_{L^2(\Omega)}}{\lvert\lvert \rho(t)\rvert\rvert_{L^2(\Omega)}}.
\end{align}
As we see in Fig.~\ref{fig:error_DMR2D}, the error of the online approximation (with calibration) does not go below a certain lower bound, even increasing the number of bases. 
We recall that, in the calibrated setting, we are interpolating the solutions to perform the transformations: for this reason, we believe that, after a certain number $N$ of modes, the interpolation error dominates the global error, and this leads to the plateau that one can observe in the figure. 
We notice that with around 30 basis for the POD in the Eulerian framework we achieve errors that are comparable with the ALE solutions.

\begin{figure}
    \centering
    \begin{minipage}{0.045\textwidth}
        \includegraphics[width=\textwidth,trim={850 0 115 0},clip]{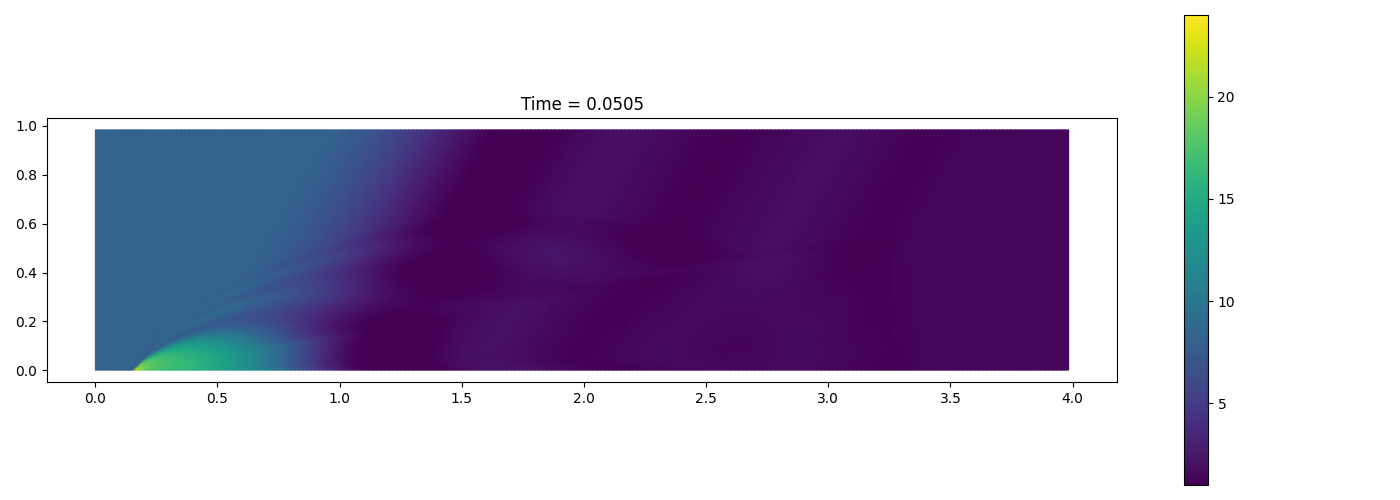}
    \end{minipage}
    \begin{minipage}{0.465\textwidth}
    \centering
    ALE $N=2$\\
        \includegraphics[width=\textwidth,trim={0 0 280 25},clip]{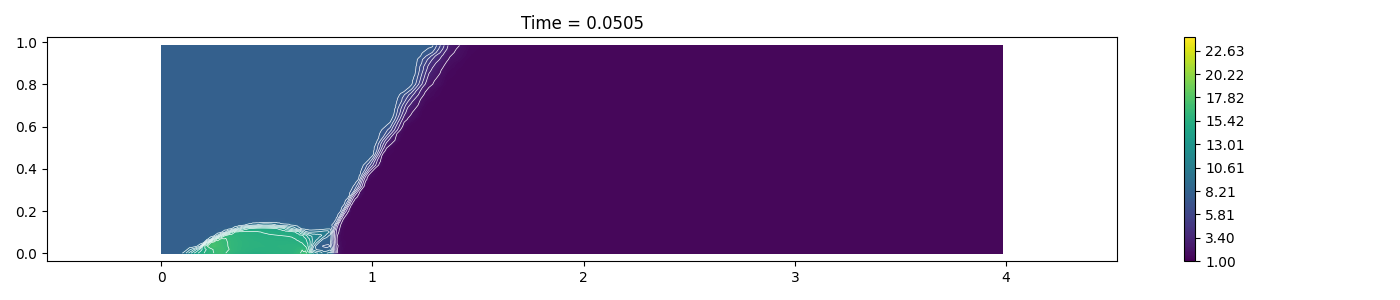}\\        
        \includegraphics[width=\textwidth,trim={0 0 280 25},clip]{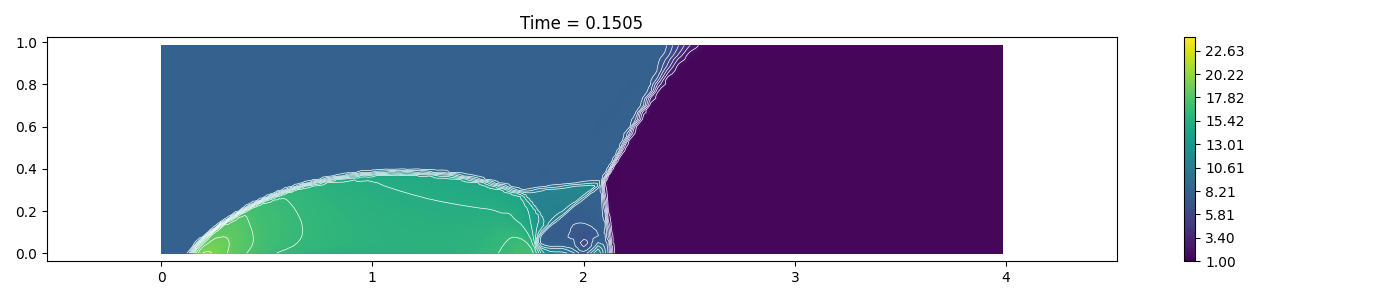}\\        
        \includegraphics[width=\textwidth,trim={0 0 280 25},clip]{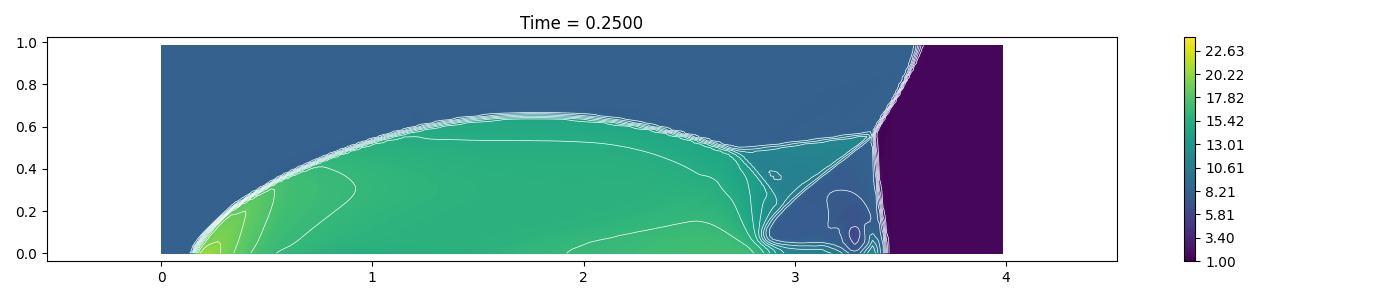}
    \end{minipage}
    \begin{minipage}{0.465\textwidth}
    \centering
    Eulerian $N=30$\\
        \includegraphics[width=\textwidth,trim={0 0 280 25},clip]{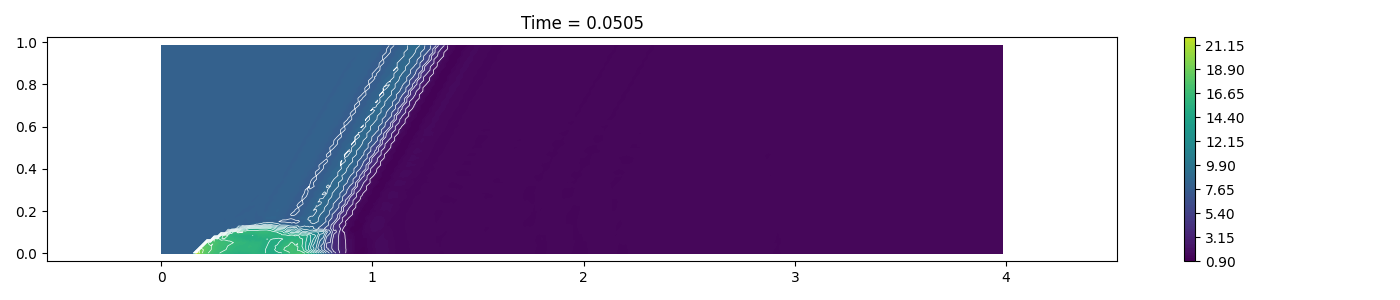}\\        
        \includegraphics[width=\textwidth,trim={0 0 280 25},clip]{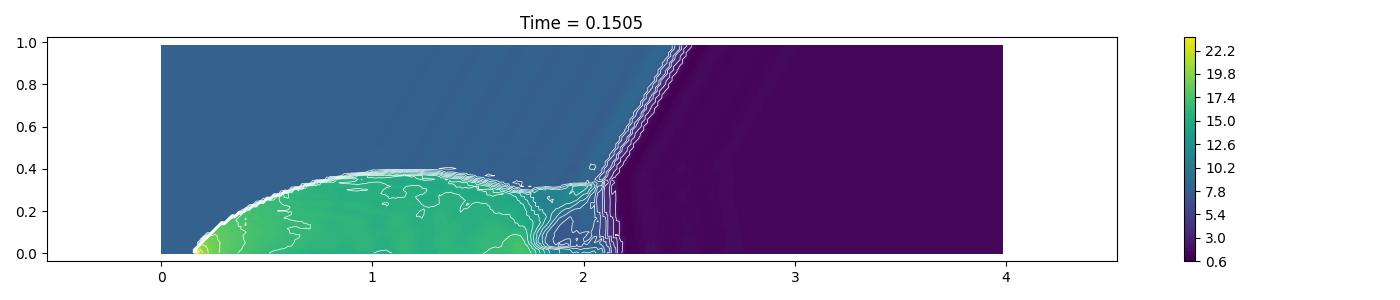}\\        
        \includegraphics[width=\textwidth,trim={0 0 280 25},clip]{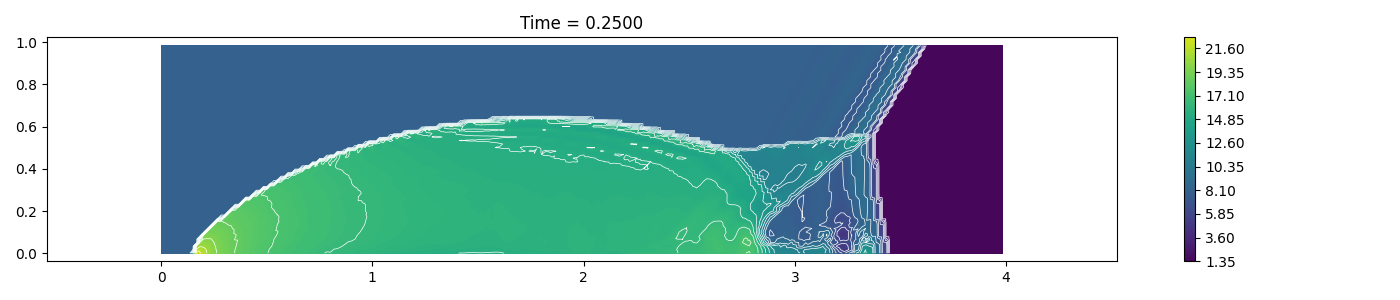}
    \end{minipage}
    \caption{DMR non parametric: ROM solutions for $\rho$ in $\Omega$ at times 0.05 (top), 0.15 (center) and 0.25 (bottom). Left column: ALE ROM solution with $N=2$. Right column: Eulerian ROM solution with $N=30$. {We plot in white 20 contour lines at equispaced values between 1 and 25}}
    \label{fig:DMR_no_param_ROM_solutions}
\end{figure}
Nevertheless, the qualitatively comparison of the Eulerian and ALE approach at the ROM level with $N=2$ for the ALE approach and $N=30$ for the Eulerian approach depicted in Fig.~\ref{fig:DMR_no_param_ROM_solutions} is still in favor of the ALE approach.
Indeed, the Eulerian ROM shows an oscillatory behavior that deteriorates the shape of the solution, the shock position and the flat areas, which are not anymore flat. 
On the contrary, the ALE ROM solutions are very similar to the FOM ones and they preserve all the original features even with a much smaller reduced basis.
So, even if the $L^2$ errors of the two approaches are comparable, the quality of the two solutions is very different.

\subsection{Parametric DMR problem}\label{sec:numerical results DMR parametric}
In this section, we consider the parametric version of the 2D DMR problem: the physical parameter is the angle $\beta$ introduced in Section~\ref{sec:numerical results DMR non-parametric}; the physical parameter interval is $\mathcal{P}_{\text{phys}}=[0.1,0.675]$, and the time interval is $[0, 0.2]$. Also in this case, the FOM snapshots have been obtained with the same FV scheme, on a mesh $600\times 150$ (computational time of 2 hours each) then downsampled to $200\times 50$ for reduction of computational time of the offline phase. {In the training set, we include for each physical parameter 51 snapshots every $\Delta t=0.004$, and we use all the snapshots for the calibration.}
\begin{figure}
    \centering
    Parameter $\beta =0.225$\\
    \begin{minipage}{0.045\textwidth}
        \includegraphics[width=\textwidth,trim={850 0 115 0},clip]{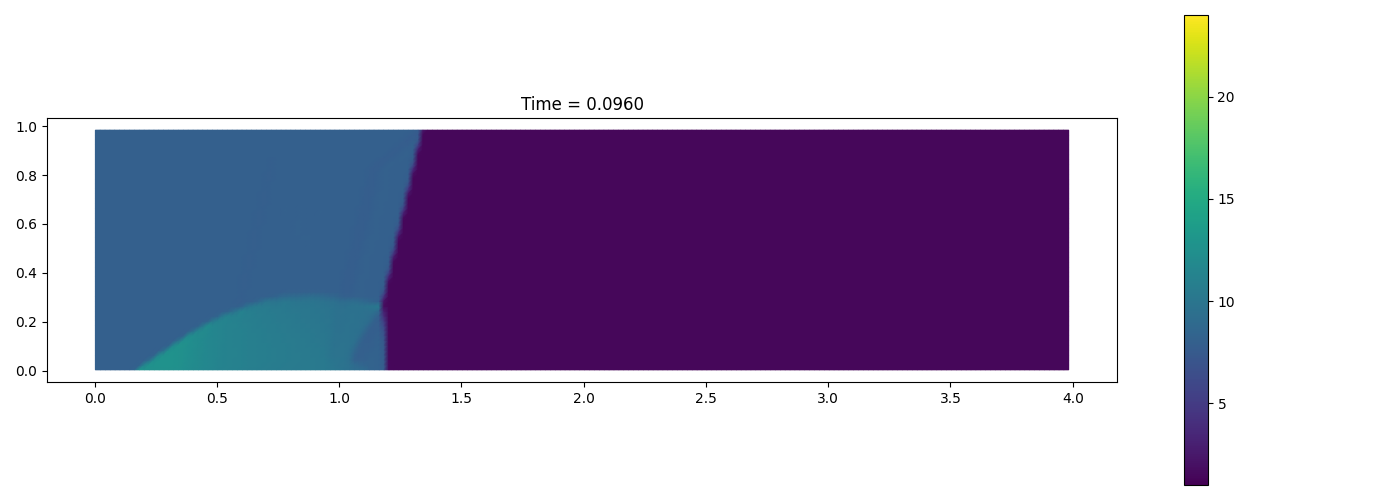}
    \end{minipage}
    \begin{minipage}{0.465\textwidth}
        \includegraphics[width=\textwidth,trim={0 0 280 25},clip]{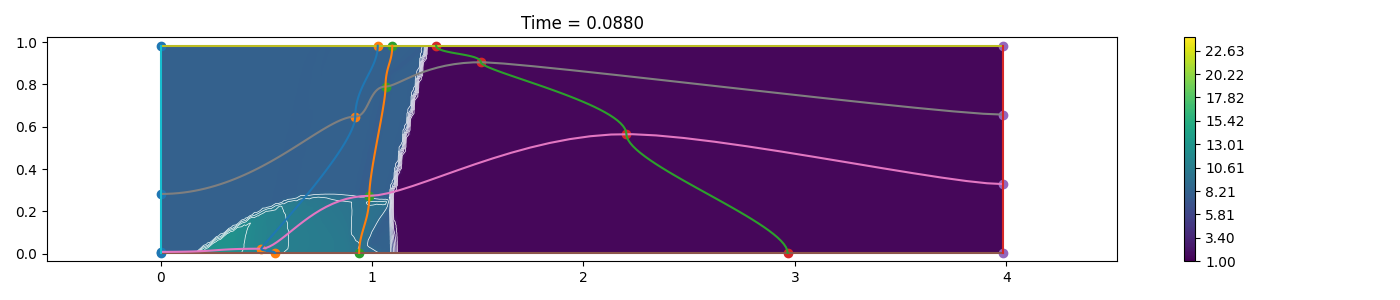}\\        
        \includegraphics[width=\textwidth,trim={0 0 280 25},clip]{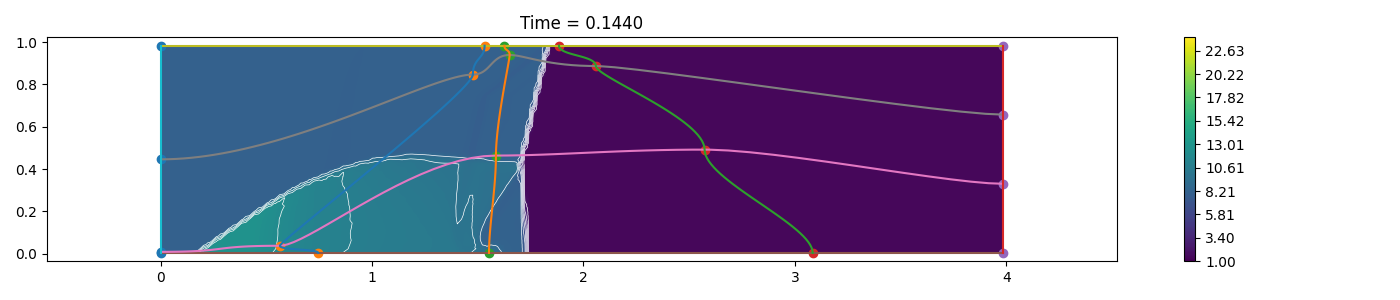}\\        
        \includegraphics[width=\textwidth,trim={0 0 280 25},clip]{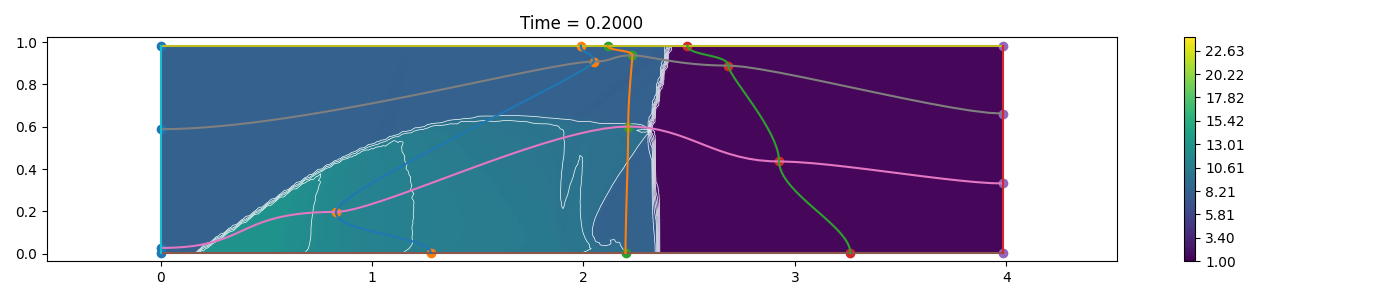}
    \end{minipage}
    \begin{minipage}{0.465\textwidth}
        \includegraphics[width=\textwidth,trim={0 0 280 25},clip]{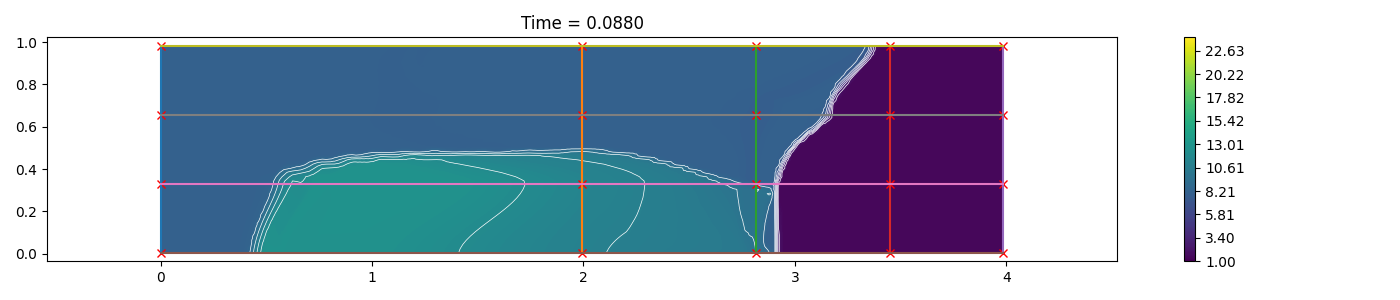}\\        
        \includegraphics[width=\textwidth,trim={0 0 280 25},clip]{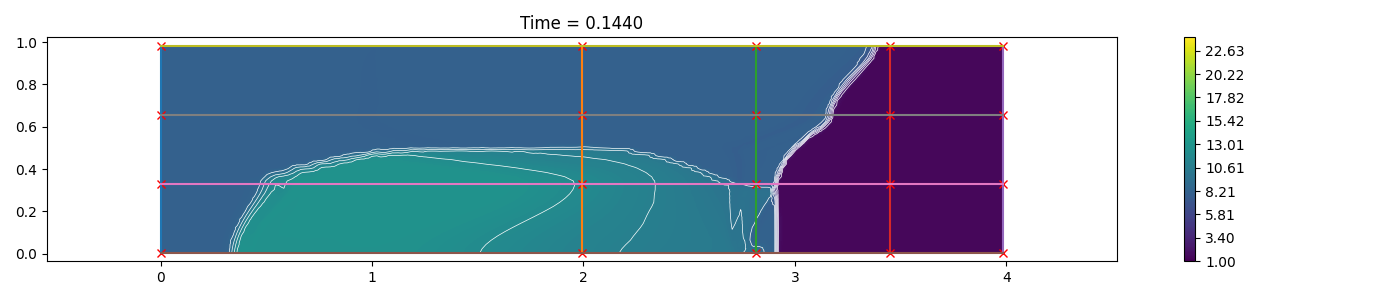}\\        
        \includegraphics[width=\textwidth,trim={0 0 280 25},clip]{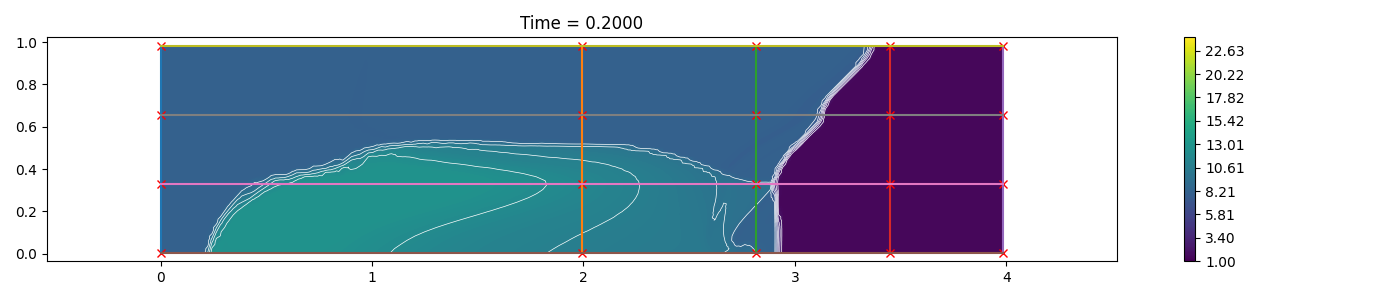}
    \end{minipage}\\[2mm]
    Parameter $\beta =0.675$\\
    \begin{minipage}{0.045\textwidth}
        \includegraphics[width=\textwidth,trim={850 0 115 0},clip]{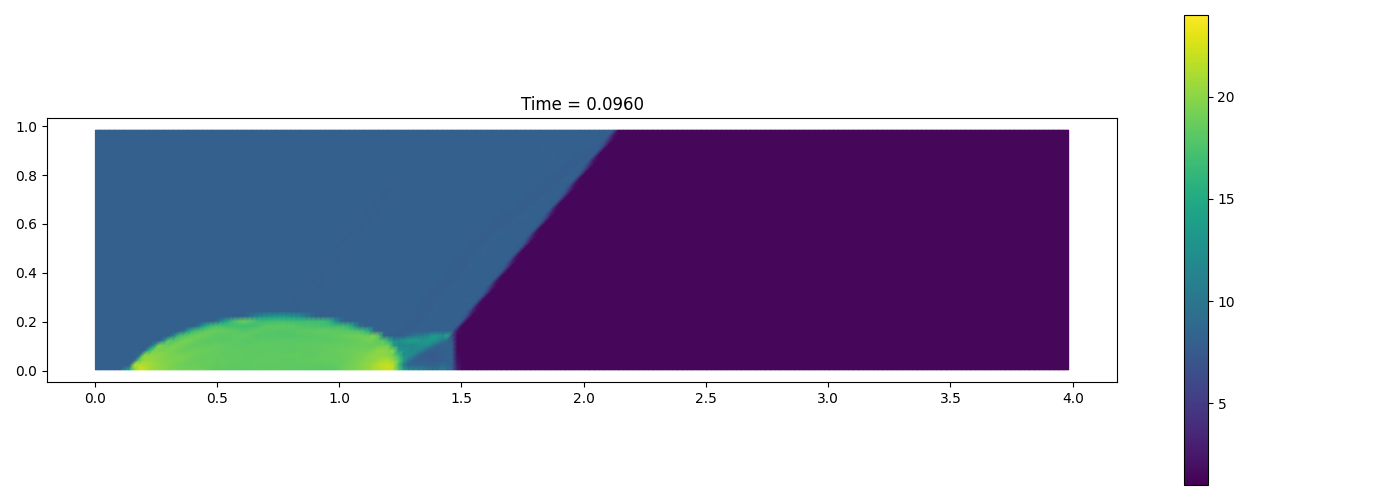}
    \end{minipage}
    \begin{minipage}{0.465\textwidth}
        \includegraphics[width=\textwidth,trim={0 0 280 25},clip]{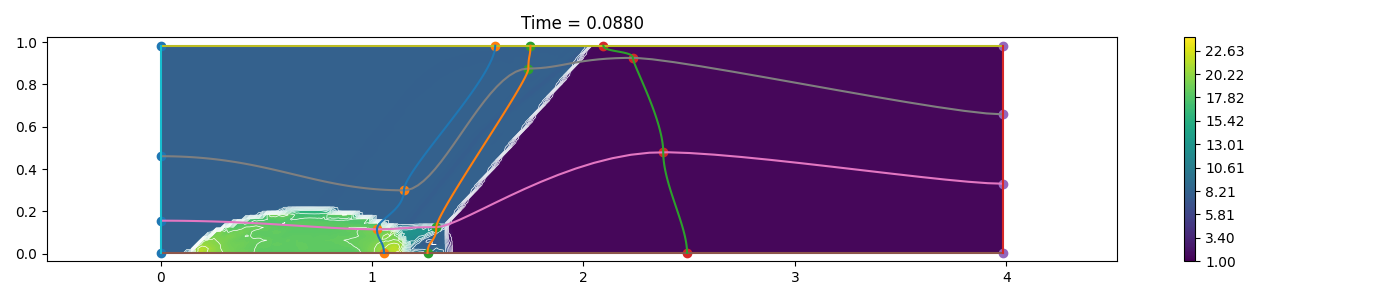}\\        
        \includegraphics[width=\textwidth,trim={0 0 280 25},clip]{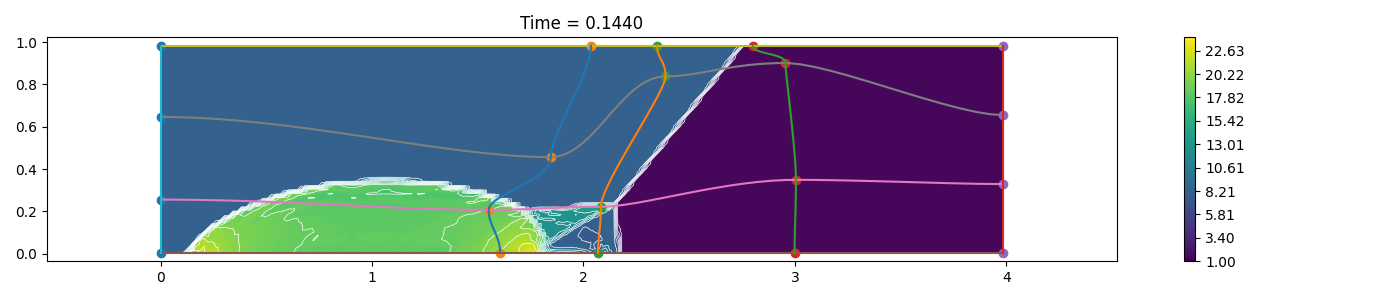}\\        
        \includegraphics[width=\textwidth,trim={0 0 280 25},clip]{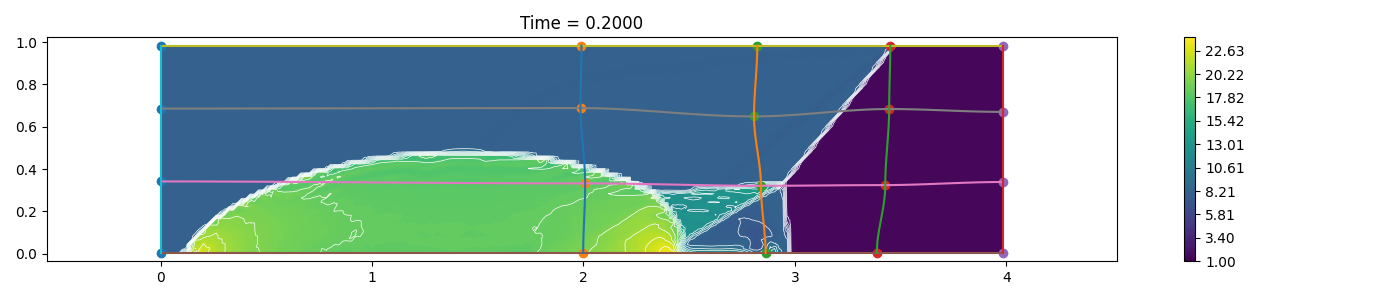}
    \end{minipage}
    \begin{minipage}{0.465\textwidth}
        \includegraphics[width=\textwidth,trim={0 0 280 25},clip]{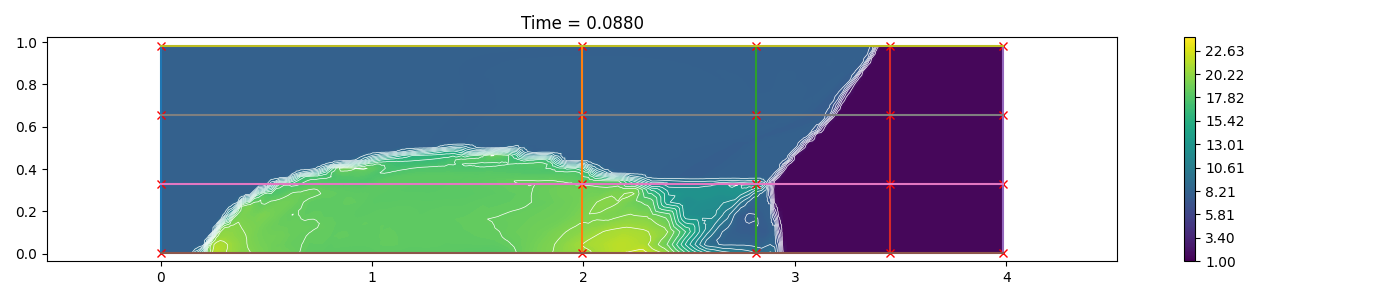}\\        
        \includegraphics[width=\textwidth,trim={0 0 280 25},clip]{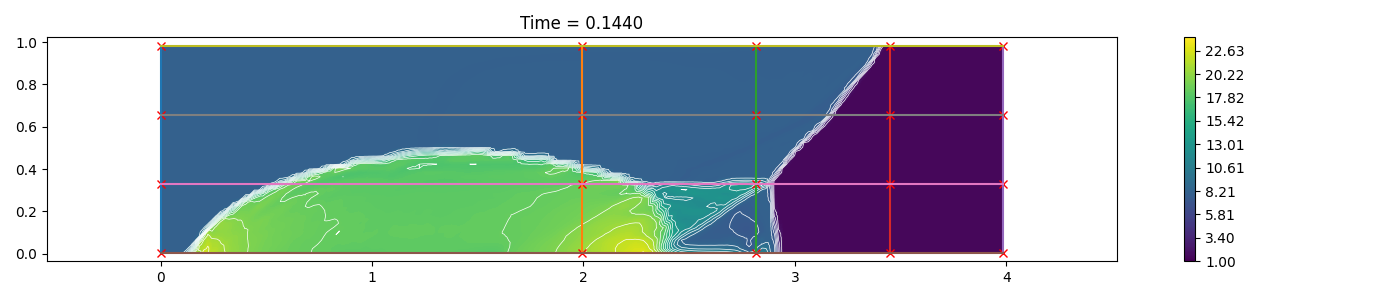}\\        
        \includegraphics[width=\textwidth,trim={0 0 280 25},clip]{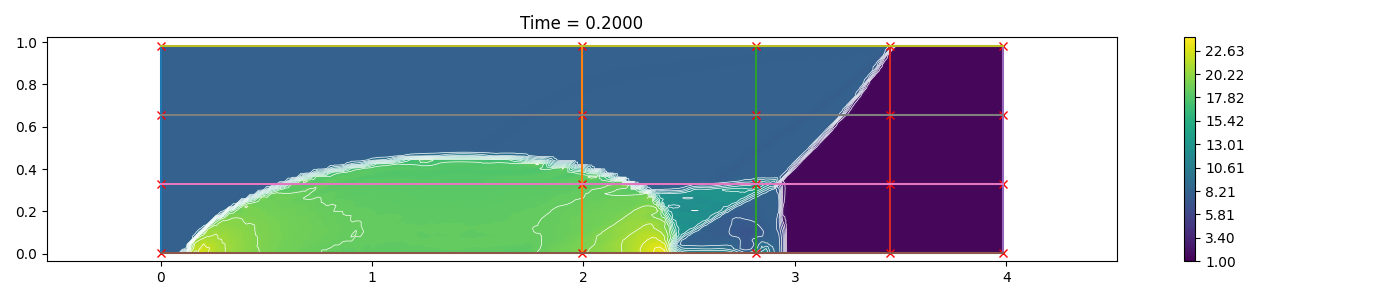}
    \end{minipage}    
    \caption{DMR parametric FOM solution for $\rho$ at times $t=0.096$, $t=0.148$ and $t=0.2$ (top to bottom) in the physical domain $\Omega$ (left) and, after calibration, in the reference configuration $\Omegaref$ (right). We mark on the plots the control points and the Cartesian grid that links them. Parameter values $\beta=0.225$ (top) and $\beta=0.675$ (bottom). {We plot in white 20 contour lines at equispaced values between 1 and 25}}
    \label{fig:DMR_param_FOM_solutions}
\end{figure}

\begin{figure}
    \centering
    \newcommand{\testfolder}{img/DMR2D_param/ROM/param_675}
    \newcommand{\whicherror}{proj_errors}
    \begin{tikzpicture}
        \begin{axis}[ymode=log,
                    xlabel={Time},
				title={Projection Error},
				width=.45\textwidth,
				height=.4\textwidth,
                    grid=both,
                    xmin = -0.01,
                    xmax= 0.21,
                    ymax=1,
                    ymin=0.015,
                    style={font=\footnotesize}]
            \addplot[mark=square,mark size=1.3pt,red] table [x=Time, y=L2Error_eul_FOM_phy, col sep=comma] {Dat045.csv};
            \addplot[mark=triangle*,mark size=1.3pt,red] table [x=Time, y=L2Error_eul_FOM_phy, col sep=comma] {Dat046.csv};
            \addplot[mark=diamond,mark size=1.3pt,red] table [x=Time, y=L2Error_eul_FOM_phy, col sep=comma] {Dat047.csv};
            \addplot[mark=otimes*,mark size=1.3pt,red] table [x=Time, y=L2Error_eul_FOM_phy, col sep=comma] {Dat048.csv};
            \addplot[mark=pentagon,mark size=1.3pt,red] table [x=Time, y=L2Error_eul_FOM_phy, col sep=comma] {Dat049.csv};

            \addplot[mark=square,mark size=1.3pt,blue] table [x=Time, y=L2Error_ALE_FOM_phy, col sep=comma] {Dat050.csv};
            \addplot[mark=triangle*,mark size=1.3pt,blue] table [x=Time, y=L2Error_ALE_FOM_phy, col sep=comma] {Dat051.csv};
            \addplot[mark=diamond,mark size=1.3pt,blue] table [x=Time, y=L2Error_ALE_FOM_phy, col sep=comma] {Dat052.csv};
            \addplot[mark=otimes*,mark size=1.3pt,blue] table [x=Time, y=L2Error_ALE_FOM_phy, col sep=comma] {Dat053.csv};
            \addplot[mark=pentagon,mark size=1.3pt,blue] table [x=Time, y=L2Error_ALE_FOM_phy, col sep=comma] {Dat054.csv};
        \end{axis}
    \end{tikzpicture}
    \renewcommand{\whicherror}{L2_errors}
    \begin{tikzpicture}
        \begin{axis}[ymode=log,
                    xlabel={Time},
				width=.45\textwidth,
				height=.4\textwidth,
                    legend pos=outer north east,
                    yticklabel=\empty,
                    grid=both,
                    xmin = -0.01,
                    xmax= 0.21,
                    ymax=1,
                    ymin=0.015,
                    title={POD-NN Error},
                    style={font=\footnotesize}]
            \addplot[mark=square,mark size=1.3pt,red] table [x=Time, y=L2Error_eul_FOM_phy, col sep=comma] {Dat055.csv};
            \addlegendentry{Eul $N=2$}
            \addplot[mark=triangle*,mark size=1.3pt,red] table [x=Time, y=L2Error_eul_FOM_phy, col sep=comma] {Dat056.csv};
            \addlegendentry{Eul $N=6$}            
            \addplot[mark=diamond,mark size=1.3pt,red] table [x=Time, y=L2Error_eul_FOM_phy, col sep=comma] {Dat057.csv};
            \addlegendentry{Eul $N=12$}           
            \addplot[mark=otimes*,mark size=1.3pt,red] table [x=Time, y=L2Error_eul_FOM_phy, col sep=comma] {Dat058.csv};
            \addlegendentry{Eul $N=20$}   
            \addplot[mark=pentagon,mark size=1.3pt,red] table [x=Time, y=L2Error_eul_FOM_phy, col sep=comma] {Dat059.csv};
            \addlegendentry{Eul $N=40$}

            \addplot[mark=square,mark size=1.3pt,blue] table [x=Time, y=L2Error_ALE_FOM_phy, col sep=comma] {Dat060.csv};
            \addlegendentry{ALE $N=2$}
            \addplot[mark=triangle*,mark size=1.3pt,blue] table [x=Time, y=L2Error_ALE_FOM_phy, col sep=comma] {Dat061.csv};
            \addlegendentry{ALE $N=6$}
            \addplot[mark=diamond,mark size=1.3pt,blue] table [x=Time, y=L2Error_ALE_FOM_phy, col sep=comma] {Dat062.csv};
            \addlegendentry{ALE $N=12$}
            \addplot[mark=otimes*,mark size=1.3pt,blue] table [x=Time, y=L2Error_ALE_FOM_phy, col sep=comma] {Dat063.csv};
            \addlegendentry{ALE $N=20$}
            \addplot[mark=pentagon,mark size=1.3pt,blue] table [x=Time, y=L2Error_ALE_FOM_phy, col sep=comma] {Dat064.csv};
            \addlegendentry{ALE $N=40$}
        \end{axis}
    \end{tikzpicture}
    \vspace{-1mm}
    \caption{DMR parametric: Error in time of reduced methods with different number $N$ of modes. Parameter in test set $\beta=0.675$}
    \label{fig:error_DMR2D_param}
\end{figure}
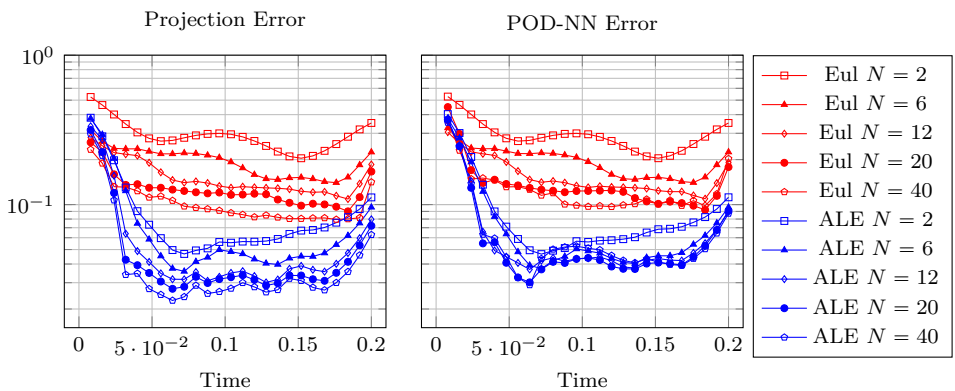
\begin{figure}
    \centering
    Parameter $\beta=0.225$\\
    \begin{minipage}{0.045\textwidth}
        \includegraphics[width=\textwidth,trim={850 0 115 0},clip]{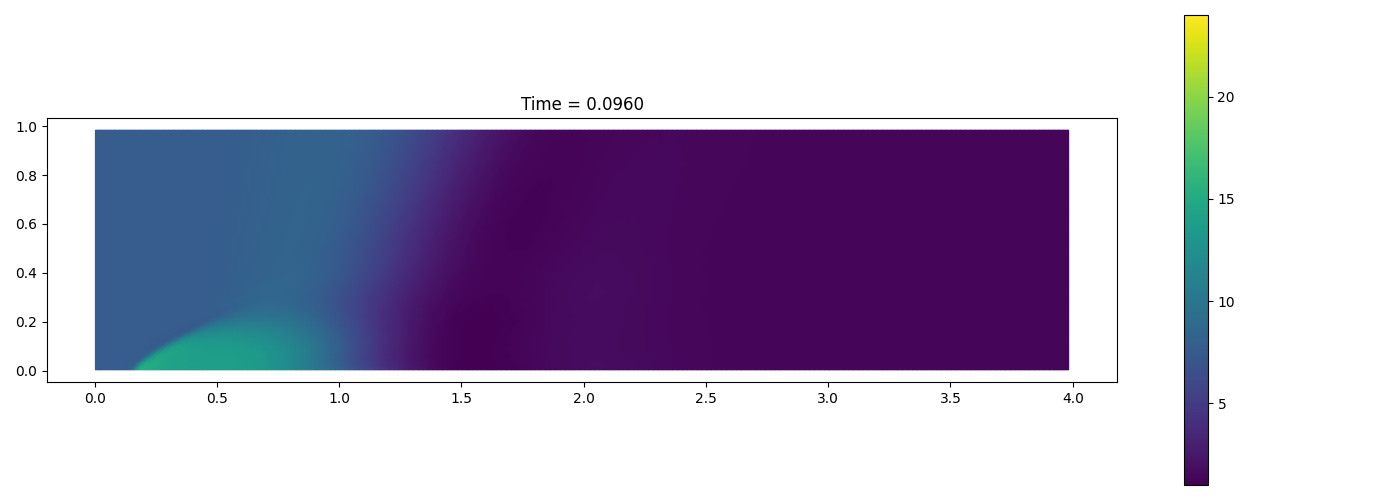}
    \end{minipage}
    \begin{minipage}{0.465\textwidth}
    \centering
        ALE POD-NN\\
        \includegraphics[width=\textwidth,trim={0 0 280 25},clip]{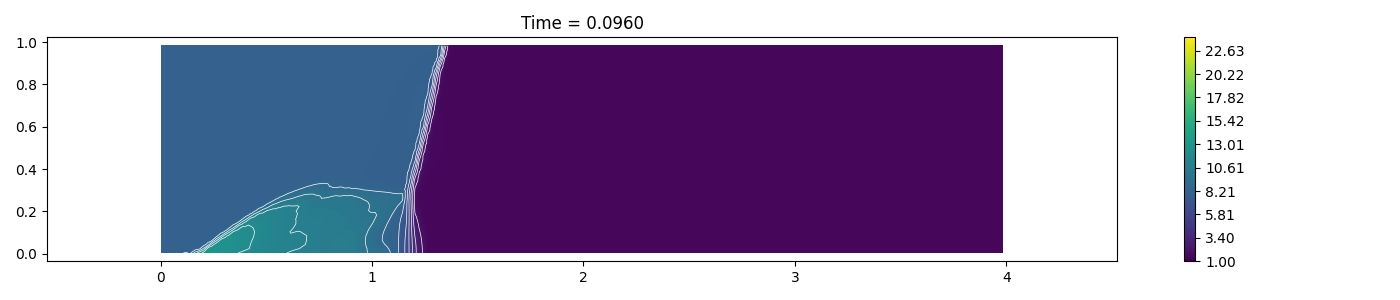}\\        
        \includegraphics[width=\textwidth,trim={0 0 280 25},clip]{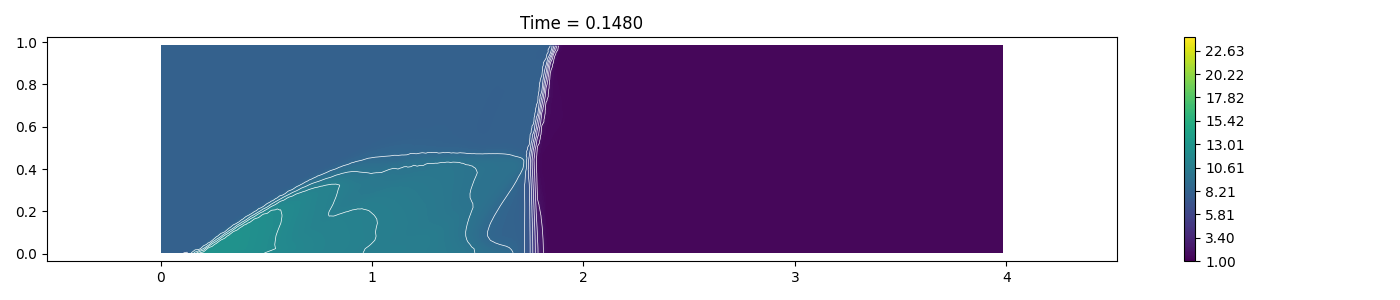}\\        
        \includegraphics[width=\textwidth,trim={0 0 280 25},clip]{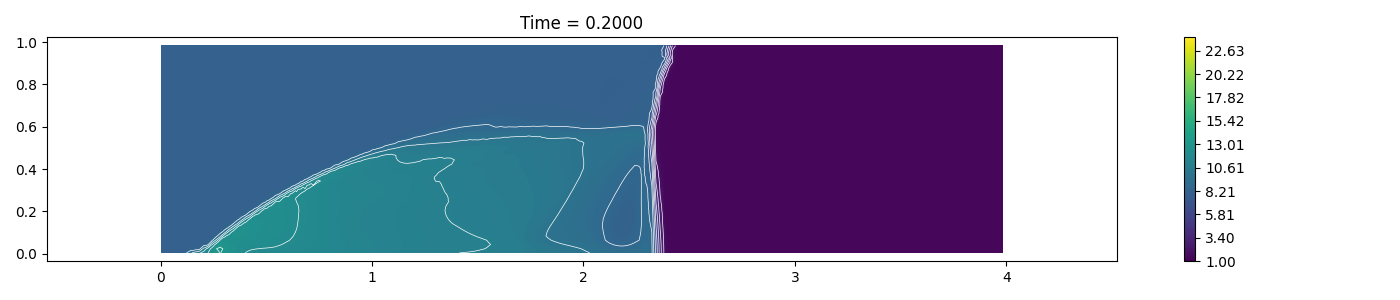}
    \end{minipage}
    \begin{minipage}{0.465\textwidth}
    \centering
        Eulerian POD-NN\\
        \includegraphics[width=\textwidth,trim={0 0 280 25},clip]{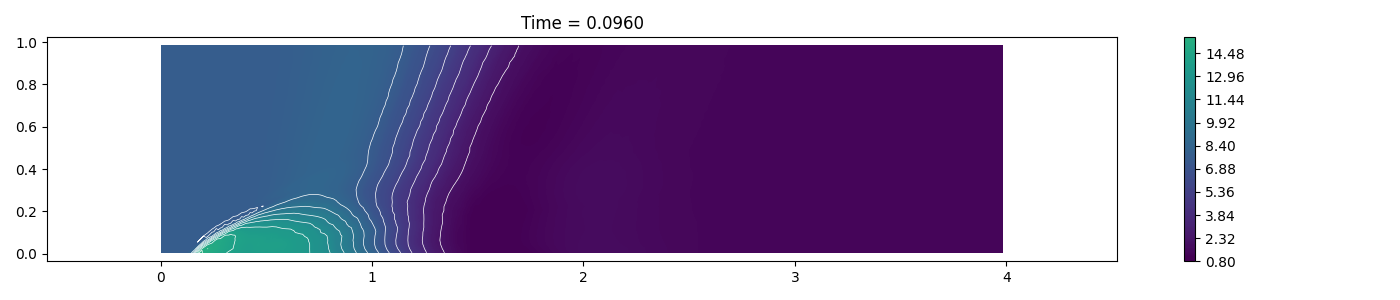}\\        
        \includegraphics[width=\textwidth,trim={0 0 280 25},clip]{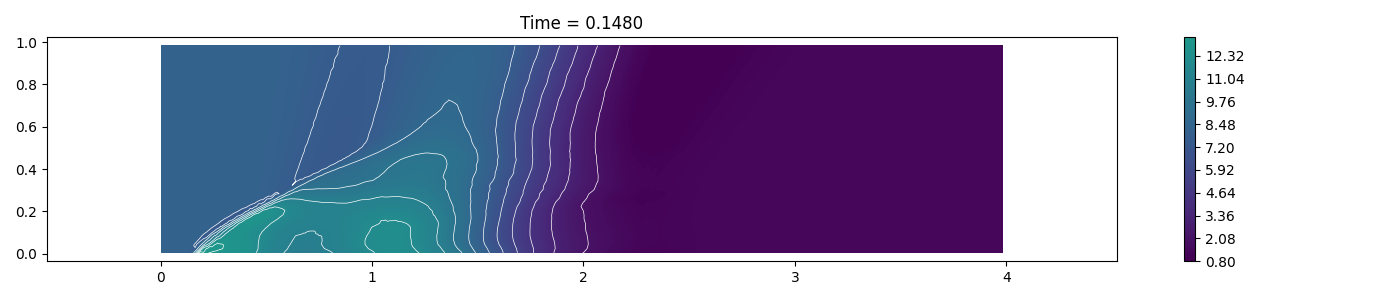}\\        
        \includegraphics[width=\textwidth,trim={0 0 280 25},clip]{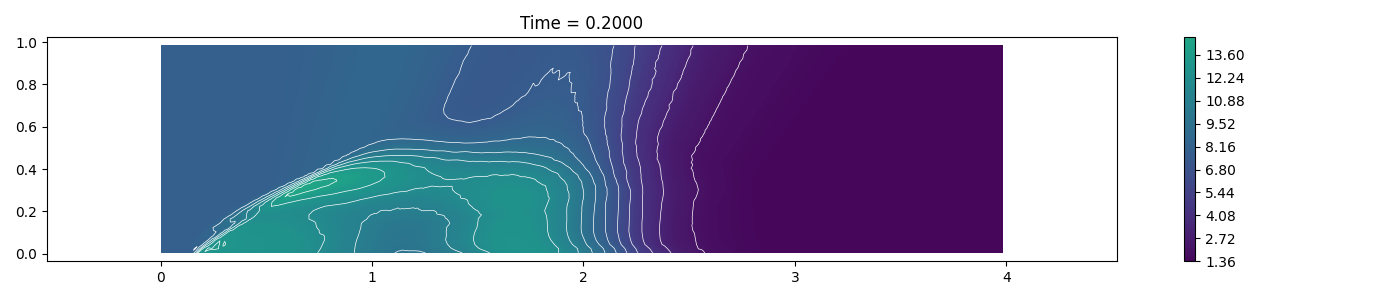}
    \end{minipage}\\[2mm]
    Parameter $\beta=0.675$\\
    \begin{minipage}{0.045\textwidth}
        \includegraphics[width=\textwidth,trim={850 0 115 0},clip]{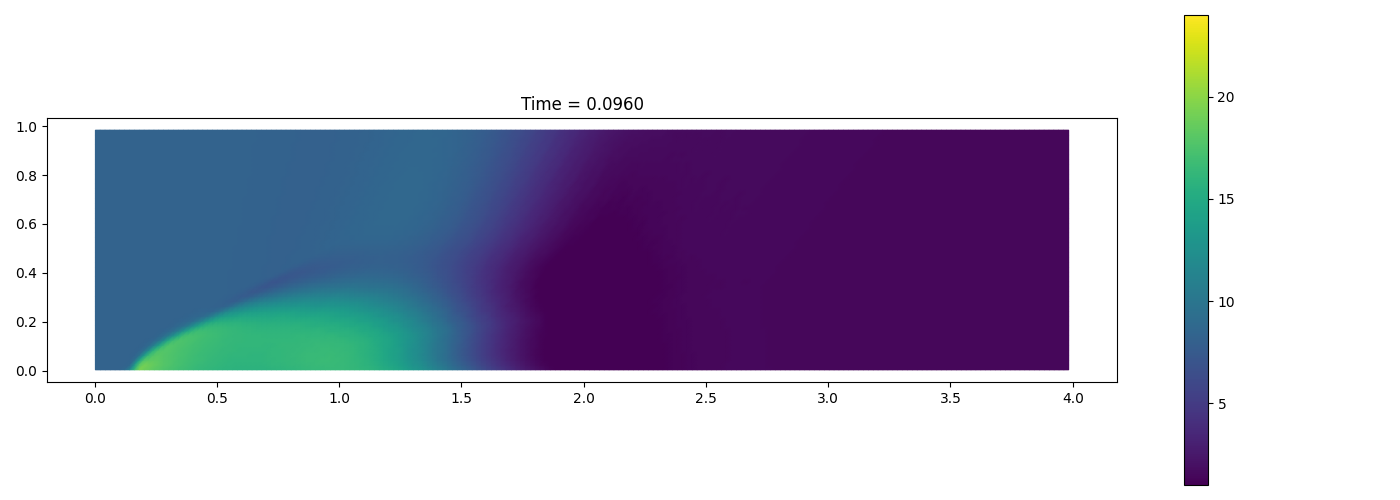}
    \end{minipage}
    \begin{minipage}{0.465\textwidth}
    \centering
        ALE POD-NN\\
        \includegraphics[width=\textwidth,trim={0 0 280 25},clip]{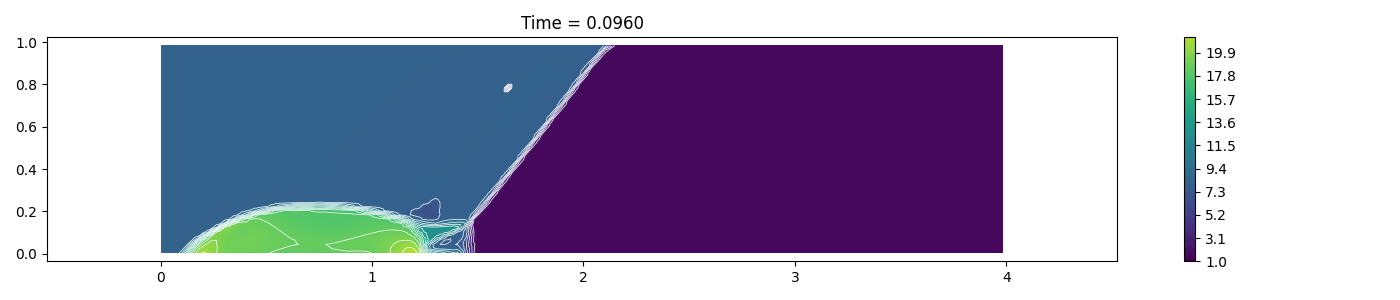}\\        
        \includegraphics[width=\textwidth,trim={0 0 280 25},clip]{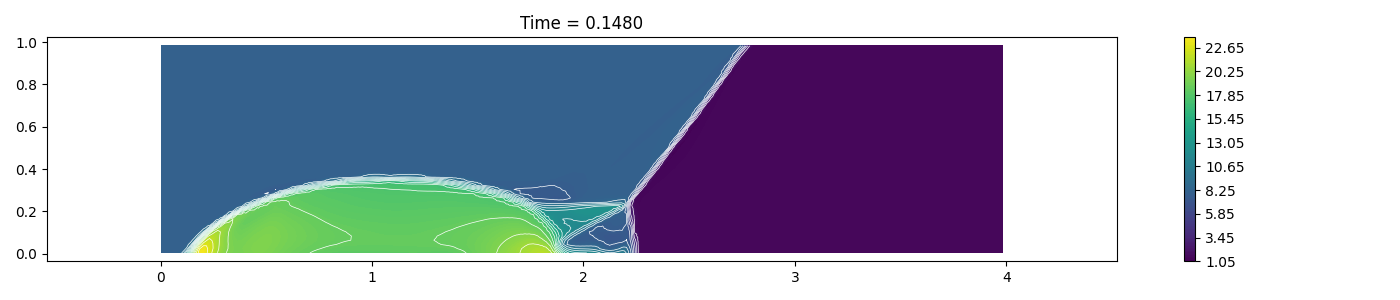}\\        
        \includegraphics[width=\textwidth,trim={0 0 280 25},clip]{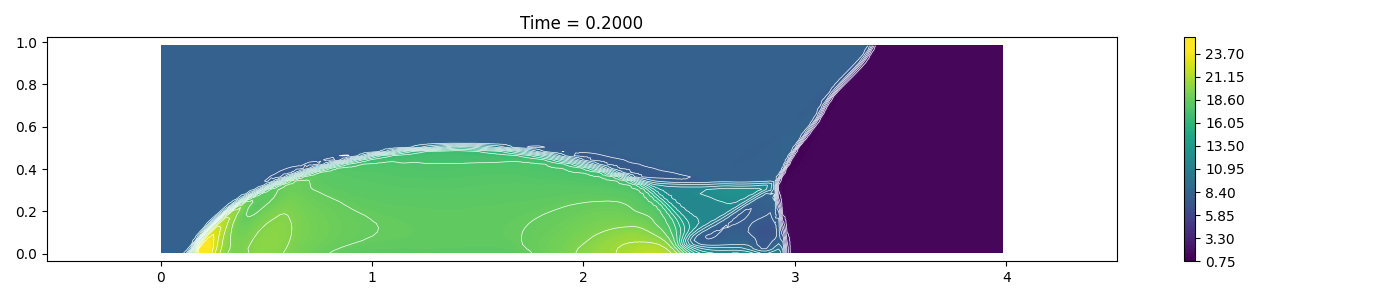}
    \end{minipage}
    \begin{minipage}{0.465\textwidth}
    \centering
        Eulerian POD-NN\\
        \includegraphics[width=\textwidth,trim={0 0 280 25},clip]{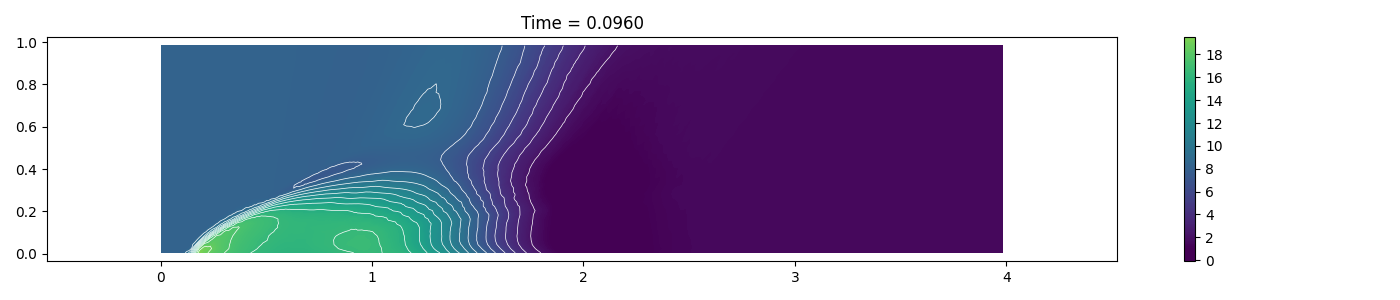}\\        
        \includegraphics[width=\textwidth,trim={0 0 280 25},clip]{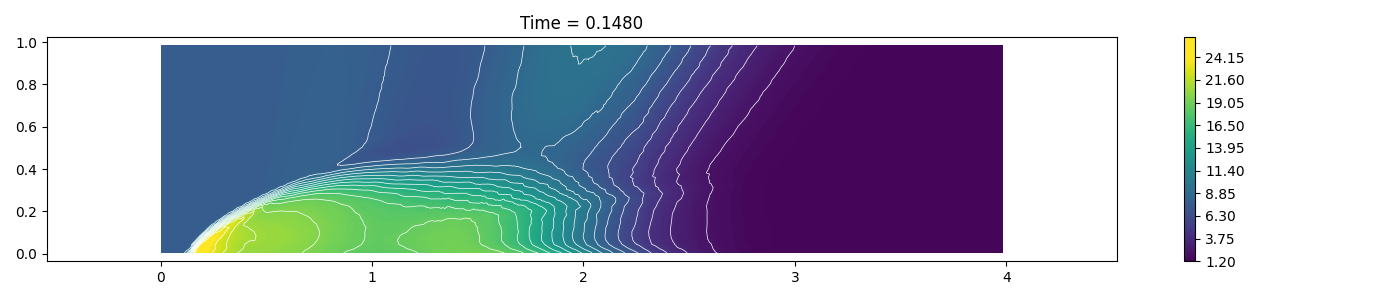}\\        
        \includegraphics[width=\textwidth,trim={0 0 280 25},clip]{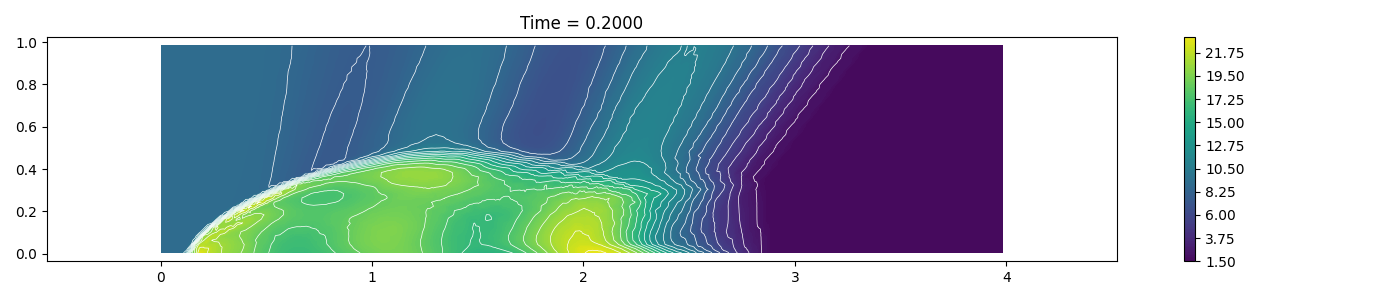}
    \end{minipage}
    \caption{POD-NN solutions with $N=6$ of $\rho$ for DMR parametric at times $t=0.096$, $t=0.148$ and $t=0.2$ (top to bottom) in the physical domain $\Omega$. Left column: with the calibration of the manifold. Right column: without calibration. Parameter values $\beta=0.225$ (top) and $\beta=0.675$ (bottom). {We plot in white 20 contour lines at equispaced values between 1 and 25}}
    \label{fig:DMR_param_ROM_solutions}
\end{figure}

Fig.~\ref{fig:DMR_param_FOM_solutions} shows some snapshots, for two different values of $\beta$ and at different times, before and after the calibration: all the details for the calibration procedure are summarized in Table~\ref{tab: calibration DMR non-parametric}. 
We also depicted the control points grid on the reference domain and its transformation onto the physical one, showing how the tracking of the interesting point is done and how much distortion we can get with such transformations.
As we can see from Fig.~\ref{fig:eig_DMR_param}, also in the parametric case the calibration procedure improves significantly the rate of the decay of the eigenvalues returned by the POD and hence, ultimately, the Kolmogorov $n$-width of the problem under consideration.
In Fig.~\ref{fig:error_DMR2D_param}, we plot the behavior of the relative error on the physical domain, as explained in Section~\ref{sec:numerical results DMR non-parametric} in~\eqref{eq:error_definition}, varying time and for different number of modes used in the reduced spaces. 
On the left, we plot the error between the FOM solution and the $L^2$ projection onto the reduced space; on the right, we have the error obtained using the POD-NN to predict the online solution.
We see that both the Eulerian and the ALE projection errors improve as we increment the number of POD modes, with the Eulerian being always much larger. In the POD-NN error, on the other side, the decay of the error is slower and it seems to stagnate at some bottleneck values, in particular for the ALE case. 
That is why we aim at extending this work in the future with a hyper-reduced Galerkin projection approach, to reintroduce some mathematical rigorousness hoping to decrease the online error.
Finally, in Fig.~\ref{fig:DMR_param_ROM_solutions} we represent the online solutions for $\beta=0.225$ and $\beta=0.675$, both with the Eulerian and the ALE approach with $N=6$. 
{Similarly to what happens in the non parametric test case, the Eulerian approach struggles to reproduce the FOM solution, providing an approximation that sometimes even loses the main features (the shape of the solution, the shocks, the flat areas). On the contrary, with the ALE approach, the online approximation preserves all these features. 
The two parameters shown validate the ability of this ROM approach to work in strongly nonlinear parametric context, where the parameters changes the solution's feature geometry, the values of the solution and vaguely the structure of the features. On the other hand, we remark that this approach works only for quasi-self-similar solutions, where we can recognize a similar structure along the parameter domain.} 

\subsection{Triple point non parametric}
{The applications of the proposed algorithm are various, we have showed some self-similar solutions, but this class includes a huge amount of test problems: for example, many hypersonic problems where a shock or multiple shocks are present in the solution as airfoil simulations, water waves propagation or acoustic waves. We want to solve another test here that involves a more complicated solution structure, which is a triple point shock interaction test. We consider a physical domain $\Omega=[0,7]\times[0,3]$. The initial conditions are:
\begin{equation}
    \begin{cases}
        (\rho_W , u_W, v_W, p_W) = (1,20,0,1)   & \bbx \in [0,1]\times[0,3],\\
        (\rho_{NE} , u_{NE}, v_{NE}, p_{NE}) = (0.125,0,0,0.1)     & \bbx \in [1,7]\times[1.5,3], \\
        (\rho_{SE} , u_{SE}, v_{SE}, p_{SE}) = (1,0,0,0.1)      & \bbx \in [1,7]\times[0,1.5].
    \end{cases}
\end{equation}
Boundary conditions are transmissive on the right, Dirichlet with state $(\rho_W , u_W, v_W, p_W)$ on the left and reflective at the top and bottom of the domain. Final time is set to $t_f=0.25$. We solve the problem on a grid $2800\times 1200$ and then we downsample it to $280\times 120$ for the reduction of computational time in the offline phase.
We have used in the training set of both calibration and reduced algorithms $N_\mu=100$ snapshots at $\Delta t=0.0025$. 

\begin{figure}
    \centering
    \begin{minipage}{0.045\textwidth}
        \includegraphics[width=\textwidth,trim={850 0 110 0},clip]{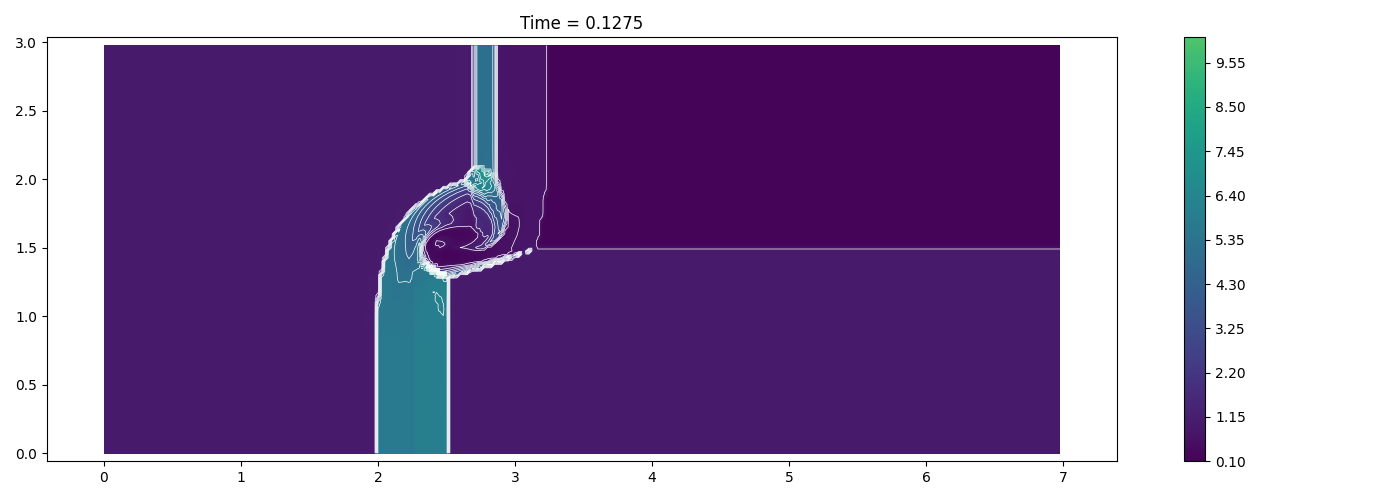}
    \end{minipage}
    \begin{minipage}{0.465\textwidth}
    \centering
        Time = 0.1275\\
        \includegraphics[width=\textwidth,trim={0 0 190 25},clip]{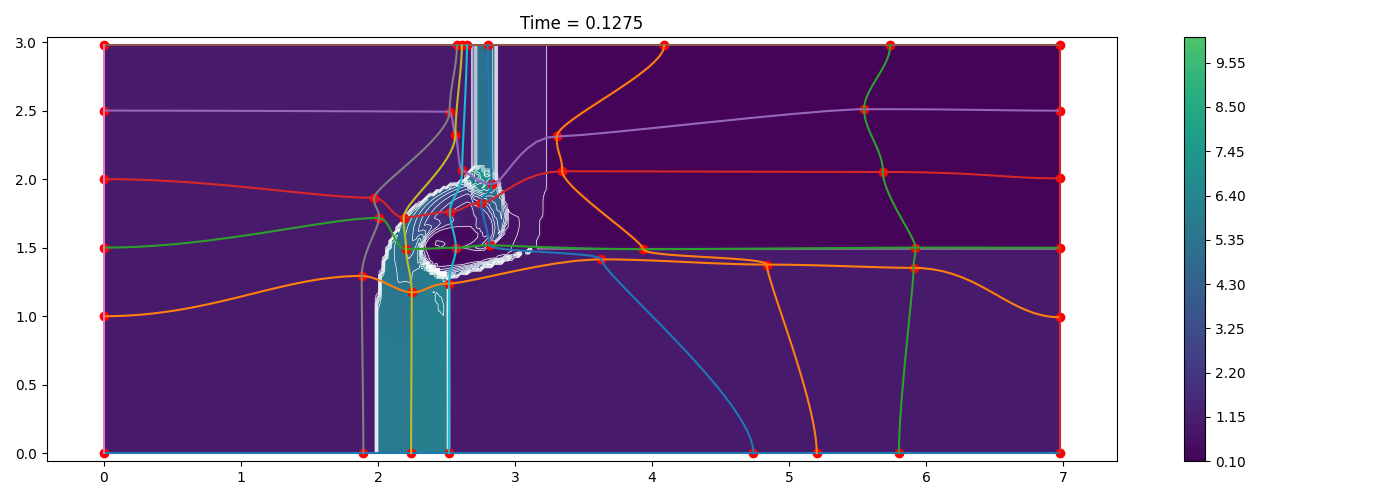}\\  
        \includegraphics[width=\textwidth,trim={0 0 190 25},clip]{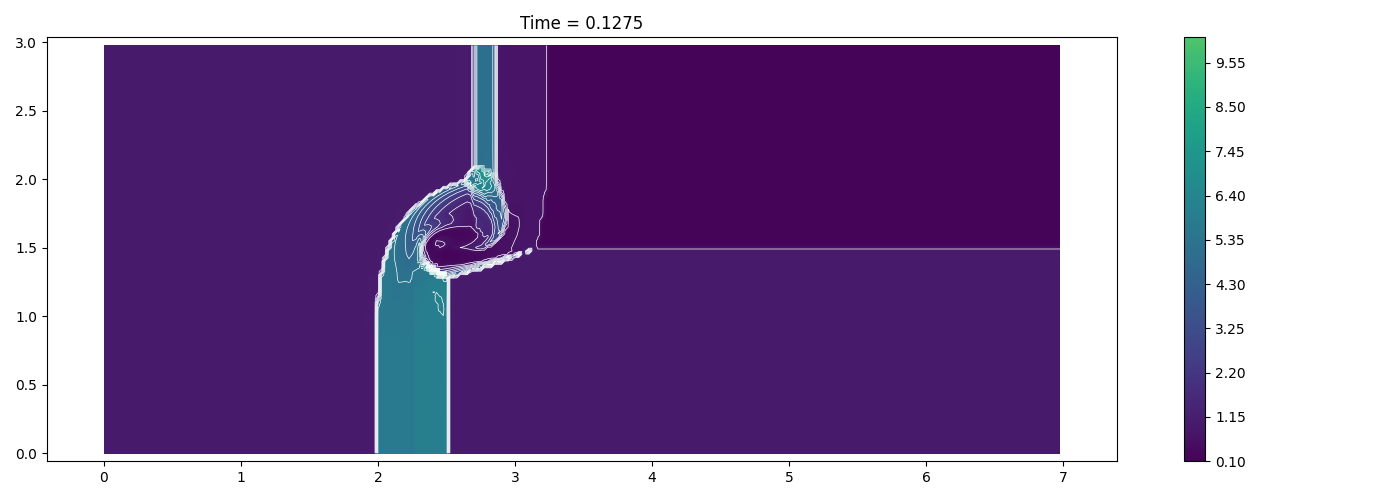}\\        
        \includegraphics[width=\textwidth,trim={0 0 190 25},clip]{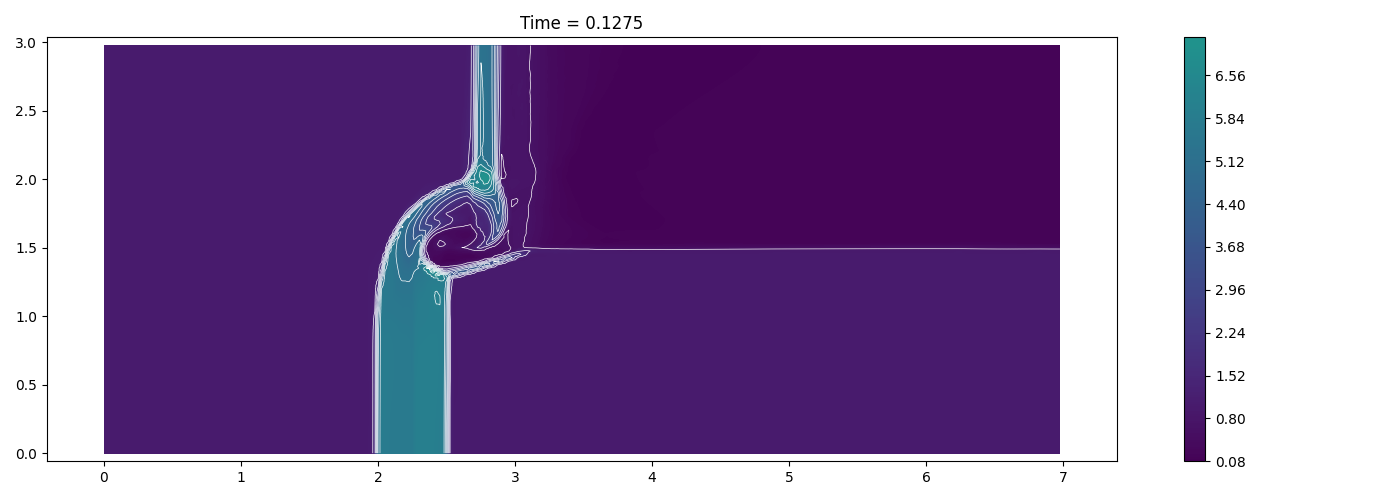}\\        
        \includegraphics[width=\textwidth,trim={0 0 190 25},clip]{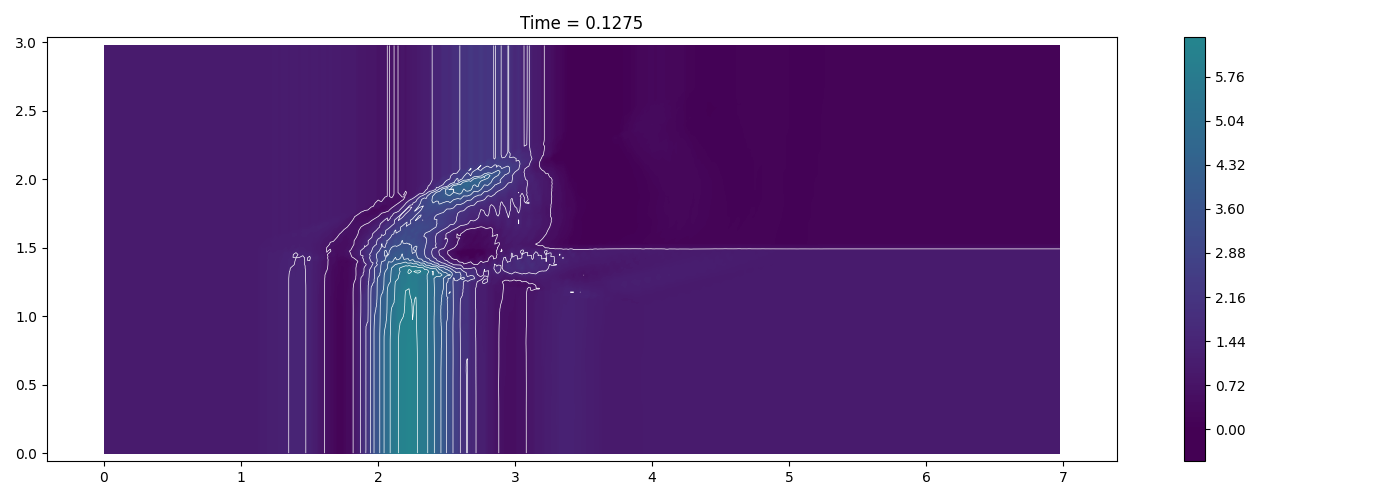}
    \end{minipage}
    \begin{minipage}{0.465\textwidth}
    \centering
        Time = 0.25\\
        \includegraphics[width=\textwidth,trim={0 0 190 25},clip]{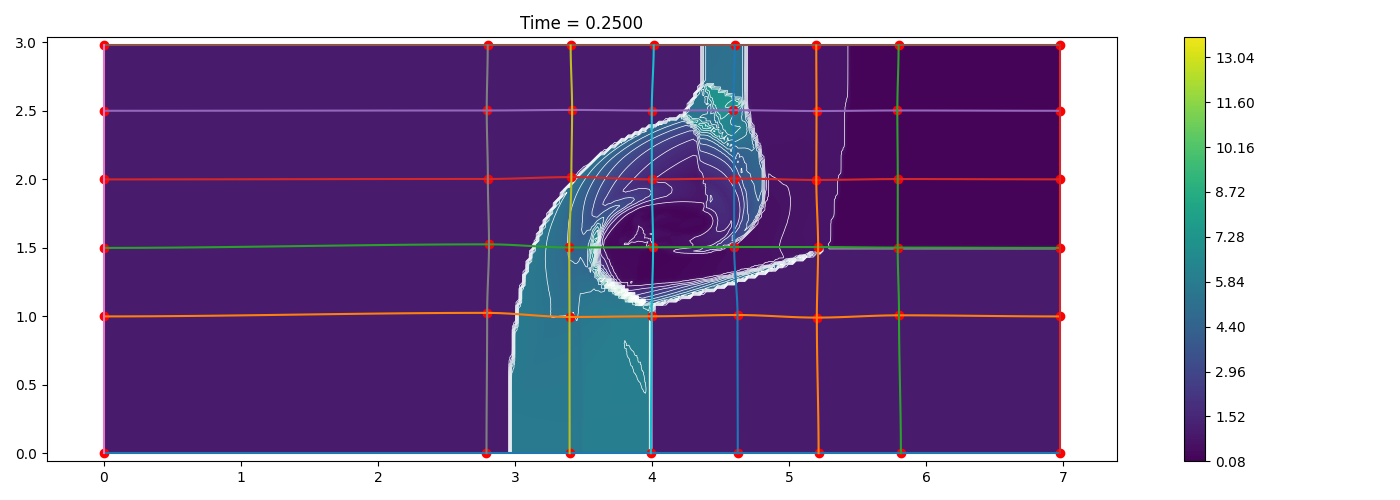}\\
        \includegraphics[width=\textwidth,trim={0 0 190 25},clip]{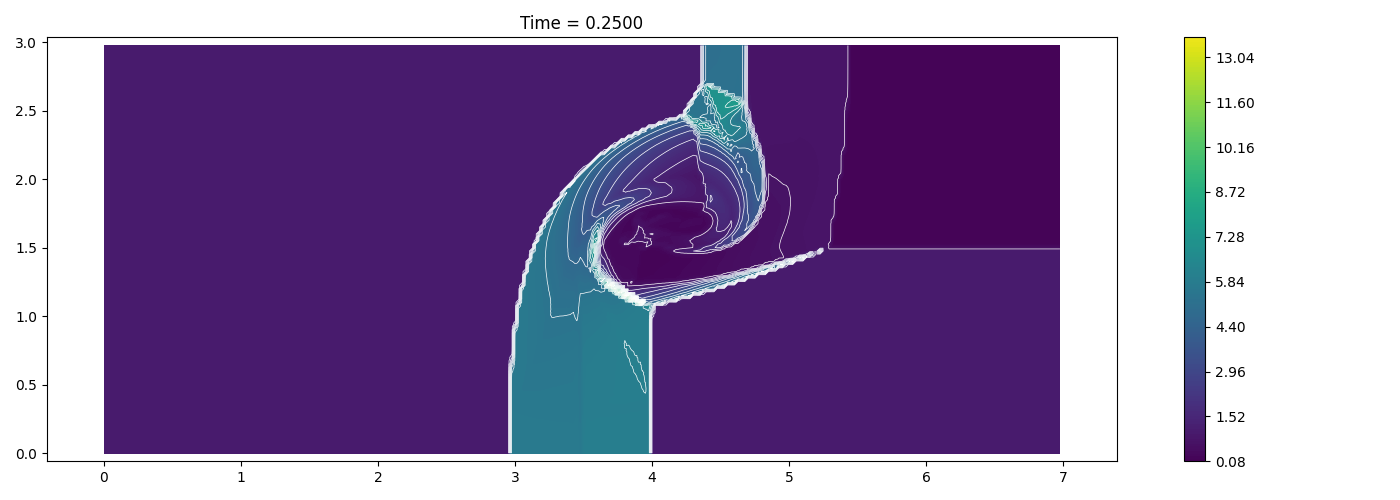}\\        
        \includegraphics[width=\textwidth,trim={0 0 190 25},clip]{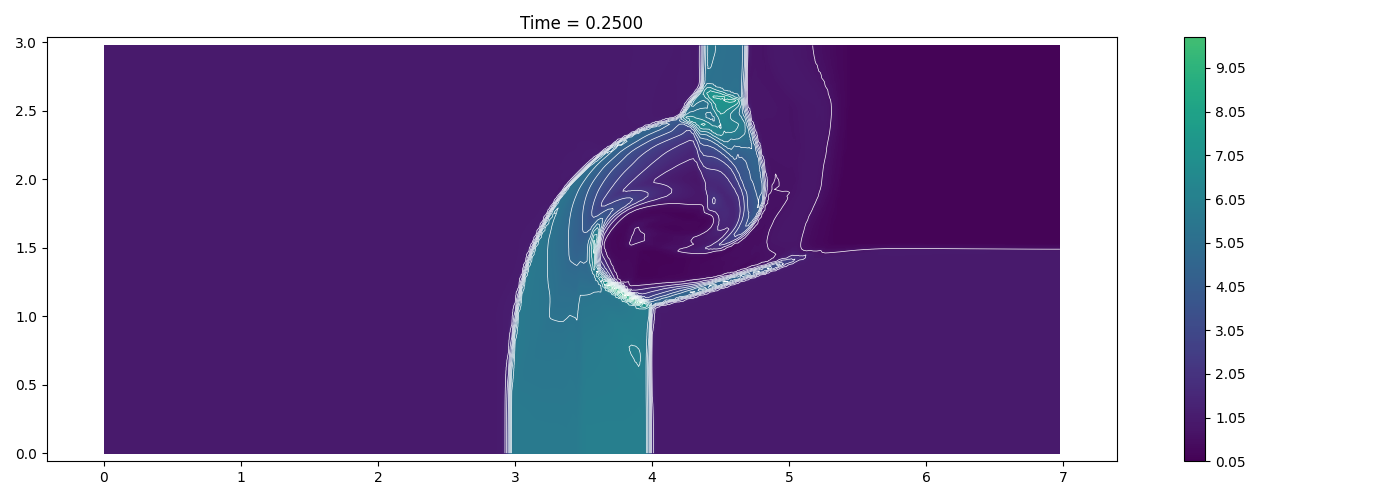}\\        
        \includegraphics[width=\textwidth,trim={0 0 190 25},clip]{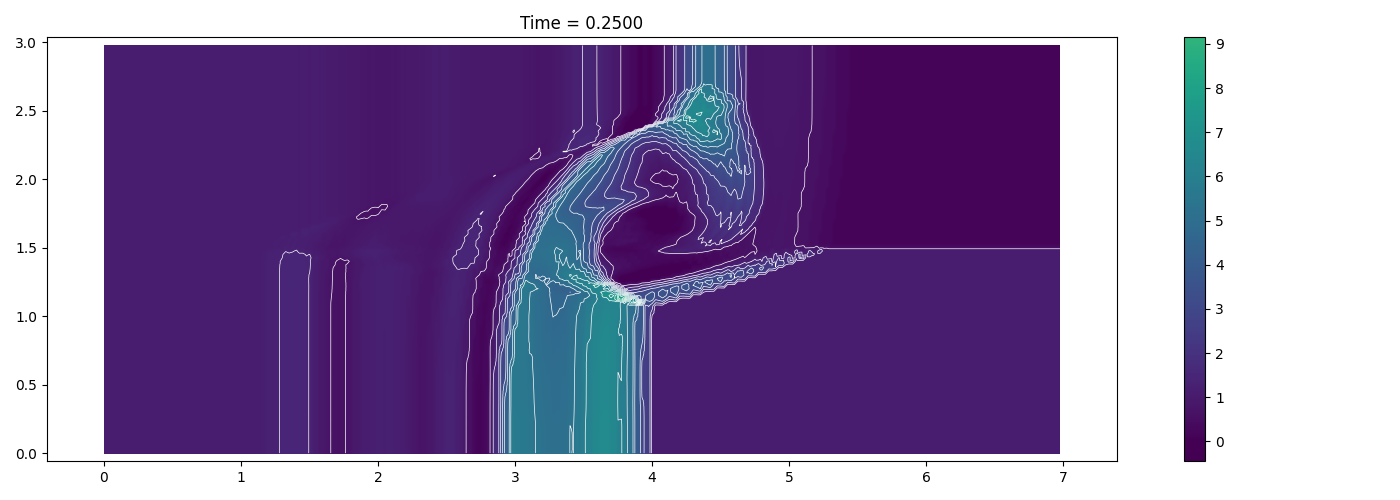}
    \end{minipage}
    \caption{FOM (first and second line), ALE POD-NN (third line) and Eulerian POD-NN (bottom) solutions with $N=7$ of $\rho$ at times $t=0.1275$ (left) and $t=0.25$ (right) in the physical domain $\Omega$. The first line shows also the calibration points. {We plot in white 20 contour lines at equispaced values between 1 and 15}}
    \label{fig:triple_ROM_solutions}
\end{figure}
In Figure~\ref{fig:triple_ROM_solutions}, we show the solutions of the calibration points, FOM, ALE ROM and Eulerian ROM at times $t=0.1275$ and $t=0.25$. The ROM solutions are obtained with a reduced basis of dimension $N=7$.
In the first FOM plots, we also plot the $M_1\cdot M_2=8\cdot 6$ optimal calibration points (the one at the final time are essentially the reference one) that vaguely surround the most dynamical area.
Again, the ALE ROM performs much better than the Eulerian ROM in catching the right position of the waves and to sharply represent them. 
On the other side, there is a slight mistake in the calibration in catching the most right shock that is represented by a vertical contour line, which is not perfectly represented in the ALE ROM.

\begin{figure}
    \centering
    \newcommand{\whicherror}{L2_errors}
    \begin{tikzpicture}
        \begin{axis}[ymode=log,
                    xlabel={Time},
				ylabel={Error},
				width=.68\textwidth,
				height=.45\textwidth,
                    legend pos=outer north east,
                    grid=both,
                    xmin = -0.01,
                    xmax= 0.26,
                    ymax=1,
                    style={font=\footnotesize}]
            \addplot[mark=square,mark size=1.3pt,red] table [x=Time, y=L2Error_eul_FOM_phy, col sep=comma] {Dat065.csv};
            \addlegendentry{Eul POD $N=2$}
            \addplot[mark=triangle*,mark size=1.3pt,red] table [x=Time, y=L2Error_eul_FOM_phy, col sep=comma] {Dat066.csv};
            \addlegendentry{Eul POD $N=4$}            
            \addplot[mark=diamond,mark size=1.3pt,red] table [x=Time, y=L2Error_eul_FOM_phy, col sep=comma] {Dat067.csv};
            \addlegendentry{Eul POD $N=7$}           
            \addplot[mark=otimes*,mark size=1.3pt,red] table [x=Time, y=L2Error_eul_FOM_phy, col sep=comma] {Dat068.csv};
            \addlegendentry{Eul POD $N=12$}
            \addplot[mark=pentagon,mark size=1.3pt,red] table [x=Time, y=L2Error_eul_FOM_phy, col sep=comma] {Dat069.csv};
            \addlegendentry{Eul POD $N=30$}

            \addplot[mark=square,mark size=1.3pt,blue] table [x=Time, y=L2Error_ALE_FOM_phy, col sep=comma] {Dat070.csv};
            \addlegendentry{ALE POD $N=2$}
            \addplot[mark=triangle*,mark size=1.3pt,blue] table [x=Time, y=L2Error_ALE_FOM_phy, col sep=comma] {Dat071.csv};
            \addlegendentry{ALE POD $N=4$}
            \addplot[mark=diamond,mark size=1.3pt,blue] table [x=Time, y=L2Error_ALE_FOM_phy, col sep=comma] {Dat072.csv};
            \addlegendentry{ALE POD $N=7$}
            \addplot[mark=otimes*,mark size=1.3pt,blue] table [x=Time, y=L2Error_ALE_FOM_phy, col sep=comma] {Dat073.csv};
            \addlegendentry{ALE POD $N=12$}
            \addplot[mark=pentagon,mark size=1.3pt,blue] table [x=Time, y=L2Error_ALE_FOM_phy, col sep=comma] {Dat074.csv};
            \addlegendentry{ALE POD $N=30$}
        \end{axis}
    \end{tikzpicture}    
    \vspace{-1mm}
    \caption{Triple point non parametric: Error in time of reduced methods with different number $N$ of modes}
    \label{fig:error_triple}
\end{figure}
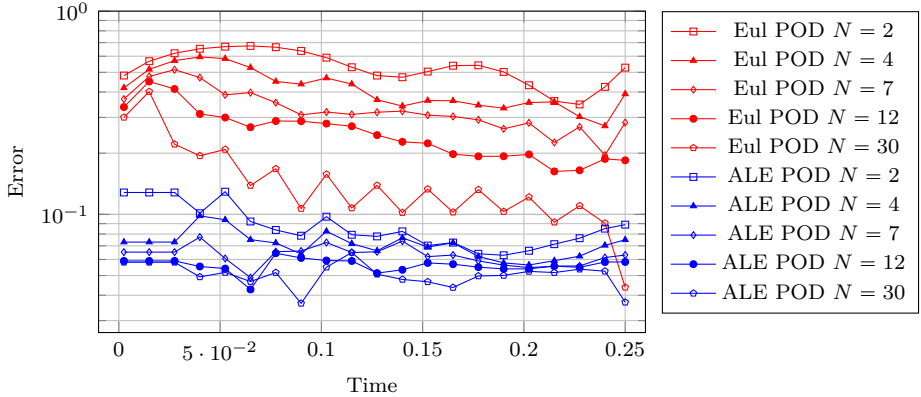
In Figure~\ref{fig:error_triple}, we observe that also here the error of the ALE approach is much smaller than the Eulerian, but that the interpolation error is probably a lower bound for the ALE ROM that keeps it from decreasing with the increasing of the reduced basis dimension. 
On the other side, for the Eulerian ROM the error decays increasing the reduced basis dimensions, even if with irregularity due to the spurious oscillations.
All the parameters not specified are as in the DMR non-parametric case.}

\section{Conclusions}\label{sec:conclusion}
We presented a novel, optimization-based calibration technique suited for hyperbolic conservation laws with (quasi) self-similar solutions that present multiple travelling structures, such as discontinuities. 
We combined the calibration technique with an ANN based Model Order Reduction, in order to obtain a non intrusive approximation setting that is able to provide satisfying results both in the non parametric and in the parametric framework, without the use of implicit shock tracking techniques, which additionally translates into a much more limited computational effort during the offline phase. 
To test the proposed methodology and to show its broad range of applicability, we considered three time-dependent problems of interest: the 1D Sod shock tube problem (non parametric and parametric), the 2D DMR problem (non parametric and parametric) and the non parametric triple point problem. 
In all of our tests, we have shown the benefits of using the proposed calibration based MOR: this is confirmed not only by the comparison on the rate of decay of the eigenvalues returned by the PODs, but also by the behavior in time of the relative $L^2$-errors (in the physical domain $\Omega$) obtained with the two approaches. To conclude, a qualitative comparison on the FOM solutions and the ROM solutions (with and without the calibration approach) is provided, in order to highlight that, by using a smaller number of modes, our proposed methodology is able to correctly capture all the important features of the full order solutions.
Indeed, classical ROMs produce oscillations, smear the shocks and cannot preserve flat areas, while the presented calibrated version does, even in the context of multiple intersecting shocks and waves (such as in the triple point test).\\
{We also showed the robustness of the calibration algorithm with respect to the choice of the reference solution $\rhoref$, the initial guess $\bbtheta^{(0)}(\bbmu)$ and the order with which we perform the calibration: all the tests have been performed for the $1D$ problem, for which an analytical solution is available. The results show that the calibration algorithm provides good results, almost independently on the choice of the reference control points. Nevertheless, we are aware that a more in depth study has to be carried out on the number of control points to choose: we envision this as a future development of the proposed work.}
The replacement of the Neural Networks with a purely ALE approach for the online system is a work in progress and a future extension of this present work. {At the time being, the strongest limitations of our method to get more physically complicated solutions are two: the fact that the domain $\Omega$ needs to be rectangular, and the fact that the configuration of the features should not vary too much, in a topological sense.
We intend to address both points in the future, introducing techniques that can map non-polyhedral shapes into rectangles \cite{taddei2023compositional}, and using local ROMs for different parameter/time zones.}\\
The proposed approximation setting is based on the use of piecewise cubic Hermite interpolating polynomials (or on some tensorial product of them), and works well with rectangular domains and Cartesian meshes: the extension of this approach to more complex geometries and other kinds of meshes (i.e. triangular ones) is envisioned as another future direction of this work. 
We also remark that, so far, we only worked with FV approximations of the full order solution. We expect to generalize the whole methodology to other discretization techniques.

\subsubsection*{Acknowledgments}
D.T. is member of the Gruppo Nazionale Calcolo Scientifico~-~Istituto Nazionale di Alta Matematica (GNCS-INdAM). Part of the research has been carried out while D.T. was a scientific guest of prof. Ilaria Perugia at the Faculty of Mathematics, University of Vienna, and while M.N. was a scientific guest of prof. Gianluigi Rozza at SISSA mathLab. The authors thank the two institutions for the hospitality. D.T. was funded by the Ateneo Sapienza projects 2022 ``Approssimazione numerica di modelli differenziali e applicazioni'' and 2023 ``Modeling, numerical treatment of hyperbolic equations and optimal control problems''.

\bibliographystyle{siamplain}
\bibliography{references}

\nomenclature[01]{\(\Omega\)}{the physical domain}
\nomenclature[02]{\(\mathcal{R}\)}{the reference domain}
\nomenclature[03]{\(d\)}{the dimension of the physical domain}
\nomenclature[04]{\(t\)}{the time variable}
\nomenclature[05]{\(\mu\)}{the physical parameter(s)}
\nomenclature[06]{\(\bbmu\)}{$\bbmu=(t,\mu)$ (parametric case) or $\bbmu=t$ (non-parametric case)}
\nomenclature[07]{\(\mathcal{P}_{\small \text{phys}}\)}{the domain of $\mu$ the physical parameter}
\nomenclature[08]{\(\mathcal{P}\)}{the domain of $\bbmu$ the time and physical parameters}
\nomenclature[09]{$M$}{the number of control points}
\nomenclature[10]{\(\bbw(\bbmu)\)}{the vector of $M$ control points in $\Omega$}
\nomenclature[11]{\(\bbwref(\bbmu)\)}{the vector $M$ reference control points in $\Omegaref$}
\nomenclature[12]{\(\bbalpha\)}{the multiindex $\bbalpha=(\alpha_1,\dots,\alpha_d)$}
\nomenclature[13]{\(\bbw^k_{\bbalpha}\)}{the $k$-th coordinate of the control point $\bbw(\bbmu)_{\bbalpha}$}
\nomenclature[14]{\(\bbtheta(\bbmu)\)}{the vector of the free coordinates of the control point $\bbw(\bbmu)$}
\nomenclature[15]{\(\thetavec\)}{the matrix of the free coordinates. $\thetavec\in\mathbb{R}^{N_{few}}\times{Q}$, and $\thetavec[i,:]=\bbtheta(\tildebbmu_i)$, $i=1,\dots,N_{few}$.}
\nomenclature[16]{\(T[\cdot]\)}{the family of geometrical transformation maps. $T[\cdot]\in \mathcal{C}^1(\Omega^M, \mathcal{C}^1(\Omegaref, \Omega))$.}
\nomenclature[17]{\(\rhoref(\cdot)\)}{the reference solution for the calibration in the self-similar setting}
\nomenclature[18]{\(\rhohat(\cdot)\)}{the calibrated snapshot}
\nomenclature[19]{\(N^{\text{POD}}_{\text{few}}\)}{the dimension of the linear space $V_{\text{POD}}$ used for the calibration in the quasi-self-similar setting}
\nomenclature[20]{\(V_{\text{POD}}\)}{the linear space considered in the residual function to perform the projection error for the minimization in the quasi-self-similar case, instead of the use of the reference solution}
\nomenclature[21]{\(\mathcal{M}_{\rho}\)}{the original solution manifold (for the density $\rho$)}
\nomenclature[22]{\(\mathcal{\hat{M}}_{\rho}\)}{the calibrated solution manifold (for  the density $\rho$)}
\nomenclature[23]{\small \(\{\phiref_i\}_{i=1}^N\)}{the set of $N$ reduced basis functions obtained by POD-compression of $\mathcal{\hat{M}}_{\rho}$}
\nomenclature[24]{\(\underline{\rhohat}_N(\bbmu)\)}{the vector of the $L^2$-projection coefficients of the calibrated snapshot $\rhohat(\bbmu)$ onto the linear space spanned by $\{\phiref_i\}_{i=1}^N$}
\printnomenclature
\end{document}